\documentstyle[11pt,amssymb]{article}
\begin{document}

{\Large

\noindent{\bf Duality of $q$-polynomials, orthogonal}

\noindent{\bf on countable sets of points} }
\bigskip

\noindent{\sc N. M. Atakishiyev${}^1$ and A. U. Klimyk${}^{2}$}
\bigskip

\noindent ${}^1$Instituto de Matem\'aticas, UNAM, CP 62210
Cuernavaca, Morelos, M\'exico

\noindent ${}^2$Bogolyubov Institute for Theoretical Physics,
03143 Kiev, Ukraine

\medskip

E-mail: natig@matcuer.unam.mx and aklimyk@bitp.kiev.ua
\bigskip

\begin{abstract}
We review properties of $q$-orthogonal polynomials, related to
their orthogonality, duality and connection with the theory of
symmetric (self-adjoint) operators, represented by a Jacobi
matrix. In particular, we show how one can naturally interpret the
duality of families of polynomials, orthogonal on countable sets
of points. In order to obtain orthogonality relations for dual
sets of polynomials, we propose to use two symmetric
(self-adjoint) operators, representable (in some distinct bases)
by Jacobi matrices. To illustrate applications of this approach,
we apply it to several pairs of dual families of $q$-polynomials,
orthogonal on countable sets, from the $q$-Askey scheme. For each
such pair, the corresponding operators, representable by Jacobi
matrices, are explicitly given. These operators are employed in
order to find explicitly sets of points, on which the polynomials
are orthogonal, and orthogonality relations for them.

\end{abstract}

\bigskip

\medskip

\bigskip
\bigskip

\begin{flushright}
\begin{minipage}{5in}
{\it "We mathematicians are particularly fond of duality theorems;
translating mathematical statements from one category to another
often gives us new and unexpected insight", M.Harris, "Postmodern
at an Early Age", Notices of the American Mathematical Society,
Vol.50, No.7, p.792, 2003.}

\end{minipage}
\end{flushright}

\bigskip

\noindent{\bf 1. INTRODUCTION}
\bigskip

It is well known that each family $p_n(x)$, $n=0,1,2,\cdots$, of
orthogonal polynomials of one variable corresponds to the
determinate or indeterminate moment problem. If a polynomial
family corresponds to the determinate moment problem, then there
exists only one positive orthogonality measure $\mu$ for these
polynomials and they constitute a complete orthogonal set in the
Hilbert space $L^2(\mu)$. If a family corresponds to the
indeterminate moment problem, then there exists infinitely many
orthogonality measures $\mu$ for these polynomials and these
measures are divided into two parts: extremal measures and
non-extremal measures. If a measure $\mu$ is extremal, then the
corresponding set of polynomials constitute a complete orthogonal
set in the Hilbert space $L^2(\mu)$. If a measure $\mu$ is not
extremal, then the corresponding family of polynomials is not
complete in the Hilbert space $L^2(\mu)$ (see [ST]).

It is also well known that there exists a close relation of the
theory of orthogonal polynomials with the theory of symmetric
(self-adjoint) operators, representable by a Jacobi matrix. The
point is that with each family of orthogonal polynomials one can
associate a closed symmetric (or self-adjoint) operator $A$,
representable by a Jacobi matrix. If the corresponding moment
problem is indeterminate, then the operator $A$ is not
self-adjoint and it has infinitely many self-adjoint extensions.
If the operator $A$ has a physical meaning, then these
self-adjoint extensions are especially important. These extensions
correspond to extremal orthogonality measures for the same set of
polynomials and can be constructed by means of these measures
(see, for example, [Ber], Chapter VII, and [Sim]). If the family
of polynomials corresponds to the determinate moment problem, then
the corresponding operator $A$ is self-adjoint and its spectrum is
determined by an orthogonality relation for the polynomials.
Moreover, the spectral measure for the operator $A$ is constructed
by means of the orthogonality measure for the corresponding
polynomials (see [Ber], Chapter VII).

In section 2 we briefly review the relations between the theory of
orthogonal polynomials, the theory of operators, representable by
a Jacobi matrix, and the theory of moment problem. This
information is basic for the exposition in further sections. In
section 2 we also discuss how one can naturally extend the
conventional notion of duality to sets of polynomials, orthogonal
on countable sets of points.

In order to find orthogonality measures for dual sets of
polynomials, we use two symmetric (or self-adjoint) operators,
representable (with respect to different bases) by Jacobi
matrices. This approach is applied to several  sets of dual
$q$-orthogonal polynomials from the Askey scheme.

Pairs of operators $(A,I)$, employed for studying some sets of
$q$-orthogonal polynomials and their duals, belong to the discrete
series representations of the quantum algebra $U_q({\rm
su}_{1,1})$ (see, for example, [AK1] and [AK6]). However, in order
to facilitate ease of comprehending to a larger number of readers
we have not exploited this deep algebraic fact; that is, we
exhibit explicit forms of these operators without using the
representation theory of the quantum algebra $U_q({\rm
su}_{1,1})$. These pairs of operators are, in fact, a
generalization of Leonard pairs, introduced by P. Terwilliger
[Ter1] (for the definition and references see section 3).

When one considers dual sets of $q$-polynomials, orthogonal
on countable sets of points, then one member of these sets
corresponds to the determinate moment problem and another to
the indeterminate moment problem. One of the two operators
$(A,I)$ (that is, the operator $A$) for a given dual pair of
sets of $q$-orthogonal polynomials corresponds to a three-term
recurrence relation for the set of polynomials, which corresponds
to the determinate moment problem. This operator is bounded and
self-adjoint; moreover, it has the discrete spectrum. We diagonalize
this self-adjoint bounded operator and find its spectrum with the
aid of the second operator $I$, which corresponds to a $q$-difference
equation for the same set of polynomials. An explicit form of all
eigenfunctions for the operator $A$ is found for each  dual set
of polynomials, considered by us. They are expressed in terms of
$q$-polynomials, which belong to the set, associated with the
determinate moment problem. Since the spectrum of $A$ is simple,
its eigenfunctions form an orthogonal basis in the Hilbert space.
One can normalize this basis. This normalization is effected by
means of the second operator $I$ from the corresponding pair. As
a result of this normalization, two orthonormal bases in the Hilbert
space emerge: the canonical (or the initial) basis and the basis of
eigenfunctions of the operator $A$. They are interrelated by a
unitary matrix $U$, whose entries $u_{mn}$ are explicitly expressed
in terms of polynomials $P_m(x)$, corresponding to the determinate
moment problem. Since the matrix $U$ is unitary (and in fact it is
real in our case), there are two orthogonality relations for its
elements, namely
$$
\sum _{n} u_{mn}u_{m'n}=\delta_{mm'},\ \ \ \sum _{m}
u_{mn}u_{mn'}=\delta_{nn'}.   \eqno (1.1)
$$
The first relation expresses the orthogonality relation for the
polynomials $P_m(x)$, which correspond to the determinate moment
problem. So, the orthogonality of $U$ yields an algebraic proof of
orthogonality relation for these polynomials. In order to
interpret the second relation, we consider the polynomials
$P_m(x_n)$ ($\{ x_n\}$ is the set of points, on which the
polynomials are orthogonal) as functions of $m$. In this way
one obtains one or two sets of orthogonal functions, which are
expressed in terms of a dual set of $q$-orthogonal polynomials
(which corresponds to the indeterminate moment problem). The
second relation in (1.1) leads to the orthogonality relations
for these dual $q$-orthogonal polynomials.

Since this set of $q$-orthogonal polynomials corresponds to
the indeterminate moment problem, there are infinitely many
orthogonality relations. Using the pair of operators $(A,I)$
and the notion of duality, one is able to find only one
orthogonality relation (which is dual to the orthogonality
relation for the corresponding set of polynomials, associated
with the determinate moment problem). Sometimes a measure,
which corresponds to this orthogonality relation, is extremal
and sometimes it is not extremal. It depends on a concrete pair
of dual sets of polynomials.

Throughout the sequel we always assume that $q$ is a fixed
positive number such that $q<1$. We use (without additional
explanation) notations of the theory of special functions
and the standard $q$-analysis (see, for example, [GR] and
[AAR]). We shall also use the well-known shorthand notation
$(a_1,\cdots,a_k;q)_n:=(a_1;q)_n\cdots (a_k;q)_n$.
\bigskip

\noindent{\bf 2. ORTHOGONALITY MEASURES AND DUALITY}
 \bigskip

\noindent{\bf 2.1. Orthogonal polynomials, Jacobi matrices and the
moment problem}
 \bigskip

Orthogonal polynomials are closely related to operators
represented by a Jacobi matrices. In what follows we shall use
only symmetric Jacobi matrices and the word "symmetric" will be
often omitted. By a symmetric Jacobi matrix we mean a (finite or
infinite) symmetric matrix of the form
$$
M =\left( \matrix{
 b_{0}&  a_{0} &0     &0&0&  \cdots   \cr
 a_{0}&  b_1   &a_{1} &0&0&  \cdots   \cr
 0   &   a_1   &b_{2} &a_2&0&  \cdots  \cr
 0   &   0    &a_{2} &b_3&a_3&  \cdots  \cr
 \cdots&\cdots &\cdots &\cdots &\cdots &\cdots  \cr}
 \right) . \eqno (2.1)
$$
We assume below that all $a_i\ne 0$, $i=0,1,2,\cdots$; then $a_i$
are real. Let $L$ be a closed symmetric operator on a Hilbert
space ${\cal H}$, representable by a Jacobi matrix $M$. Then there
exists an orthonormal basis $e_n$, $n=0,1,2,\cdots$, in ${\cal
H}$, such that
 $$
Le_n=a_ne_{n+1}+b_ne_n+a_{n-1}e_{n-1},
 $$
where $e_{-1}\equiv 0$. Let $|x\rangle =\sum_{n=0}^\infty p_n(x)
e_n$ be an eigenvector\footnote{Observe that eigenvectors of $L$
may belong to either the Hilbert space ${\cal H}$ or to some
extension of ${\cal H}$. (For example, if ${\cal H}=L^2(-\infty,
\infty)$ and in place of $L$ we have the operator $d/dx$, then the
functions $e^{{\rm i}xp}$, which do not belong to $L^2(-\infty,
\infty)$, are eigenfunctions of $d/dx$.) Below we act freely with
eigenvectors, which do not belong to ${\cal H}$, but this can be
easily made mathematically strict.} of $L$ with an eigenvalue $x$,
that is, $L| x\rangle =x|x\rangle$. Then
 $$
L|x\rangle =\sum_{n=0}^\infty [  p_n(x) a_n e_{n+1}+p_n(x) b_n
e_{n}+p_n(x) a_{n-1} e_{n-1}]=x \sum_{n=0}^\infty p_n(x) e_n.
 $$
Equating coefficients of the vector $e_n$, one comes to a
recurrence relation for the coefficients $p_n(x)$:
 $$
a_np_{n+1}(x)+b_np_{n}(x)+a_{n-1}p_{n-1}(x)=xp_n(x).  \eqno (2.2)
 $$
Since $p_{-1}(x)=0$, by setting $p_0(x)\equiv 1$ we see that
$p_1(x)=a_0x-b_0/a_0$. Similarly we can find uniquely $p_n(x)$,
$n=2,3,\cdots$. Thus, the relation (2.2) completely determines the
coefficients $p_n(x)$. Moreover, the recursive computation of
$p_n(x)$ shows that these coefficients $p_n(x)$ are polynomials in
$x$ of degree $n$. Since the coefficients $a_n$ and $b_n$ are real
(because the matrix $M$ is symmetric),  all coefficients of the
polynomials $p_n(x)$ themselves are real.

Well-known Favard's characterization theorem for polynomials
$P_n(x)$, $n=0,1,2,\cdots$, of degree $n$ states that if these
polynomials satisfy a recurrence relation
 $$
A_nP_{n+1}(x)+B_nP_{n}(x)+C_{n}P_{n-1}(x)=xP_n(x)
 $$
and the conditions $A_nC_{n+1}>0$ are satisfied, then these
polynomials are orthogonal with respect to some positive measure.
It is clear that the conditions of Favard's theorem are satisfied
for the polynomials $p_n(x)$ because in this case the requirements
simply reduce to inequalities $a_n^2>0$ for $n=0,1,2,\cdots$. This
means that the polynomials $p_n(x)$ from (2.2) are orthogonal with
respect to some positive measure $\mu(x)$. It is known that
orthogonal polynomials admit orthogonality with respect to either
unique positive measure or with respect to infinitely many
positive measures.

The polynomials $p_n(x)$ are very important for studying
properties of the closed symmetric operator $L$. Namely, the
following statements are true (see, for example, [Ber] and [Sim]):
\medskip

I. Let the polynomials $p_n(x)$ are orthogonal with respect to a
unique orthogonality measure $\mu$,
 $$
\int p_m(x)p_n(x)d\mu(x)=\delta_{mn},
 $$
where the integration is performed over some subset (possibly
discrete) of ${\Bbb R}$, then the closed operator $L$ is
self-adjoint. Moreover, the spectrum of the operator $L$ is simple
and coincides with the set, on which the polynomials $p_n(x)$ are
orthogonal (recall that we assume that all numbers $a_n$ are
non-vanishing). The measure $\mu(x)$ determines the spectral
measure for the operator $L$ (for details see [Ber], Chapter VII).
 \medskip

II. Let the polynomials $p_n(x)$ are orthogonal with respect to
infinitely many different orthogonality measures $\mu$. Then the
closed symmetric operator $L$ is not self-adjoint and has
deficiency indices (1, 1), that is, it has infinitely many (in
fact, one-parameter family of) self-adjoint extensions. It is
known that among orthogonality measures, with respect to which the
polynomials are orthogonal, there are so-called extremal measures
(that is, such measures that a set of polynomials $\{ p_n(x)\}$ is
complete in the Hilbert space $L^2$ with respect to the
corresponding measure; see subsection 2.3 below). These measures
uniquely determine self-adjoint extensions of the symmetric
operator $L$. There exists one-to-one correspondence between
essentially distinct extremal orthogonality measures and
self-adjoint extensions of the operator $L$. The extremal
orthogonality measures determine spectra of the corresponding
self-adjoint extensions.
 \medskip

The inverse statements are also true:
 \medskip

I$'$. Let the operator $L$ be self-adjoint. Then the corresponding
polynomials $p_n(x)$ are orthogonal with respect to a unique
orthogonality measure $\mu$,
 $$
\int p_m(x)p_n(x)d\mu(x)=\delta_{mn},
 $$
where the integral is taken over some subset (possibly discrete)
of ${\Bbb R}$, which coincides with the spectrum of $L$. Moreover,
a measure $\mu$ is uniquely determined by a spectral measure for
the operator $L$ (for details see [Ber], Chapter VII).
 \medskip

II$'$. Let the closed symmetric operator $L$ be not self-adjoint.
Since it is representable by a Jacobi matrix (2.1) with $a_n\ne
0$, $n=0,1,2,\cdots$, it admits one-parameter family of
self-adjoint extensions (see [Ber], Chapter VII). Then the
polynomials $p_n(x)$ are orthogonal with respect to infinitely
many orthogonality measures $\mu$. Moreover, spectral measures of
self-adjoint extensions of $L$ determine extremal orthogonality
measures for the polynomials $\{ p_n(x)\}$ (and a set of
polynomials $\{ p_n(x)\}$ is complete in the Hilbert spaces
$L^2(\mu)$ with respect to the corresponding extremal measures
$\mu$).
 \medskip

On the other hand, with the orthogonal polynomials $p_n(x)$,
$n=0,1,2,\cdots$, the classical moment problem is associated (see
[ST] and [Akh]). Namely, with these polynomials (that is, with the
coefficients $a_n$ and $b_n$ in the corresponding recurrence
relation) real numbers $c_n$, $n=0,1,2,\cdots$, are associated,
which determine the corresponding classical moment problem. (The
numbers $c_n$ are uniquely determined by $a_n$ and $b_n$.) The
definition of the classical moment problem consists in the
following. Let a set of real numbers $c_n$, $n=0,1,2,\cdots$, be
given. We are looking for a positive measure $\mu(x)$, such that
 $$
\int x^nd\mu(x)=c_n,\ \ \ \ n=0,1,2,\cdots ,  \eqno (2.3)
 $$
where the integration is taken over ${\Bbb R}$. (In this case
we deal with the {\it Hamburger moment problem}.) There are two
principal questions in the theory of moment problem:
 \medskip

(i) Does there exist a measure $\mu(x)$, such that relations (2.3)
are satisfied?

(ii) If such a measure exists, is it determined uniquely?
 \medskip

The answer to the first question is positive, if the numbers
$c_n$, $n=0,1,2,\cdots$, are those, which correspond to a family
of orthogonal polynomials. Moreover, a measure $\mu(x)$ then
coincides with the measure, with respect to which these
polynomials are orthogonal.

If a measure $\mu$ in (2.3) is determined uniquely, then we say
that we deal with {\it the determinate moment problem}. In
particular, it is the case when the measure $\mu$ is supported on
a bounded set. If a measure, with respect to which relations (2.3)
hold, is not unique, then we say that we deal with {\it the
indeterminate moment problem}. In this case there exist infinitely
many measures $\mu(x)$ for which (2.3) take place. Then the
corresponding polynomials are orthogonal with respect to all these
measures and the corresponding symmetric operator $L$ is not
self-adjoint. In this case the set of solutions of the moment
problem for the numbers $\{ c_n\}$ coincides with the set of
orthogonality measures for the corresponding polynomials $\{
p_n(x)\}$.

Observe that not each set of real numbers $c_n$, $n=0,1,2,\cdots$,
is associated with a set of orthogonal polynomials. In other
words, there are sets of real numbers $c_n$, $n=0,1,2,\cdots$,
such that the corresponding moment problem does not have a
solution, that is, there is no positive measure $\mu$, for
which the relations (2.3) are true. But if for some set of real
numbers $c_n$, $n=0,1,2,\cdots$, the moment problem (2.3) has a
solution $\mu$, then this set corresponds to some set of
polynomials $p_n(x)$, $n=0,1,2,\cdots$, which are orthogonal with
respect to this measure $\mu$. There exist criteria indicating
when for a given set of real numbers $c_n$, $n=0,1,2,\cdots$, the
moment problem (2.3) has a solution (see, for example, [ST]).
Moreover, there exist procedures, which associate a collection of
orthogonal polynomials to a set of real numbers $c_n$, $n=0,1,2,\cdots$,
for which the moment problem (2.3) has a solution (see, [ST]).
\medskip

Thus, we see that the following three theories are closely
related:
\medskip

(i) the theory of symmetric operators $L$, representable by a
Jacobi matrix;

(ii) the theory of orthogonal polynomials in one variable;

(iii) the theory of classical moment problem.

 \bigskip

\noindent{\bf 2.2. Extremal orthogonality measures}
 \bigskip

To a set of orthogonal polynomials $p_n(x)$, $n=0,1,2,\cdots$,
associated with an indeterminate moment problem (2.3), there
correspond four entire functions $A(z)$, $B(z)$, $C(z)$, $D(z)$,
which are related to appropriate orthogonality measures $\mu$ for
the polynomials by the formula
$$
F(z)\equiv \frac{A(z)-\sigma(z)C(z)}{B(z)-\sigma(z)D(z)}
=\int_{-\infty}^\infty \frac{d\mu(t)}{z-t}    \eqno (2.4)
$$
(see, for example, [ST]), where $\sigma(z)$ is any analytic
function. Moreover, to each analytic function $\sigma(z)$
(including cases of constant $\sigma(z)$ and $\sigma(z)=\pm
\infty$) there corresponds a single orthogonality measure
$\mu(t)\equiv \mu_\sigma (t)$ and, conversely, to each
orthogonality measure $\mu$ there corresponds an analytic
function $\sigma$ such that formula (2.4) holds. There exists
the Stieltjes inversion formula, which converts the formula
(3.1). It has the form
$$
[\mu(t_1+0)+\mu(t_1-0)]-[\mu(t_0+0)+\mu(t_0-0)] \qquad\qquad\qquad
 $$  $$ \qquad\qquad
= \lim_{\varepsilon\to +0}\left( -\frac{1}{\pi{\rm i}}
\int_{t_0}^{t_1} [F(t+{\rm i}\varepsilon)-F(t-{\rm
i}\varepsilon)]dt \right) .   \eqno (2.5)
$$

Thus, orthogonality measures for a given set of polynomials
$p_n(x)$, $n=0,1,2,\cdots$, in principle, can be found. However,
it is very difficult to evaluate the functions $A(z)$, $B(z)$,
$C(z)$, $D(z)$. In [IM] they are evaluated for particular example
of polynomials, namely, for the $q^{-1}$-continuous Hermite
polynomials $h_n(x|q)$. So, as a rule, for the derivation of
orthogonality measures other methods are used.

The measures $\mu_\sigma (t)$, corresponding to constants $\sigma$
(including $\sigma=\pm \infty$), are called {\it extremal
measures} (some authors, following the book [Akh], call these
measures $N$-extremal). All other orthogonality measures are not
extremal.

The importance of extremal measures is explained by Riesz's
theorem. Let us suppose that a set of polynomials $p_n(x)$,
$n=0,1,2,\cdots$, associated with the indeterminate moment
problem, is orthogonal with respect to a positive measure $\mu$
(that is, $\mu$ is a solution of the moment problem (2.3)). Let
$L^2(\mu)$ be the Hilbert space of square integrable functions
with respect to the measure $\mu$. Evidently, the polynomials
$p_n(x)$ belong to the space $L^2(\mu)$. Riesz's theorem states
the following:
  \medskip

{\bf Riesz's theorem.} {\it The set of polynomials $p_n(x)$,
$n=0,1,2,\cdots$, is complete in the Hilbert space $L^2(\mu)$
(that is, they form a basis in this Hilbert space) if and only if
the measure $\mu$ is extremal}.
 \medskip

Note that if a set of polynomials $p_n(x)$, $n=0,1,2,\cdots$,
corresponds to a determinate moment problem and $\mu$ is an
orthogonality measure for them, then this set of polynomials is
also complete in the Hilbert space $L^2(\mu)$.

In particular, Riesz's theorem is often used in order to determine
whether a certain orthogonality measure is extremal or not.
Namely, if we know that a given set of orthogonal polynomials,
corresponding to an indeterminate moment problem, is not complete
in the Hilbert space $L^2(\mu)$, where $\mu$ is an orthogonality
measure, then this measure is not extremal.

Note that for applications in physics and in functional analysis
it is of interest to have extremal orthogonality measures. If an
orthogonality measure $\mu$ is not extremal, then it is important
to find a system of orthogonal functions $\{f_m(x)\}$, which
together with a given set of polynomials constitute a complete set
of orthogonal functions (that is, a basis in the Hilbert space
$L^2(\mu)$). Sometimes, it is possible to find such systems of
functions (see, for example, [CKK]).

Extremal orthogonality measures have many interesting properties
[ST]:
 \medskip

(a) If $\mu_\sigma(x)$ is an extremal measure, associated
(according to formula (2.4)) with a number $\sigma$, then
$\mu_\sigma(x)$ is a step function. Its spectrum (that is, the set
on which the corresponding polynomials $p_n(x)$, $n=0,1,2,\cdots$,
are orthogonal) coincides with the set of zeros of the denominator
$B(z)-\sigma D(z)$ in (2.4). The mass, concentrated at a spectral
point $x_j$ (that is, a jump of $\mu_\sigma(x)$ at the point
$x_j$), is equal to $(\sum_{n=0}^\infty | p_n(x_j)|^2)^{-1}$.
\medskip

(b) Spectra of extremal measures are real and simple. This means
that the corresponding self-adjoint operators, which are
self-adjoint extensions of the operator $L$, have simple spectra,
that is, all spectral points are of multiplicity 1.
\medskip

(c) Spectral points of two different extremal measures
$\mu_\sigma(x)$ and $\mu_{\sigma'}(x)$ are mutually separated.
\medskip

(d) For a given real number $x_0$, always exists a (unique) real
value $\sigma$, such that the measure $\mu_\sigma(x)$ has $x_0$ as
its spectral point. The points of the spectrum of $\mu_\sigma(x)$
are analytic monotonic functions of $\sigma$.
\medskip

It is difficult to find all extremal orthogonality measures for a
given set of orthogonal polynomials (that is, self-adjoint
extensions of a corresponding closed symmetric operator). As far
as we know, at the present time they are known only for one family
of polynomials, which correspond to indeterminate moment problem.
They are the $q^{-1}$-continuous Hermite polynomials $h_n(x|q)$
(see [IM]).

If extremal measures $\mu_\sigma$ are known then by multiplying
$\mu_\sigma$ by a suitable factor (depending on $\sigma$) and
integrating it with respect to $\sigma$, one can obtain infinitely
many continuous orthogonality measures (which are not extremal).
\bigskip

\noindent{\bf 2.3. Dual sets of orthogonal polynomials}
\bigskip

A notion of duality for two families of polynomials, orthogonal
on finite sets of points, is well known. Namely, let $p_n(x)$,
$n=0,1,2,\cdots,N$, be orthogonal polynomials with orthogonality
relation
$$
\sum_{m=0}^N p_n(x_m)p_{n'}(x_m)
w_m=v_n^{-1}\delta_{nn'}\sum_{s=0}^N w_s , \eqno (2.6)
$$
where $n,n'=0,1,2,\cdots,N$, $w_s>0$ is a jump of the
orthogonality measure in the point $x_s$, and
$$
v_n=\prod_{k=0}^n (a_{k-1}/c_k)
$$
($a_k$ and $c_n$ are coefficients in the three-term recurrence
relation $x p_n(x)=a_n p_{n+1}(x)+b_np_n(x)+c_np_{n-1}(x)$ for
the polynomials $p_n(x)$). Then the dual orthogonality relation
is of the form
$$
\sum_{n=0}^N p_n(x_m)p_{n}(x_{m'})
v_n=w_m^{-1}\delta_{mm'}\sum_{s=0}^N w_s , \eqno (2.7)
$$
where $m,m'=0,1,2,\cdots,N$ (see, for example, [GR], Chapter 7).
If one considers the $p_n(x_m)$ as functions of $n$, in many cases
these functions turn out to be either polynomials in $n$ or
polynomials in $\nu(n)$, where $\nu(n)$ is some function of $n$.
Then the polynomials $P_m(n):=p_n(x_m)$, $m=0,1,2,\cdots,N$
(respectively, $P_m(\nu(n)):=p_n(x_m)$, $m=0,1,2,\cdots,N$) are
orthogonal polynomials of $n$ (respectively of $\nu(n)$), for
which (2.7) is an orthogonality relation. The polynomials
$P_m(\nu(n))$, $m=0,1,2,\cdots,N$, are called {\it dual}
polynomials with respect to the $p_n(x_m)$, $n=0,1,2,\cdots,N$.
If the dual polynomials $\{ P_m(\nu(n))\}$ coincide with the
$\{p_n(m)\}$, then the polynomials $\{ p_n(m)\}$ are called
{\it self-dual}. For instance, Racah polynomials and $q$-Racah
polynomials both represent families of self-dual polynomials.

It is not obvious how to extend the notion of duality to
polynomials, orthogonal on countable sets of points. In the case
of polynomials, orthogonal on a finite set of points, the
orthogonality (2.7) readily follows from the orthogonality (2.6).
Namely, the orthogonality (2.6) means that the real $(N+1)\times
(N+1)$ matrix $(a_{mn})_{m,n=0}^N$ with matrix elements
$$
a_{mn}=c_{mn} p_n(x_m),\ \ \ \ c_{mn}=\left[ w_mv_n/\sum_{s=0}^N
w_s \right]^{1/2} ,
$$
is orthogonal. Orthogonality of its columns is equivalent to the
relation (2.6). Orthogonality by rows for the matrix
$(a_{mn})_{m,n=0}^N$ yields the relation (2.7).

In the case, when we have orthogonality of polynomials on a
countable set of points a similar conclusion can be false. Let
$p_n(x)$, $n=0,1,2,\cdots$, be a set of orthogonal polynomials
with orthogonality relation
$$
\sum_{m=0}^\infty p_n(x_m)p_{n'}(x_m) w_m=h_n \delta_{nn'}, \eqno(2.8)
$$
where $n,n'=0,1,2,\cdots $ and $h_n$ are some constants. Again,
one may consider $p_n(x_m)$ as functions of $n$. We are interested
in the cases when these functions are polynomials either in $n$ or
in some $\nu(n)$. So the question arises: When the dual relation
to (2.8), namely,
$$
\sum_{n=0}^\infty p_n(x_m)p_{n}(x_{m'})) h^{-1}_n=w^{-1}_m
\delta_{mm'}, \eqno (2.9)
$$
is an orthogonality relation for the dual polynomials
$P_m(\nu(n)):=p_n(x_m)$, $m=0,1,2,\cdots$? It follows from
Riesz's theorem that {\it this is the case, when the orthogonality
measure in (2.8) corresponds to determinate moment problem or
when this measure corresponds to indeterminate moment problem
and it is extremal}. Namely, in both these cases the matrix
$(a_{mn})_{m,n=0}^\infty$ with $a_{mn}=(h_n^{-1}w_m)^{1/2}p_n(x_m)$
is orthogonal, that is,
$$
\sum _m a_{mn}a_{mn'}=\delta_{nn'},\ \ \ \ \ \sum _n
a_{mn}a_{m'n}=\delta_{mm'}.
$$
It is natural to call the polynomials $P_m(\nu(n))$ {\it dual} to
the polynomials $p_n(m)$. The orthogonality relation for them is
$$
\sum_{n=0}^\infty P_m(\nu(n))P_{m'}(\nu(n)) h^{-1}_n=w^{-1}_m
\delta_{mm'}.
$$

However, very often a function $\nu(n)$, such that
$P_m(\nu(n)):=p_n(x_m)$, $m=0,1,2,\cdots$, are polynomials in
$\nu(n)$, does not exist. Nevertheless, sometimes it turns out
that there are some $m$-independent $b_n$, $n=0,1,2,\cdots$,
such that $P_m(\nu(n)):=b_n\, p_n(x_m)$, $m=0,1,2,\cdots$, are
polynomials in $\nu(n)$ for an appropriate function $\nu(n)$.
When the orthogonality measure in (2.8) corresponds to determinate
moment problem or when this measure corresponds to indeterminate
moment problem and it is extremal, then we have the orthogonality
relation (2.9) for the functions ${\hat P}_m(\nu(n)):=p_n(x_m)$,
which is equivalent to the orthogonality relation
$$
\sum_{n=0}^\infty P_m(\nu(n))\, P_{m'}(\nu(n))b^{-2}_n
h^{-1}_n=w^{-1}_m \delta_{mm'}   \eqno (2.10)
$$
for the polynomials $P_m(\nu(n))$. In this case it is also natural
to call the polynomials $P_m(\nu(n))=b_np_n(x_m)$ {\it dual} to
the orthogonal polynomials $p_n(m)$.

The situation can be sometimes more complicated. Namely, the
orthogonality relation for polynomials with a weight function,
supported on a countable set of points, may be of the form (for
instance, for the big $q$-Jacobi polynomials)
$$
\sum_{m=0}^\infty p_n(x_m)p_{n'}(x_m) w_m +\sum_{m=0}^\infty
p_n(y_m)p_{n'}(y_m) w'_m=h_n \delta_{nn'}. \eqno (2.11)
$$
Let $b_n$, $b'_n$, $\nu(n)$ and $\nu'(n)$ be such functions of $n$
that $P_m(\nu(n)):=b_n\, p_n(x_m)$ and
$P'_m(\nu'(n)):=b'_n\,p_n(y_m)$ are polynomials in $\nu(n)$ and
$\nu'(n)$, respectively. When the orthogonality measure in (2.11)
corresponds to determinate moment problem or when this measure
corresponds to indeterminate moment problem and it is extremal,
then the matrix $\left( {(a_{mn})_{m,n=0}^\infty\atop
(a'_{mn})_{m,n=0}^\infty}\right)$, with two infinite matrices
(placed one over another) with matrix elements
$a_{mn}=(h_n^{-1}w_n)^{1/2}p_n(x_m)$ and
$a'_{mn}=(h_n^{-1}w'_n)^{1/2}p_n(y_m)$, is orthogonal, that is,
$$
\sum _m a_{mn}a_{mn'}+\sum _m a'_{mn}a'_{mn'}=\delta_{nn'},
$$  $$
\sum _n a_{mn}a_{m'n}=\delta_{mm'}, \ \ \ \sum _n
a'_{mn}a'_{m'n}=\delta_{mm'}, \ \ \  \sum _n a_{mn}a'_{m'n}=0.
$$
Orthogonality of columns of this matrix gives the orthogonality
relation (2.11). The orthogonality of rows gives orthogonality of
the polynomials $P_m(\nu'(n))$ and $P'_m(\nu'(n))$:
$$
 \sum_{n=0}^\infty P_m(\nu(n))P_{m'}(\nu(n))b^{-2}_n
h^{-1}_n=w^{-1}_m \delta_{mm'}   \eqno (2.12)
$$  $$
 \sum_{n=0}^\infty P'_m(\nu'(n))P'_{m'}(\nu'(n)){b'}^{-2}_n
h^{-1}_n={w'}^{-1}_m \delta_{mm'},   \eqno (2.13)
$$  $$
 \sum_{n=0}^\infty P_m(\nu(n))P'_{m'}(\nu'(n))b^{-1}_n {b'}^{-1}_n
h^{-1}_n=0.   \eqno (2.14)
$$
In this case both sets of the polynomials $P_m(\nu'(n))$ and
$P'_m(\nu'(n))$ are regarded on a par as {\it duals} to the set
of orthogonal polynomials $p_n(x)$. (For the big $q$-Jacobi
polynomials, these two dual sets turn out to be polynomials of
the same type, but with different values of parameters; see
section 4.) If the orthogonality relation (2.11) would contain
$r$ terms, then we had $r$ dual sets of functions.

Thus, if we have a set of polynomials $\{ p_n(x)\}$, orthogonal on
a countable set of points, and they correspond to a determinate
moment problem or to an indeterminate moment problem and the
orthogonality measure is extremal, then one can find the
corresponding orthogonality measure for dual set of polynomials,
if they exist.

The main goal of this review is to discuss a method of constructing
orthogonality measures for a given set of polynomials and their
duals in a straightforward manner. This method is based on the
use of two closed symmetric (self-adjoint) operators, representable
(in some bases) by Jacobi matrices. In the following sections we
shall illustrate how this method works by considering families
of $q$-orthogonal polynomials from the Askey scheme.

Let us emphasize that there are already known theorems on dual
orthogonality properties of polynomials, whose weight functions
are supported on an infinite set of discrete points (see, for
example, [E], [Sz], [BC] and [Ism]). But it is essential that
in most cases (especially in the cases of $q$-polynomials) dual
objects are represented by orthogonal functions. Therefore, one
still needs to make one step further in order to single out an
appropriate family of dual polynomials from these functions (in
those cases when it turns out to be possible). We show that this
step can be made by choosing the $b_n$ for the dual polynomials
$P_m(\nu(n)):=b_n\, p_n(x_m)$. Besides, when one considers some
dual set with respect to a given family of orthogonal polynomials,
it is also necessary to investigate the problem of completeness
for this dual object. In our approach, based on the use of two
particular operators, the problem of completeness is resolved
automatically.
\bigskip

\noindent {\bf 2.4. List of dual sets of $q$-orthogonal
polynomials}

 \bigskip

In this subsection we give a list of dual sets  of
$q$-polynomials, orthogonal on countable sets of points. Each of
these dual pairs will be considered in detail in the subsequent
sections. In particular, orthogonality relations for them will be
explicitly derived.
 \bigskip
 \bigskip

${}\qquad\qquad\qquad$
\begin{tabular}{|c|c|}
 \hline
 $q$-polynomials  & their
 duals\\
 (determinate moment problem) & (indeterminate moment problem) \\
 \hline\hline
 little $q$-Jacobi  &
dual little $q$-Jacobi  \\
 \hline
 big $q$-Jacobi  &
dual big $q$-Jacobi  (two sets)\\
 \hline
 discrete $q$-ultraspherical  &
dual discrete $q$-ultraspherical \\
 \hline
 big $q$-Laguerre  &
$q$-Meixner  (two sets) \\
 \hline
 alternative $q$-Charlier  &
dual alternative $q$-Charlier  \\
\hline
 Al-Salam--Carlitz I  &
$q$-Charlier  (two sets) \\
 \hline
 little $q$-Laguerre  &
Al-Salam--Carlitz II  \\
 \hline
\end{tabular}
\bigskip
\bigskip

\noindent
Let us exhibit these dual pairs explicitly.
\medskip

 \noindent
{\bf Little $q$-Jacobi polynomials and their duals.} Little
$q$-Jacobi polynomials, given by the formula
$$
p_n(\lambda ;a,b|q):= {}_2\phi_1 (q^{-n}, abq^{n+1} ;\; aq; \;
q,q\lambda ) , \eqno (2.15)
$$
are orthogonal for $0<a<q^{-1}$ and $b<q^{-1}$. The dual little
$q$-Jacobi polynomials, corresponding to the polynomials (2.15)
with the same values of the parameters $a$ and $b$, are given as
$$
d_n(\mu (m); a,b|q):= {}_3\phi_1(q^{-m},abq^{m+1},q^{-n};\; bq; \;
q,q^{n}/a ) ,   \eqno (2.16)
$$
where $\mu (m)=q^{-m}+abq^{m+1}$. Since these polynomials are
absent in the Askey $q$-scheme [KSw], we give the orthogonality
relation for these polynomials:
$$
\sum_{m=0}^\infty
\frac{(1-abq^{2m+1})(abq,bq;q)_m}{(1-abq)(aq,q;q)_m} \,a^m\,
q^{m^2}\,d_n(\mu (m))\,d_{n'}(\mu (m))
$$   $$
=\frac{(abq^2;q)_\infty}{(aq;q)_\infty} \frac{(q;q)_n(aq)^{-n}}
{(b;q)_n} \, \delta_{nn'}.
$$
The polynomials $d_n(\mu (m))$ correspond to the indeterminate
moment problem and the ortho\-gonality measure here is extremal.
The duality of these polynomials to the set of the little
$q$-Jacobi polynomials was first observed in [AK1] (see also
[AK3]).

\medskip

\noindent
{\bf Big $q$-Jacobi polynomials and their duals.} Big $q$-Jacobi
polynomials, given by the formula
$$
P_n(\lambda ;a,b,c;q):= {}_3\phi_2 (q^{-n}, abq^{n+1}, \lambda ;\;
aq,cq; \; q,q ) , \eqno (2.17)
$$
are orthogonal for $0<a,b<q^{-1}$ and $c>0$. The dual big
$q$-Jacobi polynomials, associated with the polynomials (2.17)
with the same values of the parameters $a,b,c$, are given as
$$
D_n(\mu (m); a,b,c|q):= {}_3\phi_2(q^{-m},abq^{m+1},q^{-n};\; aq,
abq/c; \; q,aq^{n+1}/c ) ,   \eqno (2.18)
$$
where $\mu (m)=q^{-m}+abq^{m+1}$. The second set of dual
polynomials with respect to (2.17) is obtained from the
polynomials (2.18) by the replacements $a,b,c\to b,a,ab/c$,
respectively. Again, since these polynomials are absent in the
$q$-Askey scheme, we give here the orthogonality relation for the
polynomials (2.18):
$$
\sum_{m=0}^\infty
\frac{(1-abq^{2m+1})(aq,abq,abq/c;q)_m}{(1-abq)(bq,cq,q;q)_m}
\,(-c/a)^m\, q^{m(m-1)/2}\,D_n(\mu (m))\,D_{n'}(\mu (m))
$$   $$
=\frac{(abq^2,c/a;q)_\infty}{(bq,cq;q)_\infty} \frac{(aq/c,q;q)_n}
{(aq,abq/c;q)_nq^n} \delta_{nn'}.
$$
The polynomials $D_n(\mu (m))$ correspond to the indeterminate
moment problem and the orthogonality measure here is not extremal.

\medskip

\noindent {\bf Discrete $q$-ultraspherical polynomials and their
duals.} Discrete $q$-ultraspherical polynomials $C_n^{(a)}(x;q)$,
$a>0$, are a particular case of the big $q$-Jacobi polynomials
$$
C_n^{(a^2)}(x;q)=P_n(x:a,a,-a;q)={}_3\phi_2
(q^{-n},a^2q^{n+1},x;\; aq,-aq;\; q,q).   \eqno (2.19)
$$
An orthogonality relation for $C_n^{(a)}(x;q)$ follows from that
for the big $q$-Jacobi polynomials and it holds for positive
values of $a$. We can consider the polynomials $C_n^{(a)}(x;q)$
also for other values of $a$. In particular, they are orthogonal
for imaginary values of $a$ and $x$. In order to dispense with
imaginary numbers in this case, the following notation is
introduced:
$$
\tilde C_n^{(a^2)}(x;q):= (-{\rm i})^n C_n^{(-a^2)}({\rm i}x;q)
=(-{\rm i})^n P_n({\rm i}x; {\rm i}a,{\rm i}a,-{\rm i}a;q),
 \eqno (2.20)
$$
The orthogonality relation for them is of the form
$$
\sum_{s=0}^\infty \frac{(-aq^2;q^2)_s\, q^s}{(q^2;q^2)_s} \tilde
C_{2k+1}^{(a)}(\sqrt{a}\, q^{s+1};q)\tilde
C_{2k'+1}^{(a)}(\sqrt{a}\, q^{s+1};q) \qquad\qquad
 $$  $$  \qquad\qquad
 =\frac{(-aq^3;q^2)_\infty}
{(q;q^2)_\infty}\frac{(1+aq)\, a^{2k+1}}{(1+aq^{4k+3})}
\frac{(q;q)_{2k+1}}{(-aq;q)_{2k+1}}\, q^{(k+2)(2k+1)}\delta_{kk'}.
$$

The formula
$$
D_n^{(a^2)}(\mu (x;a^2)|q):=D_n(\mu (x;a^2); a,a,-a|q):=\left.
{}_3\phi_2\left({q^{-x},a^2q^{x+1},q^{-n}
  \atop aq, -aq} \right| q,-q^{n+1} \right) ,   \eqno(2.21)
$$
where $\mu(x;a^2)=q^{-x}+a^2q^{x+1}$ and $D_n(\mu (x;a^2))$ are
dual big $q$-Jacobi polynomials, gives dual discrete
$q$-ultraspherical polynomials. They correspond to indeterminate
moment problem. The dual orthogonality relation for them (when
$a^2>0$) follows from the orthogonality relation for dual big
$q$-Jacobi polynomials.

For the polynomials $D_n^{(a^2)}(\mu (x;a^2)|q)$ with imaginary
$a$ we have
$$
\tilde D_n^{(a^2)}(\mu (x;-a^2)|q):=D_n(\mu (x;-a^2);
{\rm i}a,{\rm i}a,-{\rm i}a|q)$$
$$:=\left.
{}_3\phi_2\left({q^{-x},-a^2q^{x+1},q^{-n}
  \atop {\rm i}aq, -{\rm i}aq} \right| q,-q^{n+1} \right) .\eqno(2.22)
$$
These polynomials are dual to the polynomials $\tilde C_n^{(a^2)}(x;q)$
from (2.20). In this case there are also infinitely many orthogonality
relations, which are considered in section 5.
\medskip

\noindent
{\bf Big $q$-Laguerre polynomials and $q$-Meixner polynomials.}
Big $q$-Laguerre polynomials, given by the formula
$$
P_n(\lambda ;a,b;q):= {}_3\phi_2 (q^{-n},0,\lambda;\;
aq,bq;\;q,q),
$$
are orthogonal for $0<a<q^{-1}$ and $b<0$. The dual polynomials
coincide with $q$-Meixner polynomials $M_n(q^{-x};a,-b/a;q)$ and
$M_n(q^{-x};b,-a/b;q)$, where
$$
M_n(q^{-x};a,b;q):={}_2\phi_1 (q^{-n},q^{-x};\; aq;\
q,-q^{n+1}/b).
$$
We obtain orthogonality relations for the $q$-Meixner polynomials
$M_n(q^{-x};a,b;q)$ with $b<0$ and $b>0$ in section 6.

The duality relation between big $q$-Laguerre polynomials and
$q$-Meixner polynomials was studied in [AAK]. The appearance
of $q$-Meixner polynomials as a dual family with respect to the
big $q$-Laguerre polynomials is quite natural because the
transformation $q\to q^{-1}$ interrelates these two sets of
polynomials, that is,
$$
M_n(x;b,c;q^{-1})=(q^{-n}/b;q)_n\, P_n(qx/b;1/b,-c;q).
$$
{\bf Alternative $q$-Charlier polynomials and their duals.}
Alternative $q$-Charlier polynomials are given by the formula
$$
K_n(\lambda ;a;q):={}_2\phi_1 (q^{-n}, -aq^{n};\; 0; \; q,q\lambda)\,.
\eqno (2.23)
$$
They are orthogonal for $a>0$. Their duals are the polynomials
$$
d_n(\mu (m); a;q):= {}_3\phi_0(q^{-m},-a\,q^{m},q^{-n};\; -\, ;\;
q,-q^n/a)\, , \quad\quad  \mu (m):= q^{-m}-a\,q^m \,, \eqno (2.24)
$$
which correspond to the indeterminate moment problem (see [AK5]).
They are also absent in the $q$-Askey scheme. The orthogonality
relation for these polynomials
is
$$
\sum_{m=0}^\infty \frac{(1+aq^{2m})a^m}{(-aq^m;q)_\infty
(q;q)_m}\, q^{m(3m-1)/2} d_n(\mu(m)) d_{n'}(\mu(m))=\frac{
(q;q)_n}{a^nq^{n(n+1)/2}}\, \delta_{nn'} ,\ \ \ a>0.
$$
{\bf Al-Salam--Carlitz I polynomials and $q$-Charlier
polynomials.} Al-Salam--Carlitz I polynomials, given by the
formula
$$
U_n^{(a)}(x;\, q):=(-a)^nq^{n(n-1)/2}\, {}_2 \phi_1
(q^{-n},x^{-1};\ 0 ;\ q;xq/a),
$$
are orthogonal for $a<0$. There are two sets of dual polynomials
[KK]. They coincide with two sets of $q$-Charlier polynomials
$C_n(q^{-x};\, -a;\, q)$ and $C_n(q^{-x};\, -1/a;\, q)$, where
$$
C_n(q^{-x};\, a;\, q):={}_2 \phi_1 (q^{-n},q^{-x};\ 0;\
q;-q^{n+1}/a).
$$
\noindent
{\bf Little $q$-Laguerre polynomials and Al-Salam--Carlitz II
polynomials.} The little $q$-Laguerre polynomials are given by the
formula
$$
p_n(x;\, a|q):={}_2 \phi_1 (q^{-n},0;\ aq;\
q;qx)=(a^{-1}q^{-n};q)_n^{-1}\; {}_2 \phi_0 (q^{-n},x^{-1};\ - ;\
q;x/a).
$$
They are orthogonal for $0<a<q^{-1}$. The dual polynomials with
respect to them are the Al-Salam--Carlitz II polynomials (see
[AK2])
$$
V_n^{(a)}(x;\, q):=(-a)^nq^{-n(n-1)/2}\, {}_2 \phi_0 (q^{-n},x;\ -
;\ q;q^n/a).$$

\bigskip

\noindent{\bf 3. LITTLE $q$-JACOBI POLYNOMIALS AND THEIR DUALS}
 \bigskip

\noindent{\bf 3.1. Pair of operators $(I_1,J)$}
\bigskip

Let ${\cal H}$ be a separable complex Hilbert space with an
orthonormal basis $f_n$, $n = 0,1,2,\cdots$. The basis determines
uniquely a scalar product in ${\cal H}$. In order to deal with a
Hilbert space of functions on a real line, we fix a real number
$a$ such that $0<a<q^{-1}$  and realize our Hilbert space in such
a way that basis elements $f_n$ are monomials:
$$
f_n\equiv f_n(x):= c_n\,x^n,
$$
where
$$
 c^{l}_0 = 1, \qquad c^{l}_n =
 q^{(1-2l)n/4}{(aq;q)_n^{1/2}\over
(q;q)_n^{1/2}}\,, \ \ n = 1,2,3, \cdots .
$$
Thus, in fact, our Hilbert space depends on the number $a$ and can
be denoted as ${\cal H}_a$.

We fix two real parameters $a$ and $b$ such that $b<q^{-1}$,
$0<a<q^{-1}$, and define on ${\cal H}\equiv {\cal H}_a$ two
operators. The first one, denoted as $q^{J_0}$ and taken from the
theory of representations of quantum group $U_q({\rm su}_{1,1})$,
acts on the basis elements as
$$
q^{J_0}f_n =(qa)^{1/2}q^n f_n. \eqno (3.1)
$$
The second operator, denoted as $I_1$, is given by the formula
$$
I_1\, f_n=-a_nf_{n+1}-a_{n-1}f_{n-1}+b_nf_n , \eqno (3.2)
$$
where
$$
a_n=a^{1/2}q^{n+1/2}\frac{ \sqrt{(1-q^{n+1})(1-aq^{n+1})
(1-bq^{n+1})(1-abq^{n+1})}}{(1-abq^{2n+2})\sqrt{(1-abq^{2n+1})
(1-abq^{2n+3})}} ,
$$  $$
b_n=\frac{q^n}{1-abq^{2n+1}}\left(
\frac{(1-aq^{n+1})(1-abq^{n+1})}{1-abq^{2n+2}} +a
\frac{(1-q^{n})(1-bq^{n})}{1-abq^{2n}} \right) .
$$
The expressions for $a_n$ and $b_n$ are well
defined. The operator $I_1$ is symmetric.

Since $a_n\to 0$ and $b_n\to 0$ when $n\to \infty$, the operator
$I_1$ is bounded. Therefore, we assume that it is defined on the
whole space ${\cal H}$. For this reason, $I_1$ is a self-adjoint
operator. Let us show that $I_1$ is a Hilbert--Schmidt operator
(we remind that a bounded self-adjoint operator is a Hilbert--Schmidt
operator if a sum of its squared matrix elements in an orthonormal
basis is finite; the spectrum of such an operator is discrete, with
a single accumulation point at $0$). For the coefficients $a_n$
and $b_n$ from (3.2), we have
$$
a_{n+1}/a_n \to q,\ \ b_{n+1}/b_n \to q \ \ \ {\rm when}\ \ \ n\to
\infty .
$$
Therefore, for the sum of all matrix elements of the operator
$I_1$ in the canonical basis we have $\sum _n (2a_n+b_n)< \infty$.
This means that $I_1$ is a Hilbert--Schmidt operator. Thus, the
spectrum of $I_1$ is discrete and has a single accumulation point
at 0. Moreover, a spectrum of $I_1$ is simple, since $I_1$ is
representable by a Jacobi matrix with $a_n\ne 0$ (see [Ber],
Chapter VII).

To find eigenfunctions $\xi_\lambda (x)$ of the operator $I_1$,
$I_1 \xi_\lambda (x)=\lambda \xi_\lambda (x)$, we set
$$
\xi_\lambda (x)=\sum _{n=0}^\infty \beta_n(\lambda)f_n (x).
$$
Acting by the operator $I_1$ upon both sides of this relation, one
derives that
$$
\sum _{n=0}^{\infty}\, \beta_n(\lambda)\,
(a_nf_{n+1}+a_{n-1}f_{n-1}-b_nf_n) = - \lambda
\sum_{n=0}^{\infty}\, \beta_n(\lambda)f_n ,
$$
where $a_n$ and $b_n$ are the same as in (3.2). Collecting in
this identity all factors, which multiply $f_n$ with fixed $n$,
one derives the recurrence relation for the coefficients
$\beta_n(\lambda)$:
$$
\beta_{n+1}(\lambda)a_n +\beta_{n-1}(\lambda)a_{n-1}-
\beta_{n}(\lambda)b_n= -\lambda \beta_{n}(\lambda).
$$
The substitution
$$
\beta_{n}(\lambda)=\left( \frac{(abq,aq;q)_n\,(1-abq^{2n+1})}
{(bq,q;q)_n\, (1-abq)(aq)^n}\right) ^{1/2} \beta'_{n}(\lambda)
$$
reduces this relation to the following one
$$
A_n \beta'_{n+1}(\lambda)+ C_n \beta'_{n-1}(\lambda)
-(A_n+C_n)\beta'_{n}(\lambda)=-\lambda \beta'_{n}(\lambda)
$$
with
$$
A_n=\frac{q^n(1-aq^{n+1})(1-abq^{n+1})}{(1-abq^{2n+1})
(1-abq^{2n+2})}, \ \ \ \
C_n=\frac{aq^n(1-q^{n})(1-bq^{n})}{(1-abq^{2n}) (1-abq^{2n+1})}.
$$
This is the recurrence relation for the little $q$-Jacobi
polynomials
$$
p_n(\lambda ;a,b|q):={}_2\phi_1 (q^{-n}, abq^{n+1};\; aq; \;
q,q\lambda )   \eqno (3.3)
$$
(see, for example, formula (7.3.1) in [GR]). Therefore,
$\beta'_n(\lambda )=p_n(\lambda ;a,b|q)$ and
$$
\beta_n(\lambda )= \left( \frac{(abq,aq;q)_n\,(1-abq^{2n+1})}
{(bq,q;q)_n\, (1-abq)(aq)^n}\right)^{1/2}p_n(\lambda ;a,b|q).
\eqno (3.4)
$$
For the eigenfunctions $\xi _\lambda(x)$ we have the expression
$$
\xi _\lambda(x)=\sum_{n=0}^\infty\,\left( \frac{(abq,aq;q)_n\,
(1-abq^{2n+1})} {(bq,q;q)_n\,
(1-abq)(aq)^n}\right)^{1/2}p_n(\lambda ;a,b|q)\,f_n(x)
\qquad\qquad\qquad
$$   $$   \qquad\qquad\qquad
=\sum_{n=0}^\infty\,a^{-n/4}\,\frac{(aq;q)_n}{(q;q)_n}\,
\,\left(\frac{(abq;q)_n\,(1-abq^{2n+1})} {(bq;q)_n\,
(1-abq)(aq)^n}\right)^{1/2} p_n(\lambda ;a,b|q) x^n. \eqno (3.5)
$$
Since the spectrum of the operator $I_1$ is discrete, only a
discrete set of these functions belongs to the Hilbert space
${\cal H}$. This discrete set of functions determines a spectrum
of $I_1$.

Now we look for a spectrum of the operator $I_1$ and for a set of
polynomials, dual to the little $q$-Jacobi polynomials. To this
end we use the action of the operator
$$
J:= (qa)^{1/2}q^{-J_0} +(qa)^{-1/2} ab\,q^{J_0+1}
$$
upon the eigenfunctions $\xi _\lambda(x)$, which belong to
the Hilbert space ${\cal H}$. In order to find how this
operator acts upon these functions, one can use the
$q$-difference equation
$$
(q^{-n}+abq^{n+1})\,p_n(\lambda)= a\lambda^{-1}(bq\lambda
-1)\,p_{n} (q\lambda
)+\lambda^{-1}(1+a)\,p_n(\lambda)+\lambda^{-1}
(\lambda-1)\,p_n(q^{-1}\lambda) \eqno(3.6)
$$
for the little $q$-Jacobi polynomials $p_n(\lambda)\equiv
p_n(\lambda ;a,b|q)$ (see, for example, formula (3.12.5)
in [KSw]). Multiply both sides of (3.6) by $d_n\,f_n(x)$
and sum up over $n$, where $d_n$ are the coefficients of
$p_n(\lambda ;a,b|q)$ in the expression (3.4) for the
$\beta_n(\lambda)$. Taking into account the first line
in formula (3.5) and the fact that $Jf_n(x)=(q^{-n}+abq^{n+1})f_n(x)$,
one obtains the relation
$$
J\,\xi _{\lambda}(x)= a\lambda^{-1}(bq\lambda -1)\,\xi
_{q\lambda}(x) + \lambda^{-1}(1+a)\, \xi _{\lambda}(x)+
\lambda^{-1}(\lambda-1)\, \xi_{q^{-1}\lambda}(x). \eqno (3.7)
$$

It will be shown in the next section that the spectrum of
the operator $I_1$ consists of the points $\lambda=q^n$,
$n=0,1,2,\cdots$. Thus, we see that the pair of the operators
$I_1$ and $J$ form a Leonard pair (see [Ter1], where
P.~Terwilliger has actually introduced this notion in an
effort to interpret the results of D. Leonard [Leo]; see also
[Ter3], which contains a review on how one can employ Leonard
pairs to describe properties of orthogonal polynomials). We
remind to the reader that a pair of operators $R_1$ and $R_2$,
acting on a linear space ${\cal L}$, is a Leonard pair if

(a) there exists a basis in ${\cal L}$, with respect to which
the operator $R_1$ is diagonal, and the operator $R_2$ has the
form of a Jacobi matrix;

(b) there exists another basis of ${\cal L}$, with respect to
which the operator $R_2$ is diagonal, and the operator $R_1$ has
the form of a Jacobi matrix.

Properties of Leonard pairs of operators in finite dimensional
spaces are studied in detail. Leonard pairs in infinite
dimensional spaces are more complicated and only some isolated
results are known in this case (see, for example, [Ter2]).
\bigskip

\noindent{\bf 3.2. Spectrum of $I_1$ and orthogonality of little
$q$-Jacobi polynomials}
\medskip

The aim of this section is to find, by using the Leonard pair
$(I_1,J)$, a basis in the Hilbert space ${\cal H}$, which consists
of eigenfunctions of the operator $I_1$ in a normalized form, and
to derive explicitly the unitary matrix $U$, connecting this basis
with the canonical basis $f_n$, $n=0,1,2,\cdots$, in ${\cal H}$.
This matrix directly leads the orthogonality relation for the
little $q$-Jacobi polynomials.

Let us analyze a form of spectrum of the operator $I_1$. If
$\lambda$ is a spectral point of the operator $I_1$, then (as it
is easy to see from (3.7)) a successive action by the operator $J$
upon the function (eigenfunction of $I_1$) $\xi_\lambda$ leads to
the functions
$$
\xi_{q^m\lambda}, \ \ \ m=0,\pm 1, \pm 2,\cdots , \eqno (3.8)
$$
which are eigenfunctions of $I_1$ with eigenvalues $q^m\lambda$.
However, since $I_1$ is a trace class operator, not all these
points can belong to the spectrum of $I_1$, since $q^{-m}\lambda
\to\infty$ when $m\to \infty$ if $\lambda\ne 0$. This means that
under a successive action by $I_1$ upon $\xi_\lambda$, on some
step the last term in (3.7) must vanish. Thus, under the action by
$I_1$ upon $\xi_{\lambda'}$ for some $\lambda'$ the coefficient
$\lambda' -1$ of $\xi _{q^{-1}\lambda'}(x)$ in (3.7) vanishes.
Clearly, it vanishes when $\lambda' =1$. Moreover, this is the
only possibility for the coefficient of $\xi _{q^{-1}\lambda'}(x)$
in (3.7) to vanish, that is, the point $\lambda =1$ is a spectral
point for the operator $I_1$. Let us show that the corresponding
eigenfunction $\xi _{1}(x)\equiv \xi_{q^{0}}(x)$ belongs to the
Hilbert space ${\cal H}$.

Observe that by formula (II.6) of Appendix II in [GR], one has
$$
p_n(1 ;a,b|q)={}_2\phi_1 (q^{-n}, abq^{n+1};\; aq; \; q,q
)=\frac{(b^{-1}q^{-n};q)_n}{(aq;q)_n} (abq^{n+1})^n.
$$
Since $(b^{-1}q^{-n};q)_n=(bq;q)_n
(-b^{-1}q^{-1})^nq^{-n(n-1)/2}$, this means that
$$
p_n(1 ;a,b|q)=\frac{(bq;q)_n}{(aq;q)_n}(-a)^nq^{n(n+1)/2}.
$$
Therefore, due to (3.5) for the scalar product $\langle
\xi_1(x),\xi_1(x)\rangle$ we have
$$
\langle \xi_1(x),\xi_1(x)\rangle = \sum_{n=0}^\infty
\frac{(abq,aq;q)_n\,(1-abq^{2n+1})}{(bq,q;q)_n\,
(1-abq)(aq)^n}\,p^2_n(1 ;a,b|q)  \qquad\qquad\quad
$$   $$  \qquad\qquad\quad
= \sum_{n=0}^\infty
\frac{(abq,bq;q)_n\,(1-abq^{2n+1})}{(aq,q;q)_n(1-abq)}\,
a^nq^{n^2} = \frac{(abq^2;q)_\infty}{(aq;q)_\infty} . \eqno (3.9)
$$
The last leg of this equality is obtained from formula (A.1) of
Appendix. Thus, the series (3.9) converges and, therefore, the
point $\lambda =1$ actually belongs to the spectrum of the
operator $I_1$.

Let us find other spectral points of the operator $I_1$ (recall
that a spectrum of $I_1$ is discrete). Setting $\lambda = 1$ in
(3.7), we see that the operator $J$ transforms $\xi _{q^0}(x)$
into a linear combination of the functions $\xi _{q}(x)$ and
$\xi_{q^0}(x)$. Moreover, $\xi_q(x)$ belongs to the Hilbert space
${\cal H}$, since the series
$$
\langle \xi _{q} ,\xi _{q} \rangle = \sum_{n=0}^\infty
\frac{(abq,aq;q)_n\,(1-abq^{2n+1})}{(bq,q;q)_n\, (1-abq)(aq)^n}
p^2_n(q ;a,b|q) <\infty
$$
is majorized by the corresponding series for $\xi_{q^0}(x)$,
considered above. Therefore, $\xi _{q}(x)$ belongs to the Hilbert
space ${\cal H}$ and the point $q$ is an eigenvalue of the
operator $I_1$. Similarly, setting $\lambda=q$ in (3.7), we find
that $\xi _{q^2}(x)$ is an eigenfunction of $I_1$ and the point
$q^2$ belongs to the spectrum of $I_1$. Repeating this procedure,
we find that $\xi _{q^n}(x)$, $n=0,1,2,\cdots$, are eigenfunctions
of $I_1$ and the set $q^n$, $n=0,1,2,\cdots$, belongs to the
spectrum of $I_1$. So far, we do not know yet whether other
spectral points exist or not.

The functions $\xi _{q^n}(x)$, $n=0,1,2,\cdots$, are linearly
independent elements of the space ${\cal H}$ (since they
correspond to different eigenvalues of the self-adjoint operator
$I_1$). Suppose that values $q^n$, $n=0,1,2,\cdots$, constitute a
whole spectrum of the operator $I_1$. Then the set of functions
$\xi _{q^n}(x)$, $n=0,1,2,\cdots$, is a basis in the Hilbert space
${\cal H}$. Introducing the notation $\Xi _n:=\xi_{q^n}(x)$,
$n=0,1,2,\cdots$, we find from (3.7) that
$$
J \,\Xi _n = - a q^{-n}(1-bq^{n+1})\,\Xi _{n+1} + q^{-n}(a+1)\,
\Xi _n - q^{-n}(1-q^n)\, \Xi _{n-1} . \eqno (3.10)
$$
As we see, the matrix of the operator $J$ in the basis $\Xi _n$,
$n=0,1,2,\cdots$, is not symmetric, although in the initial basis
$f_n$, $n=0,1,2,\cdots$, it was symmetric. The reason is that the
matrix $A\equiv (a_{mn})$ with entries
$$
a_{mn}:=\beta_m(q^n),\ \ \ \ m,n=0,1,2,\cdots ,
$$
where $\beta_m(q^n)$ are the coefficients (3.4) in the expansion
$\xi _{q^n}(x)=\sum _m \,\beta_m(q^n)f_m(x)$, is not unitary.
(This matrix connects the bases $\{ f_n\}$ and $\{\Xi_n\}$.) It is
equivalent to the statement that the basis $\Xi _n:
=\xi_{q^n}(x)$, $n=0,1,2,\cdots$, is not normalized. To normalize
it, one has to multiply $\Xi _n$ by corresponding numbers $c_n$
(which are not known at this moment). Let $\hat\Xi _n = c_n\Xi
_n$, $n=0,1,2,\cdots$, be a normalized basis. Then the matrix of
the operator $J$ is symmetric in this basis. Since  $J$ has in the
basis $\{ \hat\Xi _n\}$ the form
$$
J\, \hat\Xi _n = -c_{n+1}^{-1}c_naq^{-n}(1-bq^{n+1})\, \hat\Xi
_{n+1} + q^{-n}(a+1)\, \hat\Xi _n - c_{n-1}^{-1}c_n q^{-n}(1-q^n)
\,\hat\Xi_{n-1} , \eqno (3.11)
$$
then its symmetricity means that
$$
c_{n+1}^{-1}c_naq^{-n}(1-bq^{n+1})=c_{n}^{-1}c_{n+1}
q^{-n-1}(1-q^{n+1})\,,
$$
that is, $c_{n}/c_{n-1} =\sqrt{aq (1-bq^n)/(1-q^n)}$. Therefore,
$$
c_n= c(aq)^{n/2}\frac{(bq;q)_n^{1/2}}{(q;q)_n^{1/2}},
$$
where $c$ is a constant.

Now instead of the expansion (3.5) we have the expansions
$$
\hat\xi _{q^n}(x)\equiv \hat\Xi _n(x)=\sum _m
c_n\beta_m(q^n)f_m(x), \eqno (3.12)
$$
which connect two orthonormal bases in the space ${\cal H}$. This
means that the matrix $({\hat a}_{mn})$, $m,n=0,1,2,\cdots$, with
entries
$$
{\hat a}_{mn}=c_n\beta _m(q^n)= c\left( (aq)^{n-m}\,
\frac{(bq;q)_n}{(q;q)_n}\,\frac{(abq,aq;q)_m\,(1-abq^{2m+1})}
{(bq,q;q)_m\, (1-abq) }\right) ^{1/2}\, p_m(q^n ;a,b|q),
\eqno(3.13)
$$
is unitary, provided  that the constant $c$ is appropriately
chosen. In order to calculate this constant, we use the relation
$\sum_{m=0}^\infty |{\hat a}_{mn}|^2=1$ for $n=0$. Then the sum in
this relation is a multiple of the sum in (3.9) and, consequently,
$$
c=\frac{(aq;q)^{1/2}_\infty}{(abq;q)^{1/2}_\infty} .
$$
Thus the $c_n$ in (3.12) and (3.13) are real and equal to
$$
c_n=\left( \frac{(aq;q)_\infty}{(abq;q)_\infty}
\frac{(bq;q)_n(aq)^n}{(q;q)_n} \right) ^{1/2} .
$$

The matrix $({\hat a}_{mn})$ with entries (3.13) is orthogonal,
that is,
$$
\sum _n {\hat a}_{mn}{\hat a}_{m'n}=\delta_{mm'},\ \ \ \ \sum _m
{\hat a}_{mn}{\hat a}_{mn'}=\delta_{nn'} . \eqno (3.14)
$$
Substituting into the first sum over $n$ in (3.14) the expressions
for ${\hat a}_{mn}$, we obtain the identity
$$
\sum_{n=0}^\infty
\frac{(bq;q)_n(aq)^n}{(q;q)_n}\,p_m(q^n;a,b|q)\,p_{m'}(q^n;a,b|q)
$$  $$
=\frac{(abq^2;q)_\infty}{(aq;q)_\infty}
\frac{(1-abq)(aq)^m\,(bq,q;q)_m} {(1-abq^{2m+1})\,(abq,aq;q)_m}\,
\delta_{mm'}\,, \eqno (3.15)
$$
which must yield the orthogonality relation for the little
$q$-Jacobi polynomials. An only gap, which appears here, is the
following. We have assumed that the points $q^n$,
$n=0,1,2,\cdots$, exhaust the whole spectrum of the operator
$I_1$. Let us show that this is the case.

Recall that the self-adjoint operator $I_1$ is represented by a
Jacobi matrix in the basis $f_n$, $n=0,1,2,\cdots$. According to
the theory of operators of such type (see, for example, [Ber],
Chapter VII; a short explanation is given in section 2),
eigenfunctions $\xi_\lambda$ of $I_1$ are expanded into series in
the monomials $f_n$, $n=0,1,2,\cdots$, with coefficients, which
are polynomials in $\lambda$. These polynomials are orthogonal
with respect to some positive measure $d\mu (\lambda)$ (moreover,
for self-adjoint operators this measure is unique). The set (a
subset of ${\Bbb R}$), on which the polynomials are orthogonal,
coincides with the spectrum of the operator under consideration
and the spectrum is simple. Let us apply these assertions to the
operator $I_1$.

We have found that the spectrum of $I_1$ contains the points
$q^n$, $n=0,1,2,\cdots$. If the operator $I_1$ had other spectral
points $x$, then on the left-hand side of (3.15) there would be
other summands $\mu_{x_k}\,
p_m({x_k};a,b|q)\,p_{m'}({x_k};a,b|q)$, corresponding to these
additional points. Let us show that these additional summands do
not appear. To this end we set $m=m'=0$ in the relation (3.15)
with the additional summands. Since $p_0(x;a,b|q)=1$, we have the
equality
$$
\sum_{n=0}^\infty \frac{(bq;q)_n(aq)^n}{(q;q)_n} + \sum_k
\mu_{x_k} =\frac{(abq^2;q)_\infty}{(aq;q)_\infty} .
$$
According to the $q$-binomial theorem (see formula (1.3.2) in
[GR]), we have
$$
\sum_{n=0}^\infty\,\frac{(bq;q)_n\,(aq)^n}{(q;q)_n} =
\frac{(abq^2;q)_\infty}{(aq;q)_\infty}. \eqno (3.16)
$$
Hence, $\sum_k \mu_{x_k} =0$ and this means that additional
summands do not appear in (3.15) and it does represent the
orthogonality relation for the little $q$-Jacobi polynomials.

By using the operators $I_1$ and $J$, which form a Leonard pair of
infinite dimensional symmetric operators, we thus derived the
orthogonality relation for little $q$-Jacobi polynomials.

The orthogonality relation for the little $q$-Jacobi polynomials
is given by formula (3.15). Due to this orthogonality, we arrive at
the following statement: {\it The spectrum of the operator $I_1$
coincides with the set of points $q^{n}$, $n=0,1,2,\cdots$. The
spectrum is simple and has one accumulation point at 0.}
\bigskip

\noindent {\bf 3.3. Dual little $q$-Jacobi polynomials}
\bigskip

Now we consider the second identity in (3.14), which gives the
orthogonality relation for the matrix elements ${\hat a}_{mn}$,
considered as functions of $m$. Up to multiplicative factors these
functions coincide with
$$
F_n(x;a,b|q):={}_2\phi_1 (x,abq/x;\; aq;\; q,q^{n+1}), \eqno
(3.17)
$$
considered on the set $x\in \{ q^{-m}\, |\, m=0,1,2,\cdots \}$.
Consequently,
$$
{\hat a}_{mn}= \left( \frac{(aq;q)_\infty}{(abq;q)_\infty}\,
\frac{(bq;q)_n}{(q;q)_n}\,(aq)^{n-m}\,\frac{(abq,aq;q)_m\,(1-abq^{2m+1})}
{(bq,q;q)_m}\right)^{1/2}\, F_n(q^{-m} ;a,b|q)
$$
and the second identity in (3.14) gives the orthogonality relation
for the functions (3.17):
$$
\sum_{m=0}^\infty \,\frac{(1-abq^{2m+1})\,(abq,aq;q)_m}
{(1-abq)(aq)^m\,(bq,q;q)_m}\,
F_n(q^{-m};a,b|q)\,F_{n'}(q^{-m};a,b|q)
$$  $$
= \frac{(abq^2;q)_\infty}{(aq;q)_\infty}
\frac{(q;q)_n(aq)^{-n}}{(bq;q)_n} \, \delta_{nn'}. \eqno (3.18)
$$

The functions $F_n(x;a,b|q)$ can be represented in another form.
Indeed, one can use the relation (III.8) of Appendix III in [GR]
in order to obtain that
$$
F_{n}(q^{-m} ;a,b|q)=\frac{(b^{-1}q^{-m};q)_m}{(aq;q)_m}
(abq^{m+1})^m  {}_3\phi_1 (q^{-m},abq^{m+1},q^{-n};\; bq; \;
q,q^n/a)
$$  $$
=\frac{(-1)^m\,(bq;q)_m}{(aq;q)_m}\, a^m \,q^{m(m+1)/2} {}_3\phi_1
(q^{-m},abq^{m+1},q^{-n};\; bq; \; q,q^n/a).  \eqno (3.19)
$$
The basic hypergeometric function ${}_3\phi_1$ in (3.19) is a
polynomial of degree $n$ in the variable $\mu(m): = q^{-m}+
ab\,q^{m+1}$, which represents a $q$-quadratic lattice; we denote
it as
$$
d_n(\mu (m); a,b|q):= {}_3\phi_1(q^{-m},ab\,q^{m+1},q^{-n};\;
bq;\; q,q^n/a)\,.  \eqno (3.20)
$$
Then formula (3.18) yields the orthogonality relation
$$
\sum_{m=0}^\infty \,\frac{(1-abq^{2m+1})(abq,bq;q)_m}
{(1-abq)(aq,q;q)_m}\,a^m\, q^{m^2}\, d_n(\mu(m))\, d_{n'}(\mu (m))
$$  $$
= \frac{(abq^2;q)_\infty}{(aq;q)_\infty}
\,\frac{(q;q)_n(aq)^{-n}}{(bq;q)_n}\, \delta_{nn'}  \eqno (3.21)
$$
for the polynomials (3.20). We call the polynomials $d_n(\mu (m);
a,b|q)$ {\it dual little $q$-Jacobi polynomials}.

Note that these polynomials can be expressed in terms of the
Al-Salam--Chihara polynomials
$$
Q_n(x;a,b|q)=\frac{(ab;q)_n}{a^n} \left.
{}_3\phi_2\left({q^{-n},az,az^{-1} \atop ab, 0} \right| q,q
\right) , \ \ \ \ x=\frac 12 (z+z^{-1}),
$$
with the parameter $q>1$. An explicit relation between them is
$$
d_n(\mu(x); \beta /\alpha,\, 1/\alpha \beta q \, |\,
q)=\frac{q^{n(n-1)/2}}{(-\beta)^n(1/\alpha\beta ;q)_n}\,
Q_n(\alpha\mu(x)/2; \alpha,\beta | q^{-1}).
$$
Ch. Berg and M. E. H. Ismail studied this type of
Al-Salam--Chihara polynomials in [BI] and derived complex
orthogonality measures for them. But [BI] does not contain any
discussion of the duality of this family of polynomials with
respect to little $q$-Jacobi polynomials.

Observe that the dual polynomials (3.20) can be also expressed in
terms of the little $q$-Jacobi polynomials (3.3):
$$
d_n(\mu(m);
a,b|q)=\frac{(aq;q)_m}{(bq/c;q)_m}(-a)^{-m}q^{-m(m+1)/2}
p_m(q^n;a,b|q).
$$

A recurrence relation for the polynomials $d_n(\mu (m);a,b|q)$ is
derived from formula (3.6). It has the form
$$
(q^{-m}+ ab q^{m+1})\,d_n(\mu (m))= -\, a\,
q^{-n}(1-bq^{n+1})\,d_{n+1}(\mu (m))
$$   $$
+\, q^{-n}(1+a)\,d_{n}(\mu (m)) - q^{-n}(1-q^{n})\, d_{n-1}(\mu
(m)),
$$
where $d_{n}(\mu (m))\equiv d_{n}(\mu(m);a,b|q)$. Comparing this
relation with the recurrence relation (3.69) in [AI1], we see that
the polynomials (3.20) are multiple to the polynomials (3.67) in
[AI1]. Moreover, if one takes into account this multiplicative
factor, the orthogonality relation (3.21) for polynomials (3.20)
turns into relation (3.82) for the polynomials (3.67) in [AI1],
although the derivation of the orthogonality relation in [AI1] is
more complicated than our derivation of (3.21). The authors of
[AI1] do not give an explicit form of their polynomials in the
form similar to (3.20). Concerning the polynomials (3.67) in [AI1]
see also [Gro].

Let ${\frak l}^2$ be the Hilbert space of functions on the set
$m=0,1,2,\cdots$ with the scalar product
$$
\langle f_1,f_2\rangle = \sum _{m=0}^\infty\,
\frac{(1-abq^{2m+1})\, (abq,bq;q)_m}{(1-abq)\, (aq,q;q)_m}\,
a^m\,q^{m^2}\,f_1(m)\,\overline{f_2(m)} ,
$$
where weight function is taken from (3.21). The polynomials (3.20)
are in one-to-one correspondence with the columns of the unitary
matrix $({\hat a}_{mn})$ and the orthogonality relation (3.21) is
equivalent to the orthogonality of these columns. Due to (3.14)
the columns of the matrix $({\hat a}_{mn})$ form an orthonormal
basis in the Hilbert space of sequences ${\bf a}=\{ a_n\, |\,
n=0,1,2,\cdots \}$ with the scalar product $\langle {\bf a},{\bf
a}'\rangle=\sum_n a_na'_n$. This assertion is equivalent to the
following one: the set of polynomials $d_n(\mu (m);a,b|q)$,
$n=0,1,2,\cdots$, form an orthogonal basis in the Hilbert space
${\frak l}^2$. This means that {\it the point measure in (3.21) is
extremal for the dual little $q$-Jacobi polynomials} $d_n(\mu
(m);a,b|q)$.
\bigskip

\noindent{\bf 4. BIG $q$-JACOBI POLYNOMIALS AND THEIR DUALS}
 \bigskip

\noindent{\bf 4.1. Pair of operators $(I_2,J)$}
\bigskip

We fix three real numbers $a$, $b$ and $c$ such that $0<a<q^{-1}$,
$0<b<q^{-1}$, $c<0$ and consider on the Hilbert space ${\cal
H}\equiv {\cal H}_a$, introduced in subsection 3.1, the following
symmetric operator $I_2$:
$$
I_2\, f_n=a_nf_{n+1}+a_{n-1}f_{n-1}-b_nf_n , \eqno (4.1)
$$
where
$$
a_{n-1}=(-acq^{n+1})^{1/2}\,\frac{ \sqrt{(1-q^{n})(1-aq^{n})
(1-bq^{n})(1-abq^{n})\,(1-cq^{n})(1-abc^{-1}q^{n})}}
{(1-abq^{2n})\,\sqrt{(1-abq^{2n-1})(1-abq^{2n+1})}} ,
$$  $$
b_n= \frac{(1-aq^{n+1})(1-abq^{n+1})(1-cq^{n+1})}
{(1-abq^{2n+1})(1-abq^{2n+2})}-acq^{n+1}
\frac{(1-q^{n})(1-bq^{n})(1-abq^{n}/c)}{(1-abq^{2n})
(1-abq^{2n+1})} - 1\,.
$$
This operator is bounded. Therefore, we assume that it is defined
on the whole Hilbert space ${\cal H}$. This means that $I_2$ is a
self-adjoint operator. Actually, $I_2$ is a Hilbert--Schmidt
operator. To show this we note that for the coefficients $a_n$ and
$b_n$ from (4.1) one obtains that
$$
a_{n+1}/a_n \to q^{1/2},\ \ b_{n+1}/b_n \to q \ \ \ {\rm when}\ \
\ n\to \infty .
$$
Therefore, $\sum _n (2a_n+b_n)< \infty$ and this means that $I_2$
is a Hilbert--Schmidt operator. Thus, the spectrum of $I_2$ is
simple (since it is representable by a Jacobi matrix with $a_n\ne
0$), discrete and have a single accumulation point at 0.

To find eigenfunctions $\psi_\lambda (x)$ of the operator $I_2$,
$I_2 \psi_\lambda (x)=\lambda \psi_\lambda (x)$, we set
$$
\psi_\lambda (x)=\sum _{n=0}^\infty \beta_n(\lambda)f_n (x). \eqno
(4.2)
$$
Acting by the operator $I_2$ on both sides of this relation, one
derives that
$$
\sum_n \beta_n(\lambda) (a_nf_{n+1}+a_{n-1}f_{n-1}-b_nf_n)
=\lambda \sum \beta_n(\lambda)f_n ,
$$
where $a_n$ and $b_n$ are the same as in (4.1). Collecting in this
identity factors, which multiply $f^l_n$ with fixed $n$, we arrive
at the recurrence relation for the coefficients
$\beta_n(\lambda)$:
$$
a_n\beta_{n+1}(\lambda) +a_{n-1}\beta_{n-1}(\lambda)-
b_n\beta_{n}(\lambda) = \lambda \beta_{n}(\lambda).
$$
Making the substitution
$$
\beta_{n}(\lambda)=\left(\frac{(abq,aq,cq;q)_n\,
(1-abq^{2n+1})}{(abq/c,bq,q;q)_n\,
(1-abq)(-ac)^n}\right)^{1/2}q^{-n(n+3)/4} \,\beta'_{n}(\lambda)
$$
we reduce this relation to the following one
$$
A_n \beta'_{n+1}(\lambda)+ C_n \beta'_{n-1}(\lambda)
-(A_n+C_n-1)p'_{n}(\lambda)=\lambda \beta'_{n}(\lambda)
$$
with
$$
A_n=\frac{(1-aq^{n+1})(1-cq^{n+1})(1-abq^{n+1})}
{(1-abq^{2n+1})\,(1-abq^{2n+2})},
$$  $$
C_n=\frac{-acq^{n+1}(1-q^{n})(1-bq^{n})(1-abc^{-1}q^{n})}
{(1-abq^{2n})\,(1-abq^{2n+1})}.
$$
It is the recurrence relation for the big $q$-Jacobi polynomials
$$
P_n(\lambda ;a,b,c;q):= {}_3\phi_2 (q^{-n}, abq^{n+1}, \lambda ;\;
aq,cq; \; q,q )  \eqno (4.3)
$$
introduced by G. E. Andrews and R. Askey [AA] (see also formula
(7.3.10) in [GR]). Therefore, $\beta'_n(\lambda )=P_n(\lambda
;a,b,c;q)$ and
$$
\beta_n(\lambda )=\left(\frac{(abq,aq,cq;q)_n\,(1-abq^{2n+1})}
{(abq/c,bq,q;q)_n\, (1-abq) (-ac)^n} \right)^{1/2}q^{-n(n+3)/4}\,
P_n(\lambda ;a,b,c;q). \eqno (4.4)
$$
For the eigenfunctions $\psi _\lambda(x)$ we have the expansion
$$
\psi _\lambda(x)=\sum_{n=0}^\infty \left(\frac{(abq,aq,cq;q)_n\,
(1-abq^{2n+1})}{(abq/c,bq,q;q)_n\,
(1-abq)(-ac)^n}\right)^{1/2}q^{-n(n+3)/4}\, P_n(\lambda
;a,b,c;q)\, f_n(x)
$$   $$
=\sum_{n=0}^\infty \,a^{-n/4}\, \frac{(aq;q)_n}{(q;q)_n}
\left(\frac{(abq,cq;q)_n\,(1-abq^{2n+1})} {(abq/c,bq;q)_n\,
(1-abq) (-ac)^n} \right)^{1/2}q^{-n(n+3)/4}\, P_n(\lambda
;a,b,c;q)\, x^n. \eqno (4.5)
$$
Since the spectrum of the operator $I_2$ is discrete, only a
discrete set of these functions belongs to the Hilbert space
${\cal H}$.

In what follows we intend to study a spectrum of the operator
$I_2$ and to find polynomials, dual to big $q$-Jacobi polynomials.
It can be done with the aid of the operator
$$
J:= (aq)^{1/2}q^{-J_0} +(aq)^{-1/2} ab\,q^{J_0+1},
$$
which has been already used in the previous case in subsection
3.1. In order to determine how this operator acts upon the
eigenfunctions $\psi _\lambda(x)$, one can use the $q$-difference
equation
$$
(q^{-n}+abq^{n+1})\, P_n(\lambda)= aq \lambda^{-2}( \lambda
-1)(b\lambda -c)\,P_{n}(q\lambda )
$$
$$
- [\lambda^{-2}acq(1+q)-\lambda^{-1}q(ab+ac+a+c)]\, P_n(\lambda)+
\lambda^{-2}(\lambda-aq)(\lambda-cq)\,P_n(q^{-1}\lambda) ,
\eqno(4.6)
$$
for the big $q$-Jacobi polynomials $P_n(\lambda)\equiv P_n
(\lambda;a,b,c;q)$ (see, for example, formula (3.5.5) in [KSw]).
Multiply both sides of (4.6) by $k_n\,f_n(x)$, where $k_n$ are the
coefficients of $P_n(\lambda ;a,b,c;q)$ in the expression (4.4)
for the coefficients $\beta_n(\lambda)$, and sum over $n$. Taking
into account formula (4.5) and the fact that
$J\,f_n(x)=(q^{-n}+ab\,q^{n+1})\,f_n(x)$, one obtains the relation
$$
J\, \psi _{\lambda}(x)= aq\lambda^{-2}( \lambda -1)(b\lambda -c)
\,\psi_{q\lambda}(x)
$$
$$
- [\lambda^{-2}acq(1+q)-\lambda^{-1}q(ab+ac+a+c)]\,\psi
_{\lambda}(x) + \lambda^{-2}(\lambda-aq)(\lambda-cq)\,\psi
_{q^{-1}\lambda}(x) .\eqno (4.7)
$$
It will be shown in the next section that the spectrum of the
operator $I_2$ consists of the points  $aq^n$, $cq^n$,
$n=0,1,2,\cdots$. The matrix of the operator $J$ in the basis of
eigenfunctions of $I_2$ consists of two Jacobi matrices (one
corresponds to the spectral points $aq^n$, $n=0,1,2,\cdots$, and
another to the spectral points $cq^n$, $n=0,1,2,\cdots$). In this
case, the operators $I_2$ and $J$ form some generalization of
Leonard pair.
\bigskip

\noindent{\bf 4.2. Spectrum of $I_2$ and orthogonality of big
$q$-Jacobi polynomials}
\bigskip

As in subsection 3.2 one can show that for some value of $\lambda$
(which must belong to the spectrum) the last term on the right
side of (4.7) has to vanish. There are two such values of
$\lambda$: $\lambda = aq$ and $\lambda = cq$. Let us show that
both of these points are spectral points of the operator $I_2$.
Observe that, according to (4.3),
$$
P_n(aq ;a,b,c;q):={}_2\phi_1 (q^{-n}, abq^{n+1} ;\; cq; \; q,q )
=\frac{(c/abq^n;q)_n}{(cq;q)_n} (ab)^n\, q^{n(n+1)}.
$$
Therefore, since
$$
(c/abq^n;q)_n=(abq/c;q)_n(-c/ab)^n q^{-n(n+1)/2}\,,
$$
one obtains that
$$
P_n(aq ;a,b,c;q):=\frac{(abq/c;q)_n}{(cq;q)_n}(-c)^nq^{n(n+1)/2}.
\eqno (4.8)
$$
Likewise,
$$
P_n(cq ;a,b,c;q):=\frac{(bq;q)_n}{(aq;q)_n}(-a)^nq^{n(n+1)/2}.
$$
Hence, for the scalar product $\langle \psi_{aq}(x),
\psi_{aq}(x)\rangle$ we have the expression
$$
\sum_{n=0}^\infty\,\frac{(1-abq^{2n+1})\,(abq,aq,cq;q)_n}
{(1-abq)(abq/c,bq,q;q)_n\, (-ac)^n}\,q^{-n(n+3)/2}\,
P^2_n(aq;a,b,c;q) \qquad\qquad\qquad
$$  $$  \qquad\qquad\qquad
= \sum_{n=0}^\infty\,\frac{(1-abq^{2n+1})\,(abq/c,abq,aq;q)_n}
{(1-abq)(bq,cq,q;q)_n\,(-a/c)^n}\,q^{n(n-1)/2} =
\frac{(abq^2,c/a;q)_\infty}{(bq,cq;q)_\infty}, \eqno (4.9)
$$
where the relation (A.6) from Appendix has been used. Similarly,
for $\langle \psi_{cq}(x),\psi_{cq}(x) \rangle$ one has the
expression
$$
\sum_{n=0}^\infty\,\frac{(1-abq^{2n+1})(abq,aq,cq;q)_n}
{(1-abq)(-ac)^n\,(abq/c,bq,q;q)_n}\,q^{-n(n+3)/2}\,P^2_n(cq;a,b,c;q)
= \frac{(abq^2,a/c;q)_\infty}{(aq,abq/c;q)_\infty} , \eqno(4.10)
$$
where formula (A.7) from Appendix has been used. Thus, the values
$\lambda=aq$ and $\lambda=cq$ are spectral points of the operator
$I_2$.

Let us find other spectral points of $I_2$. Setting $\lambda = aq$
in (4.7), we see that the operator $J$ transforms $\psi _{aq}(x)$
into a linear combination of the functions $\psi_{aq^2}(x)$ and
$\psi _{aq}(x)$. We have to show that $\psi_{aq^2}(x)$ also
belongs to the Hilbert space ${\cal H}$, that is, that
$$
\langle \psi _{aq^2} ,\psi _{aq^2} \rangle = \sum_{n=0}^\infty\,
\frac{(abq,aq,cq;q)_n\,(1-abq^{2n+1})}{(abq/c,bq,q;q)_n\,
(1-abq)(-ac)^n} \,q^{-n(n+3)/2}\, P^2_n(aq^2;a,b,c;q)<\infty .
$$
In order to achieve this we note that since
$(aq^2;q)_k=(aq;q)_k(1-aq^{k+1})/(1-aq)$, we have
$$
P_n(aq^2;a,b,c;q)=\sum_{k=0}^n\frac{1-aq^{k+1}}{1-aq}
 \frac{(q^{-n};q)_k(abq^{n+1};q)_k(aq;q)_k}{(aq;q)_k(cq;q)_k}
\frac{q^k}{(q;q)_k}
$$   $$
\leq \frac1{1-aq} \sum_{k=0}^n
\frac{(q^{-n};q)_k(abq^{n+1};q)_k(aq;q)_k}{(aq;q)_k(cq;q)_k}
\frac{q^k}{(q;q)_k}=(1-aq)^{-1} P_n(aq;a,b,c;q).
$$
Therefore, the series for $\langle \psi _{aq^2} ,\psi _{aq^2}
\rangle$ is majorized (up to the finite constant $(1-aq)^{-1}$) by
the corresponding series for $\langle \psi _{aq} ,\psi _{aq}
\rangle$. Thus, $\psi _{aq^2}(x)$ is an eigenfunction of $I_2$ and
the point $aq^2$ belongs to the spectrum of the operator $I_2$.
Setting $\lambda = aq^2$ in (4.7) and acting similarly, one
obtains that $\psi _{aq^3}(x)$ is an eigenfunction of $I_2$ and
the point $aq^3$ belongs to the spectrum of $I_2$. Repeating this
procedure, one sees that $\psi _{aq^n}(x)$, $n=1,2,\cdots$, are
eigenfunctions of $I_2$ and the set $aq^n$, $n=1,2,\cdots$,
belongs to the spectrum of $I_2$. Likewise, one concludes that
$\psi _{cq^n}(x)$, $n=1,2,\cdots$, are eigenfunctions of $I_2$ and
the set $cq^n$, $n=1,2,\cdots$, belongs to the spectrum of $I_2$.
Note that so far we do not know whether the operator $I_2$ has
other spectral points or not. In order to solve this problem we
shall proceed as in subsection 3.2.

The functions $\psi _{aq^n}(x)$ and $\psi _{cq^n}(x)$,
$n=1,2,\cdots$, are linearly independent elements of the Hilbert
space ${\cal H}$. Suppose that $aq^n$ and $cq^n$, $n=1,2,\cdots$,
constitute the whole spectrum of the operator $I_2$. Then the set
of functions $\psi _{aq^n}(x)$ and $\psi_{cq^n}(x)$,
$n=1,2,\cdots$, is a basis of the space ${\cal H}$. Introducing
the notations $\Xi _n:=\xi_{aq^{n+1}}(x)$ and $\Xi'_n:=\xi
_{cq^{n+1}}(x)$, $n=0,1,2,\cdots$, we find from (4.7) that
$$
J\, \Xi _n = a^{-1}cq^{-2n-1}(1-aq^{n+1})(1-baq^{n+1}/c)\, \Xi
_{n+1} + d_n\,\Xi _n + a^{-1}cq^{-2n}(1-q^n)(1-aq^n/c)\, \Xi
_{n-1} ,
$$ $$
J \,\Xi'_n = c^{-1}aq^{-2n-1}(1-cq^{n+1})(1-bq^{n+1})\, \Xi _{n+1}
+ d'_n\,\Xi _n + c^{-1}aq^{-2n}(1-q^n)(1-cq^n/a)\, \Xi _{n-1} ,
$$
where
$$
d_n= \frac 1a [q^{-2n-1}c(1+q)-q^{-n}(ab+ac+a+c)] ,
$$  $$
d'_n= \frac 1c [q^{-2n-1}a(1+q)-q^{-n}(ab+ac+a+c)] .
$$

As we see, the matrix of the operator $J$ in the basis $\Xi _n
=\xi_{aq^{n+1}}(x)$, $\Xi'_n=\xi _{cq^{n+1}}(x)$,
$n=0,1,2,\cdots$, is not symmetric, although in the initial basis
$f_n$, $n=0,1,2,\cdots$, it was symmetric. The reason is that the
matrix $M:=((a_{mn})_{m,n=0}^\infty \ \;
(a'_{mn})_{m,n=0}^\infty)$ with entries
$$
a_{mn}:=\beta_m(aq^n),\ \ \  a'_{mn}:=\beta_m(cq^n),\ \ \
m,n=0,1,2,\cdots ,
$$
where $\beta_m(dq^n)$, $d=a,c$, are coefficients (4.4) in the
expansion $\psi _{dq^n}(x)=\sum _m \,\beta_m(dq^n)\,f_n(x)$ (see
above), is not unitary. This matrix $M$ is formed by adding the
columns of the matrix $(a'_{mn})$ to the columns of the matrix
$(a_{mn})$ from the right, that is,
$$
M =\left( \matrix{
 a_{11}&\cdots &a_{1k} &\cdots &a'_{11} &\cdots &a'_{1l}&\cdots\cr
 a_{21}&\cdots &a_{2k} &\cdots &a'_{21} &\cdots &a'_{2l}&\cdots\cr
 \cdots&\cdots &\cdots &\cdots &\cdots &\cdots  &\cdots &\cdots\cr
 a_{j1}&\cdots &a_{jk} &\cdots &a'_{j1} &\cdots &a'_{jl}&\cdots\cr
 \cdots&\cdots &\cdots &\cdots &\cdots &\cdots  &\cdots
 &\cdots\cr}
\right) .
$$
It maps the basis $\{ f_n\}$ into the basis $\{\psi_{aq^{n+1}},
\psi _{cq^{n+1}} \}$ in the Hilbert space ${\cal H}$. The
nonunitarity of the matrix $M$ is equivalent to the statement that
the basis $\Xi _n:=\xi _{aq^{n+1}}(x)$, $\Xi _n:=\xi
_{cq^{n+1}}(x)$, $n=0,1,2,\cdots$, is not normalized. In order to
normalize it we have to multiply $\Xi _n$ by appropriate numbers
$c_n$ and $\Xi'_n$ by numbers $c'_n$. Let $\hat\Xi _n = c_n\Xi
_n$, $\hat\Xi'_n =c'_n\Xi_n$, $n=0,1,2,\cdots$, be a normalized
basis. Then the operator $J$ is symmetric in this basis and has
the form
$$
J\,\hat\Xi _n = c_{n+1}^{-1}c_na^{-1}cq^{-2n-1}(1-aq^{n+1})
(1-abq^{n+1}/c)\,\hat\Xi_{n+1} + d_n\, \hat\Xi _n
$$  $$
+ c_{n-1}^{-1}c_n a^{-1}cq^{-2n}(1-aq^{n}/c)(1-q^n)\,\hat\Xi
_{n-1}, \eqno (4.11)
$$
$$
J\,\hat\Xi'_n = {c'}_{n+1}^{-1}{c'}_nc^{-1}aq^{-2n-1}(1-bq^{n+1})
(1-cq^{n+1})\,\hat\Xi _{n+1} + d'_n \,\hat\Xi _n
$$  $$
+{c'}_{n-1}^{-1}{c'}_n c^{-1}aq^{-2n}
(1-cq^{n}/a)(1-q^n)\,\hat\Xi_{n-1} , \eqno (4.12)
$$
The symmetricity of the matrix of the operator $J$ in the basis
$\{ \hat\Xi _n,\hat\Xi'_n \}$ means that
$$
c_{n+1}^{-1}c_n q^{-2n-1}(1-aq^{n+1})(1-abq^{n+1}/c)
=c_{n}^{-1}c_{n+1} q^{-2n-2} (1-aq^{n+1}/c)(1-q^{n+1}),
$$  $$
{c'}_{n+1}^{-1}{c'}_n q^{-2n-1}(1-bq^{n+1})(1-cq^{n+1})
={c'}_{n}^{-1}{c'}_{n+1} q^{-2n-2} (1-cq^{n+1}/a)(1-q^{n+1}).
$$
that is,
$$
\frac{c_{n}}{c_{n-1}} =\sqrt{q\frac{(1-aq^n)(1-abq^n/c)}
{(1-q^n)(1-aq^n/c)}} , \ \ \ \frac{c'_{n}}{c'_{n-1}}
=\sqrt{q\frac{(1-cq^n)(1-bq^n)} {(1-q^n)(1-cq^n/a)}} .
$$
Thus,
$$
c_n= C\left(
q^{n}\frac{(abq/c,aq;q)_n}{(aq/c,q;q)_n}\right)^{1/2},\ \ \ c'_n=
C'\left(q^{n}\frac{(bq,cq;q)_n}{(cq/a,q;q)_n}\right)^{1/2},
$$
where $C$ and $C'$ are some constants.

Therefore, in the expansions
$$
\hat\psi _{aq^n}(x)\equiv \hat\Xi _n(x)= \sum _m
\,c_n\,\beta_m(aq^n)\,f^l_m(x)=\sum _m \,{\hat a}_{mn}\, f^l_m(x),
\eqno (4.13)
$$
$$
\hat\psi _{cq^n}(x)\equiv \hat\Xi _n(x) = \sum _m
\,c'_n\,\beta_m(cq^n)\,f^l_m(x)=\sum _m \,{\hat a}'_{mn}\,
f^l_m(x), \eqno (4.14)
$$
the matrix ${\hat M}:=(({\hat a}_{mn})_{m,n=0}^\infty\ \; ({\hat
a}'_{mn})_{m,n=0}^\infty)$ with entries
$$
{\hat a}_{mn}= c_n\,\beta _m(aq^n) = C
\left(q^{n}\frac{(abq/c,aq;q)_n}{(aq/c,q;q)_n}
\frac{(abq,aq,cq;q)_m\,(1-abq^{2m+1})}{(abq/c,bq,q;q)_m\,
(1-abq)(-ac)^m }\right)^{1/2}
$$
$$
\times \, q^{-m(m+3)/4}\,P_m(aq^{n+1} ;a,b,c;q), \eqno (4.15)
$$  $$
{\hat a}'_{mn}= c_n\,\beta _m(cq^n)= C'
\left(q^{n}\frac{(bq,cq;q)_n}{((cq/a,q;q)_n}\,
\frac{(abq,aq,cq;q)_m\,(1-abq^{2m+1})}{(abq/c,bq,q;q)_m
\,(1-abq)(-ac)^m}\right)^{1/2} $$
$$ \times \,q^{-m(m+3)/4}\,P_m(cq^{n+1} ;a,b,c;q)\,,
\eqno(4.16)
$$
is unitary, provided that the constants $C$ and $C'$ are
appropriately chosen. In order to calculate these constants, one
can use the relations
$$
\sum _{m=0}^\infty |{\hat a}_{mn}|^2=1\,,\ \ \ \ \sum
_{m=0}^\infty |{\hat a}'_{mn}|^2=1 \,,
$$
for $n=0$. Then these sums are multiples of the sums in (4.9) and
(4.10), so we find that
$$
C=\frac{(bq,cq;q)^{1/2}_\infty}{(abq^2, c/a;q)^{1/2}_\infty} ,\ \
\ C'=\frac{(aq,abq/c;q)^{1/2}_\infty}{(abq^2,
a/c;q)^{1/2}_\infty}.    \eqno (4.17)
$$
The coefficients $c_n$ and $c'_n$ in (4.13)--(4.16) are thus real
and equal to
$$
c_n= \left( \frac{(bq,cq;q)_\infty}{(abq^2, c/a;q)_\infty}
\frac{(abq/c,aq;q)_n\, q^n}{(aq/c,q;q)_n} \right)^{1/2}, \ \ \ \
c'_n= \left( \frac{(aq,abq/c;q)_\infty}{(abq^2, a/c;q)_\infty}
\frac{(bq,cq;q)_n\, q^n}{(cq/a,q;q)_n} \right)^{1/2}.
$$

The orthogonality of the matrix ${\hat M}\equiv (({\hat
a}_{mn})_{m,n=0}^\infty\ \; ({\hat a}'_{mn})_{m,n=0}^\infty)$
means that
$$
\sum _m {\hat a}_{mn}{\hat a}_{mn'}=\delta_{nn'},\ \ \ \sum _m
{\hat a}'_{mn}{\hat a}'_{mn'}=\delta_{nn'},\ \ \ \sum _m {\hat
a}_{mn}{\hat a}'_{mn'}=0,       \eqno (4.18)
$$  $$
\sum _n ({\hat a}_{mn}{\hat a}_{m'n}+ {\hat a}'_{mn} {\hat
a}'_{m'n} ) =\delta_{mm'} .               \eqno (4.19)
$$
Substituting the expressions for ${\hat a}_{mn}$ and ${\hat
a}'_{mn}$ into (4.19), one obtains the relation
$$
\frac{(bq,
cq;q)_\infty}{(abq^2,c/a;q)_\infty}\,\sum_{n=0}^\infty\,
\frac{(aq,abq/c;q)_n
q^n}{(aq/c,q;q)_n}P_m(aq^{n+1})P_{m'}(aq^{n+1})
\qquad\qquad\qquad
$$  $$
+\frac{(aq,
abq/c;q)_\infty}{(abq^2,a/c;q)_\infty}\,\sum_{n=0}^\infty\,
\frac{(bq,cq;q)_n q^n}{(cq/a,q;q)_n}P_m(cq^{n+1})P_{m'}(cq^{n+1})
$$  $$ \qquad\qquad\qquad
=\frac{(1-abq)(bq,abq/c,q;q)_m}{(1-abq^{2m+1})(aq,abq,cq;q)_m}
 (-ac)^m q^{m(m+3)/2}\,\delta_{mm'}\, .        \eqno (4.20)
$$
This identity must give an orthogonality relation for the big
$q$-Jacobi polynomials $P_m(y)\equiv P_m(y;a,b,c;q)$. An only gap,
which appears here, is the following. We have assumed that the
points $aq^n$ and $cq^n$, $n=0,1,2,\cdots$, exhaust the whole
spectrum of the operator $I_2$. As in the case of the operator
$I_1$ in subsection 3.2, if the operator $I_2$ had other spectral
points $x_k$, then on the left-hand side of (4.20) would appear
other summands $\mu_{x_k}
P_m({x_k};a,b,c;q)P_{m'}({x_k};a,b,c;q)$, which correspond to
these additional points. Let us show that these additional
summands do not appear. We set $m=m'=0$ in the relation (4.20)
with the additional summands. This results in the equality
$$
\frac{(bq,
cq;q)_\infty}{(abq^2,c/a;q)_\infty}\,\sum_{n=0}^\infty\,
\frac{(aq,abq/c;q)_n q^n}{(aq/c,q;q)_n} \qquad\qquad\qquad\qquad
$$  $$   \qquad\qquad\qquad
+\frac{(aq, abq/c;q)_\infty}{(abq^2,a/c;q)_\infty}\,
\sum_{n=0}^\infty \,\frac{(bq,cq;q)_n q^n}{(cq/a,q;q)_n} +
\sum_k\,\mu_{x_k} =1.                           \eqno (4.21)
$$
In order to show that $\sum_k \mu_{x_k} = 0$, take into account
the relation
$$
\frac{(Aq/C, Bq/C;q)_\infty}{(q/C,ABq/C;q)_\infty}\, {}_2\phi_1
(A,B;C;\, q,q)  \qquad\qquad\qquad\qquad\qquad
$$   $$  \qquad\qquad\qquad
+ \, \frac{(A, B;q)_\infty}{(C/q,ABq/C;q)_\infty}\, {}_2\phi_1
(Aq/C,Bq/C;q^2/C;\, q,q) = 1
$$
(see formula (2.10.13) in [GR]). Putting here $A=aq$, $B=abq/c$ and
$C=aq/c$, we obtain relation (4.21) without the summand $\sum_k
\mu_{x_k}$. Therefore, in (4.21) the sum $\sum_k \mu_{x_k}$ does
really vanish and formula (4.20) gives an orthogonality relation
for big $q$-Jacobi polynomials.

By using the operators $I_2$ and $J$, we thus derived the
orthogonality relation for big $q$-Jacobi polynomials.

The orthogonality relation (4.20) enables one to formulate the
following statement: {\it The spectrum of the operator $I_2$
coincides with the set of points $aq^{n+1}$ and $cq^{n+1}$,
$n=0,1,2,\cdots$. The spectrum is simple and has one accumulation
point at 0.}
\bigskip

\noindent {\bf 4.3. Dual big $q$-Jacobi polynomials}
\bigskip

Now we consider the relations (4.18). They give the orthogonality
relation for the set of matrix elements ${\hat a}_{mn}$ and ${\hat
a}'_{mn}$, viewed as functions of $m$. Up to multiplicative
factors, they coincide with the functions
 $$
F_n(x;a,b,c;q):={}_3\phi_2 (x,abq/x, aq^{n+1};\; aq, cq;\; q,q),\
\ n=0,1,2,\cdots,  \eqno (4.22)
$$   $$
F'_n(x;a,b,c;q):={}_3\phi_2 (x,abq/x, cq^{n+1};\; aq, cq;\;
q,q)\equiv F_n(x;c,ab/c,a),\ \ n=0,1,2,\cdots ,  \eqno (4.23)
$$
considered on the corresponding sets of points. Namely, we have
$$
{\hat a}_{mn}\equiv {\hat
a}_{mn}(a,b,c)=C\,\left(q^{n}\frac{(abq/c,aq;q)_n}{(aq/c,q;q)_n}
\frac{(abq,aq,cq;q)_m\,(1-abq^{2m+1})}{(abq/c,bq,q;q)_m\,
(1-abq)(-ac)^m}\right)^{1/2}
$$  $$
\times q^{-m(m+3)/4}\,F_n(q^{-m};a,b,c;q),   \eqno (4.24)
$$   $$
{\hat a}'_{mn}\equiv {\hat
a}'_{mn}(a,b,c)=C'\,\left(q^{n}\frac{(bq,cq;q)_n}{(cq/a,q;q)_n}
\,\frac{(abq,aq,cq;q)_m\,(1-abq^{2m+1})}{(abq/c,bq,q;q)_m\,
(1-abq)(-ac)^m}\right)^{1/2}
$$  $$
\times q^{-m(m+3)/4}\,F'_n(q^{-m} ;a,b,c;q)\equiv {\hat
a}_{mn}(c,ab/c,a), \eqno (4.25)
$$
where $C$ and $C'$ are given by formulas (4.17). The relations
(4.18) lead to the following orthogonality relations for the
functions (4.22) and (4.23):
$$
\frac{(bq, cq;q)_\infty}{(abq^2,c/a;q)_\infty}\,\sum_{m=0}^\infty
\, \rho (m) F_n(q^{-m};a,b,c;q)\,F_{n'}(q^{-m};a,b,c;q)
=\frac{(aq/c,q;q)_n}{(aq,abq/c;q)_n q^n}\delta_{nn'}, \eqno (4.26)
$$   $$
\frac{(aq, abq/c;q)_\infty}{(abq^2,a/c;q)_\infty}\,
\sum_{m=0}^\infty\, \rho (m)
F'_n(q^{-m};a,b,c;q)\,F'_{n'}(q^{-m};a,b,c;q)
=\frac{(cq/a,q;q)_n}{(bq,cq;q)_n q^n}\,\delta_{nn'},  \eqno (4.27)
$$   $$
\sum_{m=0}^\infty\, \rho (m)
F_n(q^{-m};a,b,c;q)\,F'_{n'}(q^{-m};a,b,c;q)=0,  \eqno (4.28)
$$
where
$$
\rho (m):=\frac{(1-abq^{2m+1})(aq,abq,cq;q)_m}{(1-abq)
(bq,abq/c,q;q)_m\,(-ac)^m}\, q^{-m(m+3)/2}.
$$

There is another form for the functions $F_n(q^{-m};a,b,c;q)$ and
$F'_n(q^{-m};a,b,c;q)$. Indeed, one can use the relation (III.12)
of Appendix III in [GR] to obtain that
$$
F_n(q^{-m};a,b,c;q)=\frac{(cq^{-m}/ab;q)_m}{(cq;q)_m}
(abq^{m+1})^m \, {}_3\phi_2\left( \left. {q^{-m},abq^{m+1},q^{-n}
   \atop aq, abq/c} \right| q,aq^{n+1}/c  \right)
$$  $$
=\frac{(abq/c;q)_m}{(cq;q)_m} (-c)^{m}q^{m(m+1)/2} \,
{}_3\phi_2\left( \left. {q^{-m},abq^{m+1},q^{-n}
   \atop aq, abq/c} \right| q,aq^{n+1}/c  \right)
$$
and
$$
F'_n(q^{-m};a,b,c;q) =\frac{(bq;q)_m}{(aq;q)_m}
(-a)^{m}q^{m(m+1)/2}\, {}_3\phi_2\left( \left.
{q^{-m},abq^{m+1},q^{-n}
   \atop bq, cq} \right| q,cq^{n+1}/a  \right) .
$$
The basic hypergeometric functions ${}_3\phi_2$ in these formulas
are polynomials in $\mu (m):=q^{-m}+ab\,q^{m+1}$. So if we
introduce the notation
$$
D_n(\mu (m); a,b,c|q):=\left.
{}_3\phi_2\left({q^{-m},abq^{m+1},q^{-n}
  \atop aq, abq/c} \right| q,aq^{n+1}/c \right) ,       \eqno (4.29)
$$
then
$$
F_n(q^{-m};a,b,c;q)=\frac{(abq/c;q)_m}{(cq;q)_m}
(-c)^{m}q^{m(m+1)/2} D_n(\mu (m); a,b,c|q),
$$   $$
F'_n(q^{-m};a,b,c;q)=\frac{(bq;q)_m}{(aq;q)_m}
(-a)^{m}q^{m(m+1)/2} D_n(\mu (m); b,a,ab/c|q).
$$

Formula (4.26) directly leads to the orthogonality relation for
the polynomials $D_n(\mu(m))\equiv D_n(\mu (m); a,b,c|q)$:
$$
\sum_{m=0}^\infty
\frac{(1-abq^{2m+1})(aq,abq,abq/c;q)_m}{(1-abq)(bq,cq,q;q)_m}
\,(-c/a)^m\, q^{m(m-1)/2}\,D_n(\mu (m))\,D_{n'}(\mu (m))
$$   $$
=\frac{(abq^2,c/a;q)_\infty}{(bq,cq;q)_\infty} \frac{(aq/c,q;q)_n}
{(aq,abq/c;q)_nq^n} \delta_{nn'}.  \eqno (4.30)
$$
From (4.27) one obtains the orthogonality relation for the
polynomials $D_n(\mu (m);b,a,ab/c|q)$ (which follows also from the
relation (4.30) by interchanging $a$ and $b$ and replacing $c$ by
$ab/c$).

We call the polynomials $D_n(\mu (m);a,b,c|q)$ {\it dual big
$q$-Jacobi polynomials}. It is natural to ask whether they can be
identified with some known and thoroughly studied set of
polynomials. The answer is: they can be obtained from the
$q$-Racah polynomials $R_n(\mu (x);a,b,c,d|q)$ of Askey and Wilson
[AW] by setting $a=q^{-N-1}$ and sending $N\to \infty$, that is,
$$
D_n(\mu (x);a,b,c|q)=\lim_{N\to \infty} R_n( \mu (x); q^{-N-1},
a/c,a,b|q).   \eqno (4.31)
$$
Observe that the orthogonality relation (4.30) can be also derived
from formula (4.16) in [Ros]. But the derivation of this formula
(4.16) is rather complicated.

The dual polynomials (4.29) and the big $q$-Jacobi polynomials
(4.3) are interrelated in the following way:
$$
D_n(\mu(m);
a,b,c|q)=\frac{(cq;q)_m}{(abq/c;q)_m}(-c)^{-m}q^{-m(m+1)/2}
P_m(aq^{n+1};a,b,c|q).
$$

It is worth noting here that in the limit as $c\to 0$ the dual big
$q$-Jacobi polynomials $D_n(\mu (x);a,b,c|q)$ coincide with the
dual little $q$-Jacobi polynomials $d_n(\mu (x);b,a|q)$, defined
in section 3. The dual little $q$-Jacobi polynomials $d_n(\mu
(x);a,b|q)$ reduce, in turn, to the Al-Salam--Carlitz II
polynomials $V_n^{(a)}(s;q)$ on the $q$-linear lattice $s=q^{-x}$
(see [KSw], p. 114) in the case when the parameter $b$ vanishes,
that is,
$$
d_n(\mu (x);a,0|q)=´ç_2\phi_0(q^{-n},\, q^{-x}; - ;\; q,q^n/a)=
(-a)^{-n}q^{n(n-1)/2}\, V_n^{(a)}(q^{-x};q). \eqno (4.32)
$$
This means that we have now a complete chain of reductions
$$
R_n(\mu (x);a,b,c,d|q)\; {\mathop{\longrightarrow}_{a\to
\infty}}\; D_n(\mu(x);c,d,c/b|q)\; {\mathop{\longrightarrow}_{b\to \infty}}\;
 d_n(\mu(x);d,c|q)\; {\mathop{\longrightarrow}_{c= 0}}\;
 V_n^{(d)}(q^{-x};q)
$$
from the four-parameter family of $q$-Racah polynomials, which
occupy the upper level in the Askey-scheme of basic hypergeometric
polynomials (see [KSw], p. 62), down to the one-parameter set of
Al-Salam--Carlitz II polynomials from the second level in the same
scheme. So, the dual big and dual little $q$-Jacobi polynomials
$D_n(\mu(x);a,b,c|q)$ and $d_n(\mu(x);a,b|q)$ should occupy the
fourth and third level in the Askey-scheme, respectively.

The recurrence relations for the polynomials $D_n(\mu (m)\equiv
D_n(\mu (m);a,b,c|q)$ are obtained from the $q$-difference
equation (4.6). It has the form
$$
(q^{-m}-1)(1-abq^{m+1})D_n(\mu (m))=A_nD_{n+1}(\mu (m))
-(A_n+C_n)D_n(\mu (m))+C_nD_{n-1}(\mu (m)),
$$
where
$$
A_n=q^{-2n-1} (1-aq^{n+1})\left[(c/a)-bq^{n+1}\right],\ \ \ \
C_n=q^{-2n}(1-q^n)\left[(c/a)-q^n \right].
$$

The relation (4.28) leads to the equality (another proof of this
relation is given in Appendix)
$$
\sum_{m=0}^\infty\,(-1)^m
\frac{(1-abq^{2m+1})(abq;q)_m}{(1-abq)(q;q)_m}\,
q^{\frac{m(m-1)}{2}}\,D_n(\mu (m);a,b,c|q)\, D_{n'}(\mu (m);
b,a,ab/c|q)=0. \eqno (4.33)
$$

Note that from the expression (4.29) for the dual big $q$-Jacobi
polynomials it follows that they possess the symmetry property
$$
D_n(\mu (m); a,b,c|q)=D_n(\mu (m); ab/c,c,b|q). \eqno (4.34)
$$

The set of functions (4.22) and (4.23) form an orthogonal basis in
the Hilbert space ${\frak l}^2$ of functions, defined on the set
of points $m=0,1,2,\cdots$, with the scalar product
$$
\langle f_1,f_2\rangle = \sum_{m=0}^\infty \rho (m) \, f_1(m)\,
\overline{f_2(m)},
$$
where $\rho (m)$ is the same as in formulas (4.26)--(4.28).
Consequent from this fact, one can deduce (in the same way as in
the case of dual little $q$-Jacobi polynomials) that {\it the dual
big $q$-Jacobi polynomials $D_n(\mu (m); a,b,c|q)$ correspond to
indeterminate moment problem and the orthogonality measure for
them, given by formula (4.30), is not extremal}.

It is difficult to find extremal measures for these polynomials.

\newpage

\noindent {\bf 4.4. Generating functions for dual big $q$-Jacobi
polynomials}
\bigskip

Generating functions are known to be of great importance in the
theory of orthogonal polynomials (see, for example, [AAR] and [SM]).
For the sake of completeness, we briefly discuss in this section
some instances of linear generating functions for the dual
$q$-Jacobi polynomials $D_n(\mu(x);a,b,c|q)$ and $d_n(\mu(x);a,b|q)$.
To start with, let us consider a generating-function formula
$$
 \sum _{n=0}^\infty \frac{(aq;q)_n}{(q;q)_n}\, t^n D_n(\mu(x);a,b,c|q)
=\frac{(aqt;q)_\infty}{(t;q)_\infty}  \left. {}_2\phi_2 \left(
{q^{-x},\; abq^{x+1}
  \atop abq/c,\;  aqt} \right| q,\; aqt/c   \right) ,  \eqno (4.35)
$$
where $|t|<1$ and, as before, $\mu(x)=q^{-x}+abq^{x+1}$. To verify
(4.35), insert the explicit form (4.29) of the dual big $q$-Jacobi
polynomials
$$
D_n(\mu(x);a,b,c|q)=\sum_{k=0}^n
\frac{(q^{-x},abq^{x+1},q^{-n};q)_k}{(aq,abq/c,q;q)_k} \left(
\frac{aq^{n+1}}{c}\right) ^k
$$
into the left side of (4.35) and interchange the order of
summation. The subsequent use of the relations
$$
(a;q)_{m+k}=(a;q)_m(aq^m;q)_k=(a;q)_k(aq^k;q)_m,
$$  $$
(q^{-m-k};q)_k=(-1)^kq^{-mk-k(k+1)/2}(q^{m+1};q)_k
$$
(see [GR], Appendix I) simplifies the inner sum and enables one to
evaluate it by the $q$-binomial formula (3.16). This gives the
quotient of two infinite products in front of ${}_2\phi_2$ on the
right side of (4.35), times $(aqt;q)_k^{-1}$. The remaining sum
over $k$ yields ${}_2\phi_2$ series itself.

As a consistency check, one may also obtain (4.35) directly from
the generating function for the $q$-Racah polynomials
$R_n(\mu(x);\alpha,\beta,\gamma,\delta |q)$ (see formula (3.2.13)
in [KSw]) by setting $\alpha=q^{-N-1}$ and sending $N\to \infty$.
This results in the relation
$$
 \sum _{n=0}^\infty \frac{(aq;q)_n}{(q;q)_n}\, t^nD_n(\mu(x);a,b,c|q)
=\frac{(aq^{x+1}t;q)_\infty}{(t;q)_\infty}  \left. {}_2\phi_1
\left( {q^{-x},\; c^{-1}q^{-x}
  \atop abq/c} \right| q,\; atq^{x+1}   \right) .  \eqno (4.36)
$$
The left side of (4.36) depends on the variable $x$ by dint of the
combination $\mu(x)=q^{-x}+abq^{x+1}$. Off hand, it is not evident
that the right side of (4.36) is also a function of the lattice
$\mu(x)$. Nevertheless, this is the case. Moreover, the right
sides of (4.35) and (4.36) are equivalent: this fact is known in
the theory of special functions as Jackson's transformation
$$
{}_2\phi_1(a,b;\; c;\; q,z)=\frac{(az;q)_\infty}{(z;q)_\infty}\,
{}_2\phi_2(a,c/b;\; c,az;\; q,bz)
$$
(see, for example, [GR]).

The symmetry property (4.34) of the dual big $q$-Jacobi
polynomials $D_n(\mu(x);a,b,c|q)$, combined with (4.35), generates
another relation
$$
 \sum _{n=0}^\infty \frac{(abq/c;q)_n}{(q;q)_n} t^nD_n(\mu(x);a,b,c|q)
=\frac{(abqt/c;q)_\infty}{(t;q)_\infty}  \left. {}_2\phi_2 \left(
{q^{-x},\; abq^{x+1}
  \atop aq,\;  abqt/c} \right| q,\; aqt/c   \right)
$$  $$
=\frac{(abtq^{x+1}/c;q)_\infty}{(t;q)_\infty}  \left. {}_2\phi_1
\left( {q^{-x},\; b^{-1}q^{-x}
  \atop aq} \right| q,\; abtq^{x+1}/c   \right) .  \eqno (4.37)
$$

Similarly, a generating function for the dual little $q$-Jacobi
polynomials has the form
$$
 \sum _{n=0}^\infty \frac{(bq;q)_n}{(q;q)_n} (at)^nd_n(\mu(x);a,b|q)
=\frac{(tq^{-x},abtq^{x+1};q)_\infty}{(at,t;q)_\infty}. \eqno
(4.38)
$$
One can verify (4.38) directly by inserting the explicit form
(3.20) of $d_n(\mu(x);a,b|q)$ into the left  side of (4.38) and
repeating the same steps as in the case of deriving (4.35). This
will lead to the expression
$$
\frac{(abqt;q)_\infty}{(at;q)_\infty}  \, {}_2\phi_1 ( q^{-x},\;
abq^{x+1}; \; abqt; \; q, t )
$$
and it remains only to employ Heine's summation formula (1.5.1)
from [GR]. After a simple re-scaling of the parameters the
generating function (4.38) coincides with that, obtained earlier
in [BI].

The simplest way of obtaining (4.38) is to send $c\to 0$ in both
sides of (4.36): the ${}_2\phi_1$ series on the right side of
(4.36) reduces to ${}_1\phi_0 (q^{-x}; -;\; q,t/a)$, which is
evaluated by the $q$-binomial formula (3.16).

Finally, when the parameter $b$ vanishes, (4.38) reduces to the
known generating function
$$
\sum _{n=0}^\infty
\frac{(-1)^nq^{n(n-1)/2}}{(q;q)_n}t^nV^{(a)}_n(q^{-x};q)
=\frac{(tq^{-x};q)_\infty}{(at,t;q)_\infty}
$$
for the Al-Salam--Carlitz II polynomials (see (3.25.11) in [KSw]).
\bigskip

\noindent{\bf 5. DISCRETE $q$-ULTRASPHERICAL POLYNOMIALS AND THEIR
DUALS}
 \bigskip

\noindent{\bf 5.1. Discrete $q$-ultraspherical polynomials}
 \bigskip

For the big $q$-Jacobi polynomials $P_n(x ;a,b,c;q)$ the following
limit relation holds:
$$
\lim_{q\uparrow 1} P_n(x ;q^\alpha,q^\beta,-q^\gamma;q)
=\frac{P^{(\alpha,\beta)}_n(x)}{P^{(\alpha,\beta)}_n(1)} ,
$$
where $\gamma$ is real. Therefore, $\lim_{q\uparrow 1} P_n(x
;q^\alpha,q^\alpha,-q^\gamma;q)$ is a multiple of the Gegenbauer
(ultraspherical) polynomial $C_n^{(\alpha-1/2)}(x)$. For this
reason, we introduce the notation
$$
C_n^{(a^2)}(x;q):=P_n(x; a,a,-a;q)={}_3\phi_2
(q^{-n},a^2q^{n+1},x;\;
  aq, -aq;\;  q,q ) .    \eqno (5.1)
$$
It is obvious  from (5.1) that $C_n^{(a)}(x;q)$ is a rational
function in the parameter $a$.

{}From the recurrence relation for the big $q$-Jacobi polynomials
(see subsection 4.1) one readily verifies that the polynomials
(5.1) satisfy the following three-term recurrence relation:
$$
x\, C_n^{(a)}(x;q)=A_n(a)\, C_{n+1}^{(a)}(x;q) +C_n(a)\,
C_{n-1}^{(a)}(x;q), \eqno (5.2)
$$
where $A_n(a)=(1-aq^{n+1})/(1-aq^{2n+1})$, $C_n(a)=1-A_n(a)$, and
$C_0^{(a)}(x;q)\equiv 1$.

An orthogonality relation for $C_n^{(a)}(x;q)$, which follows from
that for the big $q$-Jacobi polynomials and is considered in the
next section, holds for positive values of $a$. We shall see that
the polynomials $C_n^{(a)}(x;q)$ are orthogonal also for imaginary
values of $a$ and $x$. In order to dispense with imaginary numbers
in this case, we introduce the following notation:
$$
\tilde C_n^{(a^2)}(x;q):= (-{\rm i})^n C_n^{(-a^2)}({\rm i}x;q)
=(-{\rm i})^n P_n({\rm i}x; {\rm i}a,{\rm i}a,-{\rm i}a;q),
 \eqno (5.3)
$$
where $x$ is real and $0<a<\infty$. The polynomials $\tilde
C_n^{(a)}(x;q)$ satisfy the recurrence relation
$$
x\tilde C_n^{(a)}(x;q)=\tilde A_n (a)\, \tilde C_{n+1}^{(a)}(x;q)
+\tilde C_n(a)\, \tilde C_{n-1}^{(a)}(x;q),  \eqno (5.4)
$$
where $\tilde A_n(a)=A_n(-a)=(1+aq^{n+1})/(1+aq^{2n+1})$, $\tilde
C_n(a)=\tilde A_n(a)-1$, and $\tilde C_0^{(a)}(x;q)\equiv 1$. Note
that $\tilde A_n(a) \ge 1$ and, hence, coefficients in the
recurrence relation (5.4) for $\tilde C_n^{(a)}(x;q)$ satisfy the
conditions $\tilde A_n(a)\tilde C_{n+1}(a)>0$ of Favard's
characterization theorem for $n=0,1,2,\cdots$ (see, for example,
[GR]). This means that these polynomials are orthogonal with
respect to a positive measure with infinitely many points of
support. An explicit form of this measure is derived in the next
section.

So, we have
$$
\tilde C_n^{(a)}(x;q)=(-{\rm i})^n C_n^{(-a)}({\rm i}x;q)=
(-{\rm i})^n {}_3\phi_2\left( \left. {q^{-n},-aq^{n+1},{\rm i}x
  \atop {\rm i}\sqrt{a}q, -{\rm i}\sqrt{a}q} \right| q,q \right) .
  \eqno (5.5)
$$
(Here and everywhere below under $\sqrt{a}$, $a>0$, we understand
a positive value of the root.) From the recurrence relation (5.4)
it follows that the polynomials (5.5) are real for $x\in {\Bbb R}$
and $0<a<\infty$. From (5.5) it is also obvious that they are
rational functions in the parameter $a$. Observe that the
situation when along with orthogonal polynomials $p_n(x)$,
depending on some parameters, the set of polynomials $(-{\rm i})^n
\,p_n({\rm i}x)$ is also orthogonal (but for other values of
parameters) is known; see, for example, [Rom], [Ask], and [CS].

We show below that the polynomials $C_n^{(a)}(x;q)$ and
$\tilde C_n^{(a)}(x;q)$, interrelated by (5.5), are orthogonal
with respect to discrete measures. For this reason, they may
be regarded [AK4] as a discrete version of $q$-ultraspherical
polynomials of Rogers (see, for example, [AI2]).
 \medskip

{\bf Proposition 5.1.} {\it The following expressions for the
discrete $q$-ultraspherical polynomials (5.1) hold:
$$
C_{2k}^{(a)}(x;q)=\frac{(q;q^2)_k\, a^{k}}{(aq^2;q^2)_k}
(-1)^kq^{k(k+1)}p_k(x^2/aq^2; q^{-1},a|q^2),  \eqno (5.6)
$$   $$
C_{2k+1}^{(a)}(x;q)=\frac{(q^3;q^2)_k\, a^{k}}{(aq^2;q^2)_k}
(-1)^kq^{k(k+1)}x\, p_k(x^2/aq^2; q,a|q^2),  \eqno (5.7)
$$
where $p_k(y;a,b|q)$ are the little $q$-Jacobi polynomials (3.3).}
\medskip

{\it Proof.} To start with (5.6), apply Singh's quadratic
transformation for a terminating ${}_3\phi_2$ series
$$
{}_3\phi_2\left( \left. {a^2,\ b^2,\ c
  \atop abq^{1/2}, -abq^{1/2}} \right| q,q \right) =
  {}_3\phi_2\left( \left. {a^{2},\ b^2,\ c^{2}
  \atop a^2b^2q,\ 0 } \right| q^2,q^2 \right) , \eqno (5.8)
$$
which is valid when both sides in (5.8) terminate (see [GR], formula
(3.10.13)), to the expression in (5.1) for $q$-ultraspherical
polynomials $C_{2k}^{(a)}(x;q)$. This results in the following:
$$
C_{2k}^{(a)}(x;q)=  {}_3\phi_2\left( q^{-2k},\ aq^{2k+1},\ x^{2};\
aq^2,\ 0 ;\ q^2,q^2 \right) .
$$
Now apply to this basic hypergeometric series ${}_3\phi_2$ the
transformation formula
$$
{}_2\phi_1\left( \left. {q^{-n},\ b
  \atop c } \right| q,z \right)=\frac{(c/b;q)_n}{(c;q)_n}\,
{}_3\phi_2\left( \left. {q^{-n},\ b,\  bzq^{-n}/c
  \atop bq^{1-n}/c,\ 0 } \right| q, q \right) \eqno (5.9)
$$
(see formula (III.7) from Appendix III in [GR]) in order to get
$$
C_{2k}^{(a)}(x;q)=\frac{(q;q^2)_k\, a^{k}}{(aq^2;q^2)_k}
(-1)^kq^{k(k+1)} {}_2\phi_1\left( q^{-2k},\ aq^{2k+1};\ q;\
q^2,x^2/a \right) .
$$
Comparing this formula with the expression for the little
$q$-Jacobi polynomials (3.3) one arrives at (5.6).

One can now prove (5.7) by induction with the aid of (5.9) and the
recurrence relation (5.2). Indeed, since $C_0^{(a)}(x;q)\equiv 1$
and $A_0(a)=1$, one obtains from (5.2) that $C_1^{(a)}(x;q)=x$. As
the next step use the fact that $C_2^{(a)}(x;q)={}_3\phi_2
(q^{-2},aq^3,x^2;\; aq^2,0;\; q^2,q^2)$ to evaluate from (5.2)
explicitly that
$$
C_{3}^{(a)}(x;q)=x\, {}_3\phi_2\left( q^{-2},\ aq^5,\ x^2;\ aq^2,\
0;\  q^2, q^2 \right) .
$$
So, let us suppose that
$$
C_{2k-1}^{(a)}(x;q)=x\, {}_3\phi_2\left( q^{-2(k-1)},\ aq^{2k+1},\
x^2;\  aq^2,\ 0;\  q^2, q^2 \right) \eqno (5.10)
$$
for $k=1,2,3,\cdots$, and evaluate a sum $A_{2k}^{-1}(a)x\,
C_{2k}^{(a)}(x;q)+(1-A_{2k}^{-1}(a))C_{2k-1}^{(a)}(x;q)$.  As
follows from the recurrence relation (5.2), this sum should be
equal to $C_{2k+1}^{(a)}(x;q)$. This is the case because it is
equal to
$$
x\, \left\{ A_{2k}^{-1}\, {}_3\phi_2\left( \left. {q^{-2k},\
aq^{2k+1},\ x^2 \atop aq^2,\ 0 } \right| q^2, q^2 \right)
 +(1-A_{2k}^{-1}){}_3\phi_2\left( \left. {q^{-2(k-1)},\
aq^{2k+1},\ x^2 \atop aq^2,\ 0 } \right| q^2, q^2 \right) \right\}
$$  $$
=x\, {}_3\phi_2\left( \left. {q^{-2k},\ aq^{2k+3},\ x^2 \atop
aq^2,\ 0 } \right| q^2, q^2 \right), \eqno (5.11)
$$
if one takes into account readily verified identities
$$
A_{2k}^{-1}(q^{-2k};q^2)_m+(1-A_{2k}^{-1})(q^{-2(k-1)};q^2)_m
=\frac{1-aq^{2(k+m)+1}}{1-aq^{2k+1}} (q^{-2k};q^2)_m ,
$$  $$
\frac{1-aq^{2(k+m)+1}}{1-aq^{2k+1}} (aq^{2k+1};q^2)_m=
(aq^{2k+3};q^2)_m ,
$$
for $m=0,1,2,\cdots ,k$. The right side of (5.11) does coincide with
$C_{2k+1}^{(a)}(x;q)$, defined by the same expression (5.10) with
$k\to k+1$. Thus, it remains only to apply the same transformation
formula (5.9) in order to arrive at (5.7). Proposition is proved.
 \medskip

{\it Remark.} Observe that in the process of proving formula
(5.7), we established a quadratic transformation
$$
{}_3\phi_2\left( \left. {q^{-2k-1},aq^{2k+2},x
  \atop \sqrt{a}q, -\sqrt{a}q} \right| q,q \right)=
x\,{}_3\phi_2\left( \left. {q^{-2k},aq^{2k+3},x^2
  \atop aq^2,\ 0}  \right| q^2,q^{2} \right)   \eqno (5.12)
$$
for the terminating basic hypergeometric polynomials ${}_3\phi_2$
with $k=0,1,2,\cdots$. The left side in (5.12) defines the
polynomials $C_{2k+1}^{(a)}(x;q)$ by (5.1), whereas the right side
follows from the expression (5.11) for the same polynomials. The
formula (5.12) represents an extension of Singh's quadratic
transformation (5.8) to the case when $a^2=q^{-2k-1}$ and,
therefore, the left side in (5.8) terminates, but the right side
does not.
 \medskip

It follows from (5.5)--(5.7) that
$$
\tilde C_{2k}^{(a)}(x;q)=\frac{(q;q^2)_k\, a^{k}}{(-aq^2;q^2)_k}
(-1)^kq^{k(k+1)}p_k(x^2/aq^2; q^{-1},-a|q^2),  \eqno (5.13)
$$   $$
\tilde C_{2k+1}^{(a)}(x;q)=\frac{(q^3;q^2)_k\,
a^{k}}{(-aq^2;q^2)_k} (-1)^kq^{k(k+1)}x\, p_k(x^2/aq^2; q,-a|q^2).
\eqno (5.14)
$$
In particular, it is clear from these formulas that the
polynomials $\tilde C_{n}^{(a)}(x;q)$ are real-valued for $x\in
{\Bbb R}$ and $a>0$.
 \bigskip

\noindent{\bf 5.2. Orthogonality relations for discrete
$q$-ultraspherical polynomials}
 \bigskip

Since the polynomials $C_n^{(a)}(x;q)$ are a particular case of
the big $q$-Jacobi polynomials, an orthogonality relation for them
follows from (4.20). Setting $a=b=-c$, $a>0$, into (4.20) and
considering the case when $m=2k$ and $m'=2k'$, one verifies that
two sums on the left of (4.20) coincide (since $ab/c=-a=c$) and we
obtain the following orthogonality relation for
$C_{2k}^{(a)}(x;q)$:
$$
\sum_{s=0}^\infty \frac{(aq^2;q^2)_s \, q^s}{(q^2;q^2)_s}\,
\,C_{2k}^{(a)}(\sqrt{a}q^{s+1};q)\,C_{2k'}^{(a)}(\sqrt{a}q^{s+1};q)
\qquad\qquad
$$ $$   \qquad\qquad
= \frac{(aq^3;q^2)_\infty} {(q;q^2)_\infty}\frac{(1{-}aq)
a^{2k}}{1{-}aq^{4k+1}} \frac{(q;q)_{2k}\,q^{k(2k+3)}}{(aq;q)_{2k}}\,
\delta_{kk'}, \eqno (5.15)
$$
where $\sqrt{a}$, $a>0$, denotes a positive value of the root.
Thus, {\it the family of polynomials $C_{2k}^{(a)}(x;q)$,
$k=0,1,2,\cdots$, with $0<a<q^{-2}$, is orthogonal on the set of
points $\sqrt{a}q^{s+1}$, $s=0,1,2,\cdots$.}

As we know, the polynomials $C_{2k}^{(a)}(x;q)$ are functions in
$x^2$, that is, $C_{2k}^{(a)}(\sqrt{a}q^{s+1};q)$ is in fact a
function in $aq^{2s+2}$. The set of functions $C_{2k}^{(a)}(x;q)$,
$k=0,1,2,\cdots$, constitutes a complete basis in the Hilbert
space ${\frak l}^2$ of functions $f(x^2)$ with the scalar product
$$
(f_1,f_2)=\sum_{s=0}^\infty \frac{(aq^2;q^2)_s\, q^s}{(q^2;q^2)_s}
f_1(aq^{2s+2})\overline{f_2(aq^{2s+2})} .
$$
This result can be obtained from the orthogonality relation for
the little $q$-Jacobi polynomials, if one takes into account
formula (5.6).

Putting $a=b=-c$, $a>0$, into (4.20) and considering the case when
$m=2k+1$ and $m'=2k'+1$, one verifies that two sums on the left of
(4.20) again coincide and we obtain the following orthogonality
relation for $C_{2k+1}^{(a)}(x;q)$:
$$
\sum_{s=0}^\infty \frac{(aq^2;q^2)_s\, q^s}{(q^2;q^2)_s}\,
C_{2k+1}^{(a)}(\sqrt{a}\, q^{s+1};q)\,C_{2k'+1}^{(a)}(\sqrt{a}\,
q^{s+1};q)  \qquad\qquad
 $$  $$  \qquad\qquad
 =\frac{(aq^3;q^2)_\infty}
{(q;q^2)_\infty}\frac{(1-aq)\, a^{2k+1}}{(1-aq^{4k+3})}
\frac{(q;q)_{2k+1}}{(aq;q)_{2k+1}}\,  q^{(k+2)(2k+1)}\,\delta_{kk'}.
\eqno (5.16)
$$
{\it The polynomials $C_{2k+1}^{(a)}(x;q)$, $k=0,1,2,\cdots$, with
$0<a<q^{-2}$, are thus orthogonal on the set of points
$\sqrt{a}\,q^{s+1}$, $s=0,1,2,\cdots$.}

The polynomials $x^{-1}C_{2k+1}^{(a)}(x;q)$ are functions in
$x^2$, that is, $x^{-1}C_{2k+1}^{(a)}(\sqrt{a}\,q^{s+1};q)$ are in
fact functions in $aq^{2s+2}$. The collection of functions
$C_{2k+1}^{(a)}(x;q)$, $k=0,1,2,\cdots$, constitute a complete
basis in the Hilbert space ${\frak l}^2$ of functions of the form
$F(x)=xf(x^2)$ with the scalar product
$$
(F_1,F_2)=\sum_{s=0}^\infty \frac{(aq^2;q^2)_s\,
q^s}{(q^2;q^2)_s}\,F_1(\sqrt{a}\, q^{s+1})\,\overline{F_2(\sqrt{a}\,
q^{s+1})} .
$$
Again, this result can be obtained from the orthogonality relation
for the little $q$-Jacobi polynomials if one takes into account
formula (5.7).

We have shown that the polynomials $C_{2k}^{(a)}(x;q)$,
$k=0,1,2,\cdots$, as well as the polynomials
$C_{2k+1}^{(a)}(x;q)$, $k=0,1,2,\cdots$, are orthogonal on the set
of points $\sqrt{a}\,q^{s+1}$, $s=0,1,2,\cdots$. However, the
polynomials $C_{2k}^{(a)}(x;q)$, $k=0,1,2,\cdots$, are not
orthogonal to the polynomials $C_{2k+1}^{(a)}(x;q)$,
$k=0,1,2,\cdots$, on this set of points. In order to obtain an
orthogonality for the whole collection of the polynomials
$C_{n}^{(a)}(x;q)$, $n=0,1,2,\cdots$, one has to consider them on
the set of points $\pm \sqrt{a}\, q^{s+1}$, $s=0,1,2,\cdots$.
Since the polynomials from the first set are even and the
polynomials from the second set are odd, for each $k,k'\in \{
0,1,2,\cdots\}$ the infinite sum
$$
I_1\equiv \sum_{s=0}^\infty \frac{(aq^2;q^2)_s\, q^s}{(q^2;q^2)_s}
\,C_{2k}^{(a)}(\sqrt{a}\,q^{s+1};q)\,C_{2k'+1}^{(a)}(\sqrt{a}\,q^{s+1};q)
$$
coincides with the following one
$$
I_2\equiv - \sum_{s=0}^\infty \frac{(aq^2;q^2)_s\,
q^s}{(q^2;q^2)_s}\, C_{2k}^{(a)}(-\sqrt{a}\,q^{s+1};q)\,
C_{2k'+1}^{(a)}(-\sqrt{a}\,q^{s+1};q).
$$
Therefore, $I_1-I_2=0$. This gives the orthogonality of
polynomials from the set $C_{2k}^{(a)}(x;q)$, $k=0,1,2,\cdots$,
with the polynomials from the set $C_{2k+1}^{(a)}(x;q)$,
$k=0,1,2,\cdots$. The orthogonality relation for the whole set of
polynomials $C_{n}^{(a)}(x;q)$, $n=0,1,2,\cdots $, can be written
in the form
$$
\sum_{s=0}^\infty \sum_{\varepsilon=\pm 1} \frac{(aq^2;q^2)_s\,
q^s}{(q^2;q^2)_s}\, C_{n}^{(a)}(\varepsilon\sqrt{a}\,q^{s+1};q)\,
C_{n'}^{(a)}(\varepsilon\sqrt{a}\,q^{s+1};q)
 $$   $$
=\frac{(aq^3;q^2)_\infty} {(q;q^2)_\infty}\frac{(1-aq)\,
a^{n}}{(1-aq^{2n+1})} \frac{(q;q)_{n}\,q^{n(n+3)/2}}{(aq;q)_{n}}
\, \delta_{nn'}. \eqno (5.17)
$$
We thus see that {\it the polynomials $C_{n}^{(a)}(x;q)$,
$n=0,1,2,\cdots $, with $0<a<q^{-2}$ are orthogonal on the set of
points $\pm \sqrt{a}\,q^{s+1}$, $s=0,1,2,\cdots$.}

An orthogonality relation for the polynomials $\tilde
C_{n}^{(a)}(x;q)$, $n=0,1,2,\cdots$, is derived by using the
relations (5.13), (5.14), and the orthogonality relation for the
little $q$-Jacobi polynomials. Writing down the orthogonality
relation (3.15) for the polynomials $p_k(x^2/aq^2; q^{-1},-a|q^2)$
and using the relation (5.13), one finds an orthogonality relation
for the set of polynomials $\tilde C_{2k}^{(a)}(x;q)$,
$k=0,1,2,\cdots$, with $a>0$. It has the form
$$
\sum_{s=0}^\infty \frac{(-aq^2;q^2)_s\, q^s}{(q^2;q^2)_s} \,
\tilde C_{2k}^{(a)}(\sqrt{a}\,q^{s+1};q)\, \tilde
C_{2k'}^{(a)}(\sqrt{a}\,q^{s+1};q)  \qquad\qquad
 $$  $$  \qquad\qquad
=\frac{(-aq^3;q^2)_\infty}
{(q;q^2)_\infty}\frac{(1+aq)\,a^{2k}}{(1+aq^{4k+1})}
\frac{(q;q)_{2k}}{(-aq;q)_{2k}}\, q^{k(2k+3)}\,\delta_{kk'}. \eqno(5.18)
$$
Consequently, {\it the family of polynomials $\tilde
C_{2k}^{(a)}(x;q)$, $k=0,1,2,\cdots$, is orthogonal on the set of
points $\sqrt{a}\, q^{s+1}$, $s=0,1,2,\cdots$.}

As in the case of polynomials $C_{2k}^{(a)}(x;q)$,
$k=0,1,2,\cdots$, the set $\tilde C_{2k}^{(a)}(x;q)$,
$k=0,1,2,\cdots$, is complete in the Hilbert space of functions
$f(x^2)$ with the corresponding scalar product.

Similarly, using formula (5.14) and the orthogonality relation for
the little $q$-Jacobi polynomials $p_k(x^2/aq^2; q,-a|q^2)$, we
find the orthogonality relation
$$
\sum_{s=0}^\infty \frac{(-aq^2;q^2)_s\, q^s}{(q^2;q^2)_s}\, \tilde
C_{2k+1}^{(a)}(\sqrt{a}\, q^{s+1};q)\, \tilde
C_{2k'+1}^{(a)}(\sqrt{a}\, q^{s+1};q)  \qquad\qquad
 $$  $$  \qquad\qquad
 =\frac{(-aq^3;q^2)_\infty}
{(q;q^2)_\infty}\frac{(1+aq)\, a^{2k+1}}{(1+aq^{4k+3})}
\frac{(q;q)_{2k+1}}{(-aq;q)_{2k+1}}\, q^{(k+2)(2k+1)}\,\delta_{kk'}
\eqno (5.19)
$$
for the set of polynomials $\tilde C_{2k+1}^{(a)}(x;q)$,
$k=0,1,2,\cdots$. We see from this relation that {\it for $a>0$
the polynomials $\tilde C_{2k+1}^{(a)}(x;q)$, $k=0,1,2,\cdots$,
are orthogonal on the same set of points $\sqrt{a}\,q^{s+1}$,
$s=0,1,2,\cdots$.}

Thus, the polynomials $\tilde C_{2k}^{(a)}(x;q)$,
$k=0,1,2,\cdots$, as well as the polynomials $\tilde
C_{2k+1}^{(a)}(x;q)$, $k=0,1,2,\cdots$, are orthogonal on the set
of points $\sqrt{a}\,q^{s+1}$, $s=0,1,2,\cdots$. However, the
polynomials $\tilde C_{2k}^{(a)}(x;q)$, $k=0,1,2,\cdots$, are not
orthogonal to the polynomials $\tilde C_{2k+1}^{(a)}(x;q)$,
$k=0,1,2,\cdots$, on this set of points. As in the previous case,
in order to prove that the polynomials $\tilde C_{2k}^{(a)}(x;q)$,
$k=0,1,2,\cdots$, are orthogonal to the polynomials $\tilde
C_{2k+1}^{(a)}(x;q)$, $k=0,1,2,\cdots$, one has to consider them
on the set of points $\pm \sqrt{a}\, q^{s+1}$, $s=0,1,2,\cdots$.
Since the polynomials from the first set are even and the
polynomials from the second set are odd, then the infinite sum
$$
I_1 \equiv \sum_{s=0}^\infty \frac{(-aq^2;q^2)_s\,
q^s}{(q^2;q^2)_s} \,\tilde C_{2k}^{(a)}(\sqrt{a}\,q^{s+1};q)\,
\tilde C_{2k'+1}^{(a)}(\sqrt{a}\,q^{s+1};q)
$$
coincides with the sum
$$
I_2\equiv - \sum_{s=0}^\infty \frac{(-aq^2;q^2)_s\,
q^s}{(q^2;q^2)_s} \,\tilde C_{2k}^{(a)}(-\sqrt{a}\,q^{s+1};q)\,
\tilde C_{2k'+1}^{(a)}(-\sqrt{a}\,q^{s+1};q).
$$
Consequently, $I_1-I_2=0$. This gives the mutual orthogonality of
the polynomials $\tilde C_{2k}^{(a)}(x;q)$, $k=0,1,2,\cdots$, to
the polynomials $\tilde C_{2k+1}^{(a)}(x;q)$, $k=0,1,2,\cdots$.
Thus, the orthogonality relation for the whole set of polynomials
$\tilde C_{n}^{(a)}(x;q)$, $n=0,1,2,\cdots $, can be written in
the form
$$
\sum_{s=0}^\infty \sum_{\varepsilon=\pm 1} \frac{(-aq^2;q^2)_s\,
q^s}{(q^2;q^2)_s} \,\tilde
C_{n}^{(a)}(\varepsilon\sqrt{a}\,q^{s+1};q)\, \tilde
C_{n'}^{(a)}(\varepsilon\sqrt{a}\,q^{s+1};q) \qquad\qquad\qquad
 $$   $$   \qquad\qquad\qquad
=\frac{(-aq^3;q^2)_\infty}
{(q;q^2)_\infty}\frac{(1+aq)\,a^{n}}{(1+aq^{2n+1})}
\frac{(q;q)_{n}}{(-aq;q)_{n}}\, q^{n(n+3)/2}\,\delta_{nn'}. \eqno(5.20)
$$
Note that the family of polynomials $\tilde C_{n}^{(a)}(x;q)$,
$n=0,1,2,\cdots $, corresponds to the determinate moment problem,
since the set of orthogonality is bounded. Thus, the orthogonality
measure in (5.20) is unique.

In fact, formula (5.20) extends the orthogonality relation for the
big $q$-Jacobi polynomials $P_n(x;a,a,-a;q)$ to a new domain of
values of the parameter $a$.
 \bigskip

\noindent{\bf 5.3. Dual discrete $q$-ultraspherical polynomials}
 \bigskip

The polynomials (4.29) are dual to the big $q$-Jacobi polynomials (4.3).
Let us set $a=b=-c$ in the polynomials (4.29), as we made
before in the polynomials (4.3). This gives the polynomials
$$
D_n^{(a^2)}(\mu (x;a^2)|q):=D_n(\mu (x;a^2); a,a,-a|q):=\left.
{}_3\phi_2\left({q^{-x},a^2q^{x+1},q^{-n}
  \atop aq, -aq} \right| q,-q^{n+1} \right) ,   \eqno (5.21)
$$
where $\mu(x;a^2)=q^{-x}+a^2q^{x+1}$. They satisfy the three-term
recurrence relation
$$
(q^{-x}+aq^{x+1}) D_n^{(a)}(\mu (x;a)|q)= -q^{-2n-1}(1-aq^{2n+2})
D_{n+1}^{(a)}(\mu (x;a)|q)   \qquad\qquad
 $$  $$   \qquad\qquad
+q^{-2n-1}(1+q)D_n^{(a)}(\mu (x;a)|q)
-q^{-2n}(1-q^{2n})D_{n-1}^{(a)}(\mu (x;a)|q),
$$
which follows from the recurrence relation for the polynomials
$D_n(\mu (x;ab); a,b,c|q)$ from section 4.3.

For the polynomials $D_n^{(a^2)}(\mu (x;a^2)|q)$ with imaginary
$a$ we introduce the notation
$$
\tilde D_n^{(a^2)}(\mu (x;-a^2)|q):=D_n(\mu (x;-a^2); {\rm
i}a,{\rm i}a,-{\rm i}a|q)$$
$$:=\left.{}_3\phi_2\left({q^{-x},-a^2q^{x+1},q^{-n}
  \atop {\rm i}aq, -{\rm i}aq} \right| q,-q^{n+1} \right)\,.\eqno (5.22)
$$
The polynomials $\tilde D_n^{(a)}(\mu (x;-a^2)|q)$ satisfy the
recurrence relation
$$
(q^{-x}-aq^{x+1}) \tilde D_n^{(a)}(\mu (x;-a)|q)=
-q^{-2n-1}(1+aq^{2n+2}) \tilde D_{n+1}^{(a)}(\mu (x;-a)|q)
\qquad\qquad
 $$  $$  \qquad\qquad
+q^{-2n-1}(1+q)\tilde D_n^{(a)}(\mu (x;-a)|q)
-q^{-2n}(1-q^{2n})\tilde D_{n-1}^{(a)}(\mu (x;-a)|q).
$$
It is obvious from this relation that the polynomials $\tilde
D_n^{(a)}(\mu (x;-a)|q)$ are real for $x\in {\Bbb R}$ and $a>0$.
For $a>0$ these polynomials satisfy the conditions of Favard's
theorem and, therefore, they are orthogonal.
 \medskip

{\bf Proposition 5.2.} {\it The following expressions for the dual
discrete $q$-ultraspherical polynomials (5.21) hold:
$$
D_n^{(a)}(\mu (2k;a)|q)=d_n(\mu (k;q^{-1}a);
q^{-1},a|q^2)\qquad\qquad
$$ $$  \qquad\qquad
={}_3\phi_1\left( \left. {q^{-2k},\ aq^{2k+1},\ q^{-2n}
  \atop aq^2 } \right| q^2, q^{2n+1} \right), \eqno (5.23)
$$  $$
D_n^{(a)}(\mu (2k+1;a)|q)=q^nd_n(\mu (k;qa); q,a|q^2) \qquad\qquad
$$ $$  \qquad\qquad
=q^n \, {}_3\phi_1\left( \left. {q^{-2k},\ aq^{2k+3},\  q^{-2n}
  \atop aq^2 } \right| q^2, q^{2n-1} \right), \eqno (5.24)
$$
where $k$ are nonnegative integers and $d_n(\mu (x;bc); b,c|q)$
are the dual little $q$-Jacobi polynomials (3.20).}
\medskip

{\it Proof.} Applying to the right side of (5.21) the formula
(III.13) from Appendix III in [GR] and then Singh's quadratic
relation (5.8) for terminating ${}_3\phi_2$ series, after some
transformations one obtains
$$
D_n^{(a^2)}(\mu (2k;a^2)|q)=a^{-2k}q^{-k(2k+1)}
 {}_3\phi_2\left( \left.
{q^{-2k},\ a^2q^{2k+1},\ a^2q^{2n+2}
  \atop a^2q^2,\ 0 } \right| q^2, q^{2} \right).
$$
Now apply the relation (0.6.26) from [KSw] in order to get
$$
D_n^{(a^2)}(\mu
(2k;a^2)|q)=\frac{(q^{-2k+1};q^2)_k}{(a^2q^2;q^2)_k}
 {}_2\phi_1\left( \left.
{q^{-2k},\ a^2q^{2k+1}
  \atop q } \right| q^2, q^{2n+2} \right).
$$
Using formula (III.8) of Appendix III from [GR], one arrives at
the expression for the polynomials $D_n^{(a^2)}(\mu (2k;a^2)|q)$
in terms of the basic hypergeometric function from (5.23),
coinciding with $d_n(\mu (k;q^{-1}a^2); q^{-1},a^2|q^2)$.

The formula (5.24) is proved in the same way by using the relation
(5.12). Proposition is proved.
 \medskip

For the polynomials $\tilde D_n^{(a)}(\mu (m;-a)|q)$ with
nonnegative integers $m$, we have the expressions
$$
\tilde D_n^{(a)}(\mu (2k;-a)|q)=d_n(\mu (k;-q^{-1}a);
q^{-1},-a|q^2) $$
$$={}_3\phi_1\left( \left. {q^{-2k},\ -aq^{2k+1},\
q^{-2n}
  \atop -aq^2 } \right| q^2, q^{2n+1} \right), \eqno (5.25)
$$  $$
\tilde D_n^{(a)}(\mu (2k{+}1;-a)|q)=q^nd_n(\mu (k;-qa); q,-a|q^2)$$
$$=q^n {}_3\phi_1\left( \left. {q^{-2k}, -aq^{2k+3}, q^{-2n}
  \atop -aq^2 } \right| q^2, q^{2n-1} \right). \eqno (5.26)
$$
It is plain from the explicit formulas that the polynomials $
D_n^{(a)}(\mu (m)|q)$ and $\tilde D_n^{(a)}(\mu (m)|q)$ are
rational functions of $a$.
 \bigskip

\noindent{\bf 5.4. Orthogonality relations for dual discrete
$q$-ultraspherical polynomials}
 \bigskip

An example of the orthogonality relation for $D_n^{(a^2)}(\mu
(x;a^2)|q)\equiv D_n( \mu (x;a^2);a,a,-a|q)$, $0<a<q^{-1}$, has
been discussed in section 4. However, these polynomials correspond
to the indeterminate moment problem and, therefore, this
orthogonality relation is not unique. Let us find other
orthogonality relations. In order to derive them we take into
account the relations (5.23) and (5.24), and the orthogonality
relation (3.21) for the dual little $q$-Jacobi polynomials. By
means of formula (5.23), we arrive at the following orthogonality
relation for $0<a<q^{-2}$:
$$
\sum_{k=0}^\infty
\frac{(1-aq^{4k+1})(aq;q)_{2k}}{(1-aq)(q;q)_{2k}}\, q^{k(2k-1)}
D_n^{(a)}(\mu (2k)|q) D_{n'}^{(a)}(\mu (2k)|q)
 =\frac{(aq^3;q^2)_\infty}{(q;q^2)_\infty}
 \frac{(q^2;q^2)_nq^{-n}}{(aq^2;q^2)_n}\,\delta_{nn'} ,
$$
where $\mu(2k)\equiv \mu(2k;a)$. The relation (5.24) leads to the
orthogonality, which can be written in the form
$$
\sum_{k=0}^\infty
\frac{(1-aq^{4k+3})(aq;q)_{2k+1}}{(1-aq)(q;q)_{2k+1}}\,
q^{k(2k+1)} D_n^{(a)}(\mu (2k+1)|q)\, D_{n'}^{(a)}(\mu (2k+1)|q)
 $$   $$
 =\frac{(aq^3;q^2)_\infty}{(q;q^2)_\infty}
 \frac{(q^2;q^2)_n\, q^{-n}}{(aq^2;q^2)_n} \delta_{nn'},
$$
where $\mu(2k+1)\equiv \mu(2k+1;a)$ and $0<a<q^{-2}$.

Thus, we have obtained two orthogonality relations for the
polynomials $ D_n^{(a)}(\mu (x;a)|q)$, $0<a<q^{-2}$, one on the
lattice $\mu (2k;a)\equiv q^{-2k}+aq^{2k+1}$, $k=0,1,2,\cdots$,
and another on the lattice $\mu (2k+1;a)\equiv
q^{-2k-1}+aq^{2k+3}$, $k=0,1,2,\cdots$. {\it The corresponding
orthogonality measures are extremal} since they are extremal for
the dual little $q$-Jacobi polynomials from formulas (5.23) and (5.24)
(see section 3).

The polynomials $\tilde D_n^{(a)}(\mu (x;a)|q)$ also correspond to
the indeterminate moment problem and, therefore, they have
infinitely many positive orthogonality measures. Some of their
orthogonality relations can be derived in the same manner as for
the polynomials $D_n^{(a)}(\mu (x)|q)$ by using the connection
(5.25) and (5.26) of these polynomials with the dual little
$q$-Jacobi polynomials (3.20). The relation (5.25) leads to the
orthogonality relation
$$
\sum_{k=0}^\infty
\frac{(1{+}aq^{4k+1})(-aq;q)_{2k}}{(1+aq)(q;q)_{2k}}\, q^{k(2k{-}1)}
\,\tilde D_n^{(a)}(\mu (2k)|q)\, \tilde D_{n'}^{(a)}(\mu (2k)|q)$$
 $$ =\frac{(-aq^3;q^2)_\infty}{(q;q^2)_\infty}
 \frac{(q^2;q^2)_n\,q^{-n}}{(-aq^2;q^2)_n} \,\delta_{nn'} ,
$$
where $\mu(2k)\equiv \mu(2k;-a)$, and the relation (5.26) gives rise
to the orthogonality relation, which can be written in the form
$$
\sum_{k=0}^\infty
\frac{(1+aq^{4k+3})(-aq;q)_{2k+1}}{(1+aq)(q;q)_{2k+1}} q^{k(2k+1)}
\tilde D_n^{(a)}(\mu (2k+1)|q)\tilde D_{n'}^{(a)}(\mu (2k+1)|q)
 $$  $$
 =\frac{(-aq^3;q^2)_\infty}{(q;q^2)_\infty}
 \frac{(q^2;q^2)_nq^{-n}}{(-aq^2;q^2)_n} \delta_{nn'},
$$
where $\mu(2k+1)\equiv \mu(2k+1;-a)$. In both cases, $a$ is any
positive number.

Thus, in the case of the polynomials $\tilde D_n^{(a)}(\mu
(x;-a)|q)$ we also have two orthogonality relations. {\it The
corresponding orthogonality measures are extremal} since they are
extremal for the dual little $q$-Jacobi polynomials from formulas
(5.25) and (5.26).

Note that the extremal measures for the polynomials $
D_n^{(a)}(\mu (x)|q)$ and $\tilde D_n^{(a)}(\mu (x)|q)$, discussed
in this section, can be used for constructing self-adjoint
extensions of the closed symmetric operators, connected with the
three-term recurrence relations for these polynomials and
representable in an appropriate basis by a Jacobi matrix (details
of such construction are given in [Ber], Chapter VII). These
operators are representation operators for discrete series
representations of the quantum algebra $U_q({\rm su}_{1,1})$ (see,
for example, [KS] for a description of this algebra). Moreover,
the parameter $a$ for these polynomials is connected with the
number $l$, which characterizes the corresponding representation
$T_l$ of the discrete series.
 \bigskip

\noindent{\bf 5.5. Other orthogonality relations}
 \bigskip

The polynomials $D_n^{(a)}(\mu (x;a)|q)$ and the polynomials
$\tilde D_n^{(a)}(\mu (x;-a)|q)$ correspond to the indeterminate
moment problems. For this reason, there exist infinitely many
orthogonality relations for them. Let us derive some set of these
relations for the polynomials $\tilde D_n^{(a)}(\mu (x;-a)|q)$, by
using orthogonality relations for the polynomials (5.18) in [BI].
These polynomials are (up to a factor) of the form
$$
u_n((e^\xi -e^{-\xi})/2; t_1,t_2|q) ={}_3\phi_1\left( \left.
{q\,e^{\xi}/t_1,\ -q\,e^{-\xi}/t_1,\ q^{-n}
  \atop -q^2/t_1t_2 } \right| q, q^{n}t_1/t_2 \right) \eqno (5.27)
$$
and orthogonality relations, parameterized by a number $d$, $q\le
d<1$, are given by the formula
$$
\sum_{n=-\infty}^\infty \frac{(-t_1q^{-n}/d,
t_1q^{n}d,-t_2q^{-n}/d, t_2q^{n}d;q)_\infty}{(-t_1t_2/q;q)_\infty}
 \frac{d^{4n}q^{n(2n-1)}(1+d^2q^{2n})}{(-d^2;q)_\infty
(-q/d^2;q)_\infty(q;q)_\infty}
 $$  $$
\times u_r\left( (d^{-1}q^{-n}-dq^n)/2; t_1,t_2\right) u_s\left(
(d^{-1}q^{-n}-dq^n)/2; t_1,t_2\right)
=\frac{(q;q)_r(t_1/t_2)^r}{(-q^2/t_1t_2;q)_r\, q^{r}} \delta_{rs}.
\eqno (5.28)
 $$
The orthogonality measure here is positive for $t_1,t_2\in {\Bbb
R}$ and $t_1t_2>0$. It is not known whether these measures are
extremal or not.

In order to use this orthogonality relation for the polynomials
$\tilde D_n^{(a)}(\mu (x;-a)|q)$, let us consider the
transformation formula
$$
{}_3\phi_2\left( \left. {q^{-2k},-a^2q^{2k+1},q^{-n}
  \atop {\rm i}aq, -{\rm i}aq} \right| q,-q^{n+1} \right)=
{}_3\phi_1\left( \left. {q^{-2k},-a^2q^{2k+1},q^{-2n}
  \atop -a^2q^2}  \right| q^2,q^{2n+1} \right), \eqno (5.29)
$$
which is true for any nonnegative integer values of $k$. It is
obtained by equating two expressions (5.22) and (5.25) for the
dual discrete $q$-ultraspherical polynomials $\tilde D_n^{(a)}(\mu
(2k;-a)|q)$. Observe that (5.29) is still valid if one replaces
numerator parameters $q^{-2k}$ and $-a^2q^{2k+1}$ in both sides of
it by $c^{-1}q^{-2k}$ and $-ca^2\, q^{2k+1}$, $c\in {\Bbb C}$,
respectively. Indeed, the left side of (5.29) represents a finite
sum:
$$
 {}_3\phi_2 (q^{-n},\alpha,\beta;\ \gamma,\delta;\ q,z):=\sum_{m=0}^n
\frac{(q^{-n},\alpha,\beta;q)_m}{(\gamma,\delta,q;q)_m} z^m.
\eqno (5.30)
$$
In the case in question $\alpha=q^{-2k}$ and $\beta=-a^2q^{2k+1}$,
so the $q$-shifted factorial $(\alpha,\beta;q)_m$ in (5.30) is equal
to
$$
(q^{-2k},-a^2q^{2k{+}1};q)_m
$$  $$
=\prod_{j=0}^{m-1}
[1-a^2q^{2j+1}-q^j (q^{-2k}-a^2q^{2k{+}1})]
=\prod_{j=0}^{m-1} [1-a^2q^{2j+1}-q^j \mu(2k;-a^2)] ,
\eqno (5.31)
$$
where, as before, $\mu(2k;-a^2)=q^{-2k}-a^2q^{2k+1}$. The left
side in (5.29) thus represents a polynomial $p_n(x)$ in the
$\mu(2k;-a^2)$ of degree $n$. In a similar manner, one easily
verifies that the right side of (5.29) also represents a polynomial
$p'_n(x)$ in the same $\mu(2k;-a^2)$ of degree $n$. In other
words, the transformation formula (5.29) states that the polynomials
$p_n(x)$ and $p'_n(x)$ are equal to each other on the infinite set
of distinct points $x_k=\mu(2k;-a^2)$. Thus, they are identical.

An immediate consequence of this statement is that (5.29) still
holds if one replaces the numerator parameters $q^{-2k}$ and
$a^2q^{2k+1}$ in both sides of (5.29) by $c^{-1}q^{-2k}$ and
$ca^2q^{2k+1}$, respectively. The point is that
 $$
(c^{-1}q^{-2k},-ca^2q^{2k+1};q)_m=\prod_{j=0}^{m-1}
[1-a^2q^{2j+1}-q^j (c^{-1}q^{-2k}-ca^2q^{2k+1})]
 $$  $$
=\prod_{j=0}^{m-1} [1-a^2q^{2j+1}-q^j \mu_c(2k;-a^2)] ,
$$
where $\mu_c(2k;-a^2)=c^{-1}q^{-2k}-ca^2q^{2k+1}$. So, the
replacements $q^{-2k}\to c^{-1}q^{-2k}$ and $a^2q^{2k+1}\to
ca^2q^{2k+1}$ change only the variable: $\mu(2k;-a^2)\to
\mu_c(2k;-a^2)$, whereas all other parameters in both sides of
(5.29) are unaltered. Thus, our statement is proved.

We are now in a position to establish other orthogonality
relations for the polynomials $\tilde D_n^{(a)}(\mu (x;-a)|q)$,
distinct from those, obtained in subsection 5.4. To achieve this, we
use the fact that the polynomials $\tilde D_n^{(a)}(\mu (x;-a)|q)$
at the points $x=x^{(d)}_k:= 2k - \ln(\sqrt{aq}/d)/\ln q $ are
equal to
$$
\tilde D_n^{(a)}(\mu (x^{(d)}_k;-a)|q)= {}_3\phi_2\left( \left.
{q^{-2k}d^{-1}\sqrt{aq},\ -q^{2k}d\sqrt{aq},\ q^{-n} \atop {\rm
i}\sqrt{a}\,q,\ -{\rm i}\sqrt{a}\,q } \right| q, -q^{n+1}
\right),\eqno (5.32)
$$
where $\mu(x^{(d)}_k;-a)= \sqrt{aq}\left( d^{-1}\,q^{-2k} - d
\,q^{2k}\right)$. From (5.27) and (5.29) (with $q^{-2k}$ and
$-a^2q^{2k+1}$ replaced by $d^{-1}q^{-2k}$ and $-da^2q^{2k+1}$,
respectively) it then follows that
$$
\tilde D_n^{(a)}(\mu(x^{(d)}_k;-a)|q)=u_n\left((d^{-1}\,q^{-2k}-
d\,q^{2k})/2;\sqrt{q^3/a},\sqrt{q/a}\,|q^2\right).
$$
Hence, from the orthogonality relations (5.28) one obtains infinite
number of orthogonality relations for the polynomials $\tilde
D_n^{(a)}(\mu(x;-a)|q)$, which are parameterized by the same $d$ as
in (5.28). They are of the form
$$
\sum_{n=-\infty}^\infty \frac{(-t_1q^{-2n}/d,t_1q^{2n}d,
-t_2q^{-2n}/d, t_2q^{2n}d;q^2)_\infty}{(-t_1t_2/q^2;q^2)_\infty}
 \frac{d^{4n}q^{2n(2n-1)}(1+d^2q^{4n})}{(-d^2;q^2)_\infty
(-q^2/d^2;q^2)_\infty(q^2;q^2)_\infty}
 $$  $$
\times \tilde D_r^{(a)}(\mu(x^{(d)}_n;-a)|q)\tilde
D_s^{(a)}(\mu(x^{(d)}_n;-a)|q)
=\frac{(q^2;q^2)_r}{(-q^2a;q^2)^2_r}\,\delta_{rs}, \eqno (5.33)
$$
where $t_1=\sqrt{q^3/a}$ and $t_2=\sqrt{q/a}$.

There exist yet another connection of the polynomials $\tilde
D_n^{(a)}(\mu(x)|q)$ with the polynomials (5.27). In order to obtain
it we consider the relation
$$
 {}_3\phi_2\left( \left. {q^{-2k-1},-a^2q^{2k+2},q^{-n}
  \atop {\rm i}aq, -{\rm i}aq} \right| q,-q^{n+1} \right)=
q^{n}\, {}_3\phi_1\left( \left. {q^{-2k},-a^2q^{2k+3},q^{-2n}
  \atop -a^2q^2}  \right| q^2,q^{2n-1} \right).  \eqno (5.34)
$$
This relation is true for nonnegative integer values of $k$.
However, it can be proved (in the same way as in the case of
formula (5.29)) that it holds also if we replace $q^{2k}$ and
$q^{-2k}$ by $cq^{2k}$ and $c^{-1}q^{-2k}$, respectively.

As in the previous case, we put $x=x^{(d)}_k:= 2k -
\ln(\sqrt{aq}/d)/\ln q $, that is, $\mu(x^{(d)}_k;-a)=
\sqrt{aq}\left( d^{-1}\,q^{-2k} - d \,q^{2k}\right)$. Then by
means of formula (5.34) (with $q^{-2k}$ and $q^{2k}$ replaced by
$d^{-1}q^{-2k}$ and $dq^{2k}$, respectively), we derive that
$$
\tilde D_n^{(a)}(\mu(x^{(d)}_k;-a)|q)=q^{n}
u_n((d^{-1}q^{-2k}-dq^{2k})/2;\sqrt{q/a},\sqrt{q^3/a}|q^2).
$$
We may apply to the polynomials
$u_n(x;\sqrt{q/a},\sqrt{q^3/a}|q^2)$ the orthogonality relations
(5.28) and obtain an infinite number of orthogonality relations for
the polynomials $\tilde D_n^{(a)}(\mu(y)|q)$. However, they
coincide with the orthogonality relations (5.33).

It is important to know whether an orthogonality measure for
polynomials is extremal. The extremality of the measures in (5.33)
for the polynomials $\tilde D_n^{(a)}(\mu_c (x;-a)|q)$ depends on
the extremality of the orthogonality measures in (5.28) for the
polynomials (5.27). If some of the measures in (5.28) are
extremal, then the corresponding measures in (5.33) are also
extremal.
 \bigskip

\noindent{\bf 6. DUALITY OF BIG $q$-LAGUERRE AND $q$-MEIXNER
POLYNOMIALS}
 \bigskip

\noindent{\bf 6.1. Operators related to big $q$-Laguerre polynomials}
\bigskip

Let ${\cal H}\equiv {\cal H}_a$ be the Hilbert space of functions,
introduced in subsection 3.1, which is fixed by a real number $a$,
such that $0<a<q^{-1}$. In this section we are interested in the
operator
$$
A\, f_n=r_{n+1} f_{n+1} + r_n f_{n-1} +d_n f_n, \eqno (6.1)
$$
where
$$
r_n=(-abq^{n+2})^{1/2}\sqrt{(1-q^{n})(1-aq^{n})(1-bq^n)}, \ \ \ \
d_n=-abq^{2n+1}(1+q)+ q^{n+1}(a+ab+b),
$$
$a$ is the same as above and $b$ is a fixed real number, such that
$b<0$.

Since $q<1$ the operator $A$ is bounded. Therefore, one can close
this operator and we assume in what follows that $A$ is a closed
(and consequently defined on the whole space ${\cal H}$) operator.
Since $A$ is symmetric, its closure is a self-adjoint operator. In
the same way as in subsection 3.1, one readily proves that $A$ is
a Hilbert--Schmidt operator. Therefore, $A$ has a discrete
spectrum.

We wish to find eigenfunctions $\xi_\lambda (x)$ of the operator
$A$, $A \xi_\lambda (x)=\lambda \xi_\lambda (x)$. We set
$$
\xi_\lambda (x)=\sum _{n=0}^\infty a_n(\lambda)f_n (x).
$$
Acting by $A$ upon both sides of this relation, one derives that
$$
\sum _{n=0}^\infty a_n(\lambda) [ r_{n+1}\,f_{n+1} +r_n\, f_{n-1}
+d_n\, f_n] =\lambda \sum _{n=0}^\infty a_n(\lambda)f_n.
$$
Comparing coefficients of a fixed $f_n$, one obtains a three-term
recurrence relation for the $a_n (\lambda)$:
$$
r_{n+1}a_{n+1}(\lambda) +r_n a_{n-1}(\lambda)+d_na_n
(\lambda)=\lambda a_n (\lambda).
$$
Making the substitution
$$
a_n (\lambda)=(-ab)^{-n/2}q^{-n(n+3)/4} \left( \frac{(aq,bq;q)_n}{(q;q)_n }
\right) ^{1/2} a'_n (\lambda)
$$
and using the explicit expressions for the $r_n$ and $d_n$, we
derive the relation
$$
(1-aq^{n+1})(1-bq^{n+1}) a'_{n+1}(\lambda)- abq^{n+1}(1-q^{n})
a'_{n-1}(\lambda) +d_na'_n(\lambda)=\lambda a'_n (\lambda),
$$
where, as before, $0<a<q^{-1}$ and $b<0$. It coincides with the
recurrence relation for the big $q$-Laguerre polynomials, which
are defined as
$$
P_n(\lambda ;a,b;q):= {}_3\phi_2 (q^{-n},0,\lambda;\; aq,bq;\;q,q)
\qquad\qquad
$$  $$   \qquad\qquad
=(q^{-n}/b;q)^{-1}_n\; {}_2\phi_1 (q^{-n},\, aq/\lambda ;\; aq;\
q,\lambda /b)   \eqno (6.2)
$$
(see formula (3.11.3) in [KSw]), that is, $
a'_n(\lambda)=P_n(\lambda;a,b;q)$. Therefore,
$$
a_n(\lambda)=(-ab)^{-n/2}q^{-n(n+3)/4} \left( \frac{(aq,bq;q)_n}
{(q;q)_n} \right) ^{1/2} P_n( \lambda ;a,b;q) .  \eqno (6.3)
$$

Thus, the eigenfunctions of the operator $A$ have the form
$$
\xi_\lambda(x)= \sum _{n=0}^\infty (-ab)^{-n/2}q^{-n(n+3)/4}\left(
\frac{(aq,bq;q)_n}{(q;q)_n} \right) ^{1/2} P_n( \lambda ;a,b;q)
f_n (x)  \eqno (6.4)
$$  $$
=\sum _{n=0}^\infty (-b)^{-n/2}a^{-3n/4}q^{-n(n+3)/4}
\frac{(aq;q)_n} {(q;q)_n} (bq;q)^{1/2}_n P_n( \lambda ;a,b;q)x^n.
$$
Since the spectrum of the operator $A$ is discrete, only a
countable set of these functions belongs to the Hilbert space
${\cal H}$. This discrete set of functions determines a spectrum
of $A$.

In order to be able to find a spectrum of the operator $A$, we
consider the linear operator $q^{-J_0}$ which is diagonal in the
basis $\{ f_n\}$ and is given by the formula
$$
q^{-J_0}f_n=(aq)^{-1/2}q^{-n} f_n
$$
(see formula (3.1)). We have to find how the operator $q^{-J_0}$
acts upon the eigenfunctions $\xi _\lambda(x)$ of the operator $A$
(which belong to the Hilbert space ${\cal H}$). In order to do
this one can use the $q$-difference equation
$$
q^{-n}(1-q^n)\lambda^2 \,p_n(\lambda)= B(\lambda)\,p_{n} (q\lambda
)-[B(\lambda)+D(\lambda)]p_n(\lambda)+D(\lambda)\,p_n(q^{-1}\lambda)
 \eqno(6.5)
$$
for the big $q$-Laguerre polynomials $p_n(\lambda)\equiv
P_n(\lambda; a,b;q)$, where
$$
B(\lambda)=abq(1-\lambda), \ \ \ \  D(\lambda)=(\lambda
-aq)(\lambda-bq).
$$
Multiply both sides of (6.5) by $k_nf_n(x)$ and sum up over $n$,
where $k_n$ are the coefficients of $P_n(\lambda; a,b;q)$ in the
expression (6.3) for the coefficients $a_n(\lambda)$. Taking into
account formula (6.4) and the form of the operator $q^{-J_0}$ in
the basis $\{ f_n\}$, one obtains the relation
$$
(aq)^{1/2}q^{-J_0} \lambda^2\xi_\lambda (x)=B(\lambda)
\xi_{q\lambda}(x)-[B(\lambda) +D(\lambda)-\lambda^2]\xi_\lambda
(x)+D(\lambda)\xi_{q^{-1}\lambda}(x),  \eqno (6.6)
$$
where $B(\lambda)$ and $D(\lambda)$ are the same as in (6.5).
 \bigskip

{\bf 6.2. Spectrum of $A$ and orthogonality of big $q$-Laguerre
polynomials}
 \medskip

The aim of this subsection is to find a basis in the Hilbert space
${\cal H}$, which consists of eigenfunctions of the operator $A$
in a normalized form, and to derive explicitly the unitary matrix
$U$, connecting this basis with the canonical basis $f_n$,
$n=0,1,2,\cdots$, in ${\cal H}$. This matrix directly leads to the
orthogonality relation for the big $q$-Laguerre polynomials.

Exactly as in section 4, one can show that for some value of
$\lambda$ (which belongs to the spectrum) the last term on the
right side of (6.6) has to vanish. There are two such values of
$\lambda$: $\lambda = aq$ and $\lambda = bq$, which are the roots
of the equation $D(\lambda)=0$. Let us show that both of these
points are spectral points of the operator $A$. Due to (6.2) we
have
$$
P_n(aq ;a,b;q)=(q^{-n}/b;q)_n^{-1}=\frac{(-bq)^n q^{n(n-1)/2}}{
(bq;q)_n},
$$
$$
P_n(bq ;a,b;q)=\frac{(-aq)^n q^{n(n-1)/2}}{ (aq;q)_n}.
$$
Hence, for the scalar product $\langle \xi_{aq}(x),
\xi_{aq}(x)\rangle$ we have the expression
$$
\sum_{n=0}^\infty\,\frac{(aq,bq;q)_n} {(-ab)^nq^{n(n+3)/2}\,
(q;q)_n}P^2_n(aq;a,b;q)
$$  $$
= \sum_{n=0}^\infty\,(-b/a )^n q^{n(n-1)/2}\, \frac{(aq;q)_n}
{(bq;q)_n(q;q)_n}={}_1\phi_1
(a;b;q,b/a)=\frac{(b/a;q)_\infty}{(b;q)_\infty} <\infty . \eqno
(6.7)
$$
Similarly, for $\langle \xi_{bq}(x),\xi_{bq}(x) \rangle$ one has
the expression
$$
\sum_{n=0}^\infty\,\frac{(aq,bq;q)_n} {(-ab)^nq^{n(n+3)/2}\,
(q;q)_n}P^2_n(bq;a,b;q) =\frac{(a/b;q)_\infty}{(a;q)_\infty}
<\infty .  \eqno (6.8)
$$
Thus, the values $\lambda=aq$ and $\lambda=bq$ are the spectral
points of the operator $A$.

Let us find other spectral points of the operator $A$. Setting
$\lambda = aq$ in (6.6), we see that the operator $q^{-J_0}$
transforms $\xi _{aq}(x)$ into a linear combination of the
functions $\xi_{aq^2}(x)$ and $\xi _{aq}(x)$. We have to show that
$\xi_{aq^2}(x)$ belongs to the Hilbert space ${\cal H}$, that is,
that
$$
\langle \xi _{aq^2} ,\xi _{aq^2} \rangle =
\sum_{n=0}^\infty\,\frac{(aq,bq;q)_n} {(-ab)^nq^{n(n+3)/2}\,
(q;q)_n}P^2_n(aq^2;a,b;q)<\infty .
$$
In order to achieve this we note that since
$(aq^2;q)_k=(aq;q)_k(1-aq^{k+1})/(1-aq)$, we have
$$
P_n(aq^2;a,b;q)=\sum_{k=0}^n\frac{1-aq^{k+1}}{1-aq}
 \frac{(q^{-n};q)_k(aq;q)_k}{(aq;q)_k(bq;q)_k}
\frac{q^k}{(q;q)_k}
$$   $$
\leq \frac1{1-aq} \sum_{k=0}^n
\frac{(q^{-n};q)_k(aq;q)_k}{(aq;q)_k(bq;q)_k}
\frac{q^k}{(q;q)_k}=(1-aq)^{-1} P_n(aq;a,b;q).
$$
Therefore, the series for $\langle \xi _{aq^2} ,\xi _{aq^2}
\rangle$ is majorized (up to the finite constant $(1-aq)^{-1}$) by
the corresponding series for $\langle \xi _{aq} ,\xi _{aq}
\rangle$. Thus, $\xi _{aq^2}(x)$ is an eigenfunction of $A$ and
the point $aq^2$ belongs to the spectrum of the operator $A$.
Setting $\lambda = aq^2$ in (6.6) and acting similarly, one
obtains that $\xi _{aq^3}(x)$ is an eigenfunction of $A$ and the
point $aq^3$ belongs to the spectrum of $A$. Repeating this
procedure, one sees that $\xi _{aq^n}(x)$, $n=1,2,\cdots$, are
eigenfunctions of $A$ and the set $aq^n$, $n=1,2,\cdots$, belongs
to the spectrum of $A$. Likewise, one concludes that $\xi
_{bq^n}(x)$, $n=1,2,\cdots$, are eigenfunctions of $A$ and the set
$bq^n$, $n=1,2,\cdots$, belongs to the spectrum of $A$. So far we
do not know whether the operator $A$ has other spectral points or
not. In order to solve this problem we shall proceed as in
subsection 3.2.

The functions $\xi _{aq^n}(x)$ and $\xi _{bq^n}(x)$,
$n=1,2,\cdots$, are linearly independent elements of the Hilbert
space ${\cal H}$. Suppose that $aq^n$ and $bq^n$, $n=1,2,\cdots$,
constitute the whole spectrum of the operator $A$. Then the set of
functions $\xi _{aq^n}(x)$ and $\xi_{bq^n}(x)$, $n=1,2,\cdots$, is
a basis in the space ${\cal H}$. Introducing the notations $\Xi
_n:=\xi_{aq^{n+1}}(x)$ and $\Xi'_n:=\xi _{bq^{n+1}}(x)$,
$n=0,1,2,\cdots$, and taking into account the relation
$B(\lambda)+D(\lambda)-\lambda^2 = abq(1+q)-\lambda q(ab+a+b)$, we
find from (6.6) that
$$
q^{-J_0}\Xi_n=a^{-3/2}bq^{-2n-3/2}(1-aq^{n+1})\Xi_{n+1}  -
a^{-3/2}q^{-2n-3/2} [b(1+q)-q^{n+1}(ab+a+b)]\Xi_n
$$  $$
+ \,a^{-3/2}bq^{-2n-1/2}(1-q^n)(1-aq^n/b)\Xi_{n-1}
$$
for $\lambda =aq^{n+1}$, that is, for $\xi_{aq^{n+1}} (x)=
\Xi_n(x)$. Similarly,
$$
q^{-J_0}\Xi'_n=a^{1/2}b^{-1}q^{-2n-3/2}(1-bq^{n+1})\Xi'_{n+1}-
a^{1/2}b^{-1}q^{-2n-3/2}
$$  $$
\times [1+q-a^{-1}q^{n+1}(ab+a+b)] \Xi'_n
+a^{1/2}b^{-1}q^{-2n-1/2}(1-q^n)(1-bq^{n}/a)\Xi'_{n-1}
$$
for $\lambda =bq^{n+1}$, that is, for $\xi_{bq^{n+1}} (x)=
\Xi'_n(x)$.

As we see, the matrix of the operator $q^{-J_0}$ in the basis $\Xi
_n =\xi_{aq^{n+1}}(x)$, $\Xi'_n=\xi _{bq^{n+1}}(x)$,
$n=0,1,2,\cdots$, is not symmetric, although in the initial basis
$f_n$, $n=0,1,2,\cdots$, it was symmetric. The reason is that the
matrix $M:=((b_{mn})_{m,n=0}^\infty\ \;
(b'_{mn})_{m,n=0}^\infty)$ with entries
$$
b_{mn}:=a_m(aq^n),\ \ \  b'_{mn}:=a_m(bq^n),\ \ \ m,n=0,1,2,\cdots
,
$$
where $a_m(dq^n)$, $d=a,b$, are coefficients (6.3) in the
expansion $\xi _{dq^n}(x)=\sum _m \,a_m(dq^n)\,f_n(x)$ (see
(6.4)), is not unitary. It maps the basis $\{ f_n\}$ into the
basis $\{\xi_{aq^{n+1}}, \xi _{bq^{n+1}} \}$ in the Hilbert space
${\cal H}\equiv {\cal H}_a$. The nonunitarity of the matrix $M$ is
equivalent to the statement that the basis $\Xi _n:=\xi
_{aq^{n+1}}(x)$, $\Xi' _n:=\xi _{bq^{n+1}}(x)$, $n=0,1,2,\cdots$,
is not normalized. In order to normalize it, we have to multiply
$\Xi _n$ by appropriate numbers $c_n$ and $\Xi'_n$ by numbers
$c'_n$. Let $\hat\Xi _n = c_n\Xi _n$, $\hat\Xi'_n =c'_n\Xi_n$,
$n=0,1,2,\cdots$, be a normalized basis. Then the operator
$q^{-J_0}$ is symmetric in this basis and has the form
$$
q^{-J_0}\hat\Xi_n=c_{n+1}^{-1}c_na^{-3/2}bq^{-2n-3/2}
(1-aq^{n+1})\hat\Xi_{n+1}
- a^{-3/2}q^{-2n-3/2} [b(1+q)-q^{n+1}(ab+a+b)]\hat\Xi_n
$$  $$
+
c_{n-1}^{-1}c_na^{-3/2}bq^{-2n-1/2}(1-q^n)(1-aq^n/b)\hat\Xi_{n-1}
,
$$
$$
q^{-J_0}\hat\Xi'_n={c'}_{n+1}^{-1}c'_na^{1/2}b^{-1}q^{-2n-3/2}
(1-bq^{n+1})\hat\Xi'_{n+1}- a^{1/2}b^{-1}q^{-2n-3/2}
$$  $$
\times [1+q-a^{-1}q^{n+1}(ab+a+b)]
\hat\Xi'_n+{c'}_{n-1}^{-1}c'_na^{1/2}
b^{-1}q^{-2n-1/2}(1-q^n)(1-bq^{n}/a)\hat\Xi'_{n-1}.
$$
The symmetricity of the matrix of the operator $q^{-J_0}$ in the
basis $\{ \hat\Xi _n,\hat\Xi'_n \}$ means that
$$
c_{n+1}^{-1}c_nq^{-2n-3/2}(1-aq^{n+1})=c_{n}^{-1}c_{n+1}q^{-2n-5/2}(1-q^{n+1})
(1-aq^{n+1}/b),
$$   $$
{c'}_{n+1}^{-1}c'_nq^{-2n-3/2}(1-bq^{n+1})=
{c'}_{n}^{-1}c'_{n+1}q^{-2n-5/2}(1-q^{n+1})(1-bq^{n+1}/a),
$$
that is,
$$
\frac{c_{n}}{c_{n-1}} =\sqrt{q\frac{(1-aq^n)} {(1-q^n)(1-aq^n/b)}}
, \ \ \ \frac{c'_{n}}{c'_{n-1}}
=\sqrt{q\frac{(1-bq^n)}{(1-q^n)(1-bq^n/a)}} .
$$
Thus,
$$
c_n=C\left( \frac{(aq;q)_n \, q^n} { (aq/b,q;q)_n }\right )^{1/2}
\eqno (6.9)
$$
and
$$
c'_n=C'\left( \frac{(bq;q)_n q^n }{(bq/a,q;q)_n }\right )^{\frac
12}, \eqno (6.10)
$$
where $C$ and $C'$ are some constants.

Therefore, in the expansions
$$
\hat\xi _{aq^n}(x)\equiv \hat\Xi _n(x)= \sum _m
\,c_n\,a_m(aq^n)\,f_m(x)=\sum _m \,{\hat b}_{mn}\, f_m(x), \eqno
(6.11)
$$
$$
\hat\xi' _{bq^n}(x)\equiv \hat\Xi _n(x) = \sum _m
\,c'_n\,a_m(bq^n)\,f_m(x)=\sum _m \,{\hat b}'_{mn}\, f_m(x), \eqno
(6.12)
$$
the matrix ${\hat M}:=(({\hat b}_{mn})_{m,n=0}^\infty\ \; ({\hat
b}'_{mn})_{m,n=0}^\infty)$ with entries
$$
{\hat b}_{mn}= c_n\,a_m(aq^n) = C
\left(\frac{(aq;q)_nq^n}{(aq/b,q;q)_n} \frac{(a,b;q)_m}{(q;q)_m\,
(-ab)^m }\right)^{\frac 12}  q^{-m(m+3)/4}\,P_m(aq^{n+1} ;a,b;q),
\eqno (6.13)
$$  $$
{\hat b}'_{mn}{=} c'_n\,a_m(bq^n){=} C'
\left(\frac{(bq;q)_nq^n}{((bq/a,q;q)_n}\, \frac{(a,b;q)_m}{(q;q)_m
\,(-ab)^m}\right)^{\frac 12} q^{-m(m+3)/4}\,P_m(bq^{n+1} ;a,b;q),
\eqno(6.14)
$$
is unitary, provided that the constants $C$ and $C'$ are
appropriately chosen. In order to calculate these constants, one
can use the relations
$$
\sum _{m=0}^\infty |{\hat b}_{mn}|^2=1\,,\ \ \ \ \sum
_{m=0}^\infty |{\hat b}'_{mn}|^2=1
$$
for $n=0$. Then these sums are multiples of the sums in (6.7) and
(6.8), so we find that
$$
C=\frac{(bq;q)^{1/2}_\infty}{(b/a;q)^{1/2}_\infty} ,\ \ \ \ \
C'=\frac{(aq;q)^{1/2}_\infty}{(a/b;q)^{1/2}_\infty}. \eqno (6.15)
$$
The coefficients $c_n$ and $c'_n$ in (6.11)--(6.14) are thus real
and equal to
$$
c_n=\left( \frac{(aq;q)_n(bq;q)_\infty \, q^n} { (aq/b,q;q)_n
(b/a;q)_\infty}\right )^{\frac 12} \eqno (6.16)
$$
and
$$
c'_n=\left( \frac{(bq;q)_n(aq;q)_\infty q^n} {(bq/a,q;q)_n
(a/b;q)_\infty }\right )^{\frac 12}. \eqno (6.17)
$$
The orthogonality of the matrix ${\hat M}\equiv ({\hat b}_{mn}\
{\hat b}'_{mn})$ means that
$$
\sum _m {\hat b}_{mn}{\hat b}_{mn'}=\delta_{nn'},\ \ \ \sum _m
{\hat b}'_{mn}{\hat b}'_{mn'}=\delta_{nn'},\ \ \ \sum _m {\hat
b}_{mn}{\hat b}'_{mn'}=0,       \eqno (6.18)
$$  $$
\sum _n ({\hat b}_{mn}{\hat b}_{m'n}+ {\hat b}'_{mn} {\hat
b}'_{m'n} ) =\delta_{mm'} .               \eqno (6.19)
$$
Substituting the expressions for ${\hat b}_{mn}$ and ${\hat
b}'_{mn}$ into (6.19), one obtains the relation
$$
\frac{(bq;q)_\infty}{(b/a;q)_\infty}  \sum _{n=0}^\infty
\frac{(aq;q)_n q^n}{(aq/b;q)_n(q;q)_n} P_m( aq^{n+1}
;a,b;q)P_{m'}( aq^{n+1} ;a,b;q)  \qquad\qquad\qquad
$$  $$
+\, \frac{(aq;q)_\infty}{(a/b;q)_\infty} \sum _{n=0}^\infty
\frac{(bq;q)_n q^n}{(bq/a;q)_n(q;q)_n} P_m( bq^{n+1}
;a,b;q)P_{m'}( bq^{n+1} ;a,b;q)
$$  $$ \qquad\qquad\qquad\qquad\qquad
 =\frac{(q;q)_m}
{(aq,bq;q)_m}(-ab)^{m}q^{m(m+3)/2}\delta_{mm'} . \eqno (6.20)
$$
This identity must give an orthogonality relation for the big
$q$-Laguerre polynomials $P_m(y)\equiv P_m(y;a,b;q)$. An only gap,
which appears here, is the following. We have assumed that the
points $aq^n$ and $bq^n$, $n=0,1,2,\cdots$, exhaust the whole
spectrum of the operator $A$. As in the case of the operator $I_2$
in section 4, if the operator $A$ had other spectral points $x_k$,
then on the left-hand side of (6.20) would appear other summands
$\mu_{x_k} P_m({x_k};a,b;q)P_{m'}({x_k};a,b;q)$, which correspond
to these additional points. Let us show that these additional
summands do not appear. To this end we set $m=m'=0$ in the
relation (6.20) with the additional summands. This results in the
equality
$$
\frac{(bq;q)_\infty}{(b/a;q)_\infty}  \sum _{n=0}^\infty
\frac{(aq;q)_n q^n}{(aq/b;q)_n(q;q)_n}+
\frac{(aq;q)_\infty}{(a/b;q)_\infty} \sum _{n=0}^\infty
\frac{(bq;q)_n q^n}{(bq/a;q)_n(q;q)_n}+ \sum_k\,\mu_{x_k}  =1  .
\eqno (6.21)
$$
In terms of the ${}_2\phi_1$ basic hypergeometric series this
identity can be written as
$$
\frac{(bq;q)_\infty}{(b/a;q)_\infty}\; {}_2\phi_1 (aq,0;\; aq/b;\;
q,q) +\frac{(aq;q)_\infty}{(a/b;q)_\infty}  {}_2\phi_1 (bq,0;\;
bq/a;\; q,q)+ \sum_k\,\mu_{x_k} =1. \eqno (6.22)
$$
It represents a particular case of Sears' three-term
transformation formula for ${}_2\phi_1(a,0;c;q,q)$ series (see
[GR], formula (3.3.5)) if $\mu_{x_k}=0$ for all values of $k$.
Therefore, in (6.20) the sum $\sum_k \mu_{x_k}$ does really vanish
and formula (6.20) gives an orthogonality relation for big
$q$-Laguerre polynomials.

By using the operators $A$ and $q^{-J_0}$, we thus derived the
orthogonality relation for big $q$-Laguerre polynomials.

The orthogonality relation (6.20) for big $q$-Laguerre polynomials
enables one to formulate the following statement: {\it The
spectrum of the operator $A$ coincides with the set of points
$aq^{n+1}$ and $bq^{n+1}$, $n=0,1,2,\cdots$. The spectrum is
simple and has one accumulation point at 0.}
\bigskip

\noindent{\bf 6.3. Dual polynomials and functions}
\bigskip

The matrix $M\equiv (\hat b_{mn}\ \hat b'_{mn})$ with entries
(6.13) and (6.14) is unitary and it connects two orthonormal bases
in the Hilbert space ${\cal H}$. The relations (6.18) for its
matrix elements is the orthogonality relation for the functions,
which are dual to the big $q$-Laguerre polynomials and are defined
as
$$
f_n(q^{-m}; a,b|q):= P_m(aq^{n+1};a,b;q), \ \ \ n=0,1,2,\cdots ,
\eqno (6.23)
$$  $$
g_n(q^{-m}; a,b|q):= P_m(bq^{n+1};a,b;q), \ \ \ n=0,1,2,\cdots .
\eqno (6.24)
$$
Taking into account the expressions for the entries $\hat b_{mn}$
and $\hat b'_{mn}$, the first two relations in (6.18) can be
written as
$$
\sum _{m=0}^\infty a_m(aq^{n+1})a_m({aq^{n'+1}})=c_n^{-2} \,\delta
_{nn'},\ \ \  \sum _{m=0}^\infty
a_m(bq^{n+1})a_m({bq^{n'+1}})={c'}_n^{-2}\, \delta _{nn'}.
$$
Substituting the explicit expressions for the coefficients
$a_m(\lambda_n)$, we derive the following orthogonality relations
for the functions (6.23) and (6.24):
$$
\sum _{m=0}^\infty \frac{(aq,bq;q)_m} {(q;q)_m(-abq^{2})^{m}}
\,q^{-m(m-1)/2} f_n(q^{-m}; a,b|q)f_{n'}(q^{-m}; a,b|q)
=c_n^{-2}\,\delta_{nn'}, \eqno (6.25)
$$  $$
\sum _{m=0}^\infty \frac{(aq,bq;q)_m}
{(q;q)_m(-abq^{2})^{m}}\,q^{-m(m-1)/2}\, g_n(q^{-m}; a,b|q)\,
g_{n'}(q^{-m}; a,b|q) ={c'}_n^{-2}\,\delta_{nn'}, \eqno (6.26)
$$  $$
\sum _{m=0}^\infty \frac{(aq,bq;q)_m}
{(q;q)_m(-abq^{2})^{m}}\,q^{-m(m-1)/2}\, f_n(q^{-m};
a,b|q)\, g_{n'}(q^{-m}; a,b|q)=0,
 \eqno (6.27)
$$
where $c_n$ and $c'_n$ are given by the formulas (6.16) and
(6.17).

Comparing the expression
$$
M_n(q^{-x};a,b;q):={}_2\phi_1 (q^{-n},q^{-x};\; aq;\ q,-q^{n+1}/b)
$$
for the $q$-Meixner polynomials with the explicit form (6.2) of
the big $q$-Laguerre polynomials $P_m (x;a,b;q)$, we see that
$$
f_n(q^{-m};a,b|q)=(q^{-m}/b;q)^{-1}_m M_n(q^{-m};a,-b/a;q).
$$
Since $(q^{-m}/b;q)_m=(bq;q)_m(-b)^{-m}q^{-m(m+1)/2}$, the
orthogonality relation (6.25) leads to the orthogonality relation
for the $q$-Meixner polynomials $M_n(q^{-m})\equiv M_n(q^{-m};
a,-b/a;q)$:
$$
\sum _{m=0}^\infty \frac{(aq;q)_m (-b/a)^m \,q^{m(m-1)/2}}
{(bq,q;q)_m} \,M_n(q^{-m})\,M_{n'}(q^{-m})
$$
$$
=\frac{(b/a;q)_\infty}{(bq;q)_\infty}
\frac{(aq/b,q;q)_n}{(aq;q)_n} \,q^{-n} \,\delta_{nn'} ,
 \eqno (6.28)
$$
where, as before, $0<a<q^{-1}$ and $b<0$. This orthogonality
relation coincides with the known formula for the $q$-Meixner
polynomials (see, for example, (3.13.2) in [KSw]).

The functions (6.24) are also expressed in terms of $q$-Meixner
polynomials. Indeed, we have
$$
g_n(q^{-m};a,b|q)={}_3\phi_2 (q^{-m},0,bq^{n+1};\; aq,bq;\ q,q)
$$  $$
=(q^{-m}/a;q)_m^{-1}\, {}_2\phi_1 (q^{-n},q^{-m};\; bq;\
q,bq^{n+1}/a) =(q^{-m}/a;q)_m^{-1}\, M_n(q^{-m}; b,-a/b;\; q) ,
$$
where $b<0$, that is, one of the parameters in these $q$-Meixner
polynomials is negative.

Substituting this expression for $g_n(q^{-m};a,b|q)$ into (6.26),
we obtain the orthogonality relation for $q$-Meixner polynomials
$M_n(q^{-m})\equiv M_n(q^{-m}; b,-a/b;\; q)$ with negative $b$:
$$
\sum _{m=0}^\infty \frac{(bq;q)_m (-a/b)^m}{(aq,q;q)_m}
\, q^{m(m-1)/2}\,M_n(q^{-m})\,M_{n'}(q^{-m})
$$
$$
=\frac{(a/b;q)_\infty}{(aq;q)_\infty}
\frac{(bq/a,q;q)_n}{(bq;q)_n} \,q^{-n}\,\delta_{nn'}. \eqno (6.29)
$$
Observe that this orthogonality relation is of the same form as
for $b>0$. As far as we know, this type of orthogonality relation
for negative values of the parameter $b$ has been first discussed
in [AAK].

The relation (6.27) can be written as the equality
$$
\sum_{m=0}^\infty \frac{(-1)^m q^{m(m-1)/2}}{(q;q)_m}
 M_n(q^{-m}; a,-b/a;\; q)M_{n'}(q^{-m}; b,-a/b;\; q)=0,
 $$
which holds for $n,n'=0,1,2,\cdots$. The validity of this identity
for arbitrary nonnegative integers $n$ and $n'$ can be verified
directly by using Jackson's $q$-exponential function
$$
E_q(z):= \sum_{n=0}^\infty \frac{q^{n(n-1)/2}}{(q;q)_n}\,z^n=
(-z;q)_\infty
$$
and the fact that $E_q(z)$ has zeroes at the points $z_j =
-q^{-j}$, $j=0,1,2,\cdots$.

Notice that the appearance of the $q$-Meixner polynomials
here as a dual family with respect to the big $q$-Laguerre
polynomials is quite natural because the transformation
$q\to q^{-1}$ interrelates these two sets of polynomials,
that is,
$$
M_n(x;b,c;q^{-1})=(q^{-n}/b;q)_n\, P_n(qx/b;1/b,-c;q).
$$

Let us introduce the Hilbert space ${\frak l}_{a,b}^2$ of
functions $F(q^{-m})$ on the set $m\in \{0,1,2,\cdots \}$ with a
scalar product given by the formula
$$
\langle F_1,F_2\rangle_b = \sum _{m=0}^\infty \rho(m)
\,F_1(q^{-m})\,\overline{F_2(q^{-m})}, \eqno (6.30)
$$
where the weight function
$$
\rho (m)=\frac{(aq,bq;q)_m} {(q;q)_m(-abq^2)^{m}}\, q^{-m(m-1)/2}
$$
is the same as in (6.25)--(6.27). Now we can formulate the
following statement.
\medskip

{\bf Theorem 6.1.} {\it The functions (6.23) and (6.24) constitute
an orthogonal basis in the Hilbert space ${\frak l}_{a,b}^2$.}
\medskip

{\it Proof.} To show that the system of functions (6.23) and
(6.24) constitutes a complete basis in the space ${\frak
l}_{a,b}^2$ we take in ${\frak l}_{a,b}^2$ the set of functions
$F_k$, $k=0,1,2, \cdots$, such that $F_k(q^{-m})=\delta _{km}$. It
is clear that these functions constitute a basis in the space
${\frak l}_{a,b}^2$. Let us show that each of these functions
$F_k$ belongs to the closure $\bar V$ of the linear span $V$ of
the functions (6.23) and (6.24). This will prove the theorem. We
consider the functions
$$
{\hat F}_k(q^{-m})=\sum _{n=0} ^\infty \hat b_{kn}\hat b_{mn}+\sum
_{n=0} ^\infty \,\hat b'_{kn}\,\hat b'_{mn} ,\ \ \ \ k=0,1,2,\cdots ,
$$
where $\hat b_{jn}$ and $\hat b'_{kn}$ are the same as in (6.19).
Then $\rho^{-1}(m){\hat F}_k(q^{-m})$ is an infinite linear
combination of the functions (6.23) and (6.24). Moreover, $\hat
F_k(q^{-m})\in \bar V$ and, due to (6.19), $\hat F_k$,
$k=0,1,2,\cdots$, coincide, up to a constant, with the
corresponding functions $F_k$, introduced above. The theorem is
proved.
\medskip

The weight function $\rho(m)$ in (6.30) does not coincide with the
orthogonality measure for $q$-Meixner polynomials. Multiplying
this weight function by $[(bq;q)_m(-b)^{-m}q^{-m(m+1)/2}]^{-2}$,
we obtain the measure in (6.28). Let ${\frak l}^2_{(1)}$ be the
Hilbert space of functions $F(q^{-m})$ on the set $m\in
\{0,1,2,\cdots \}$ with the scalar product
$$
\langle F_1,F_2\rangle_{(1)} = \sum _{m=0}^\infty \frac{(aq;q)_m
(-b/a)^m}{(bq,q;q)_m}\, q^{m(m-1)/2}
\,F_1(q^{-m})\,\overline{F_2(q^{-m})} ,
$$
where the weight function coincides with the measure in (6.28).

Taking into account the modification of the measure and the
statement of Theorem 6.1, we conclude that the $q$-Meixner
polynomials $M_n(q^{-m};\, a,-b/a;\; q)$ and the functions
$$
(bq;q)_m\,(-b)^{-m}\,q^{-m(m+1)/2}\, g_n(q^{-m}; a,b|q)
$$
constitute an orthogonal basis in the space ${\frak l}^2_{(1)}$.
\medskip

{\bf Proposition 6.1.} {\it The $q$-Meixner polynomials
$M_n(q^{-m};a,c;q)$, $n=0,1,2,\cdots$, with the parameters $a$ and
$c=-b/a$ do not constitute a complete basis in the Hilbert space
${\frak l}^2_{(1)}$, that is, the $q$-Meixner polynomials are
associated with the indeterminate moment problem and the measure
in (6.28) is not an extremal measure for these polynomials.}
\medskip

{\it Proof.} In order to prove this proposition we note that if
the $q$-Meixner polynomials were associated with the determinate
moment problem, then they would constitute a basis in the space of
square integrable functions with respect to the measure from
(6.28). However, this is not the case. By the definition of an
extremal measure, if the measure in (6.28) were extremal, then
again the set of the $q$-Meixner polynomials would be a basis in
that space. Therefore, the measure is not extremal. Proposition is
proved.
\medskip

Let now ${\frak l}^2_{(2)}$ be the Hilbert space of functions
$F(q^{-m})$ on the set $m\in \{0,1,2,\cdots \}$, with the scalar
product
$$
\langle F_1,F_2\rangle_{(2)} = \sum _{m=0}^\infty \frac{(bq;q)_m
(-a/b)^m}{(aq,q;q)_m}\, q^{m(m-1)/2}
\,F_1(q^{-m})\,\overline{F_2(q^{-m})}.
$$
The measure here coincides with the orthogonality measure in
(6.29) for $q$-Meixner polynomials $M_n(q^{-m}; b,-a/b;\; q)$,
$b<0$. The following proposition is proved in the same way as
Proposition 6.1.
\medskip

{\bf Proposition 6.2.} {\it The $q$-Meixner polynomials
$M_n(q^{-m};b,-a/b;q)$, $n=0,1,2,\cdots$, with $b<0$ do not
constitute a complete basis in the Hilbert space ${\frak
l}^2_{(2)}$, that is, these $q$-Meixner polynomials are associated
with the indeterminate moment problem and the measure in (6.29) is
not an extremal measure for them.}
 \bigskip

\noindent{\bf 6.4. Generating function for big $q$-Laguerre polynomials}
\bigskip

The aim of this section is to derive a generating function for the
big $q$-Laguerre polynomials
$$
G(x,t; a,b;\, q):= \sum _{n=0}^\infty
\frac{(aq,bq;q)_n}{(q;q)_n}q^{-n(n-1)/2} P_n(x;a,b;q)t^n , \eqno
(6.31)
$$
which will be used in the next subsection. Observe that this
formula is a bit more general than each of the three instances of
generating functions for big $q$-Laguerre polynomials, given in
section 3.11 of [KSw].

Employing the explicit expression
$$
P_n(x;a,b;q)=(b^{-1}q^{-n};q)_n^{-1}{}_2\phi_1(q^{-n},aqx^{-1};\;
aq;\ q, x/q)
$$
for big $q$-Laguerre polynomials, one obtains
$$
G(x,t; a,b;\, q)=\sum_{n=0}^\infty \frac{(aq;q)_n}{(q;q)_n}
(-bqt)^n \sum_ {k=0}^n \frac{(q^{-n},aqx^{-1};q)_k}{(aq,q;q)_k}
\left( \frac xb \right) ^k
$$
$$
=\sum_{n=0}^\infty (aq;q)_n
(-bqt)^n \sum_ {k=0}^n \frac{(-x/b)^k(aqx^{-1};q)_k}
{(aq,q;q)_k(q;q)_{n-k}} q^{-nk+k(k-1)/2}
$$  $$
=\sum_{k=0}^\infty \frac{(aqx^{-1};q)_k(-x/b)^k}
{(aq,q;q)_k} q^{k(k-1)/2}
 \sum_ {m=0}^\infty \frac{(aq;q)_{m+k}}
{(q;q)_m} (-bqt)^{m+k} q^{-(k+m)k}
$$   $$
=\sum_{k=0}^\infty \frac{(aqx^{-1};q)_k}
{(q;q)_k} (xt)^k q^{-k(k-1)/2}
 \sum_ {m=0}^\infty \frac{(aq^{k+1};q)_m}
{(q;q)_m} (-bq^{1-k}t)^m .
$$
By the $q$-binomial theorem, the last sum equals to
$(-abq^2;q)_\infty/(-bq^{1-k}t;q)_\infty$.
Since
$$
(-bq^{1-k}t;q)_\infty =q^{-k(k-1)/2}(-q/bqt;q)_k(-bqt;q)_\infty ,
$$
then
$$
\frac{(-abq^2;q)_\infty}{(-bq^{1-k}t;q)_\infty}
=\frac{(-abq^2;q)_\infty}{(-bqt;q)_\infty}
\frac{q^{k(k-1)/2}}{(bt)^k(-1/bt;q)_k}.
$$
Thus,
$$
G(x,t; a,b;\, q)=\frac{(-abq^2;q)_\infty}{(-bqt;q)_\infty}
\sum_{k=0}^\infty \frac{(aqx^{-1};q)_k} {(-1/bt,q;q)_k}\left(
\frac xb \right) ^k
$$   $$
 =\frac{(-abq^2;q)_\infty}{(-bqt;q)_\infty}\,
{}_2\phi_1 (aqx^{-1},\, 0;\, -1/bt;\; q,x/b). \eqno (6.32)
$$
This gives a desired generating function for big $q$-Laguerre
polynomials.
\bigskip

\noindent{\bf 6.5. Biorthogonal systems of functions}
\bigskip

Note that the operator $A$ from formula (6.1) can be written in
the form
$$
A =\alpha q^{J_0/4}\,\left(\sqrt{1-bQ}J_+ \,Q^{1/2} + Q^{1/2} J_-
\sqrt{1-bQ}\,\right)\, q^{J_0/4} - \beta_1 q^{2J_0} + \beta_2 Q ,
\eqno (6.33)
$$
where
$$
\alpha =(-abq)^{1/2}(1-q),\ \ \ \beta_1=b(1+q),\ \ \
\beta_2=bq+aq(b+1) ,
$$
and $J_\pm$, $Q$ and $q^{J_0}$ are the operators on ${\cal H}$
given as
$$
J_+\, f_n =\frac{a^{-1/4}q^{-n/2}}{1-q}
\sqrt{(1-q^{n+1})(1-aq^{n+1})} f_{n+1} ,
$$ $$
J_-\, f_n =\frac{a^{-1/4}q^{-(n-1
)/2}}{1-q}
\sqrt{(1-q^{n})(1-aq^{n})} f_{n-1},
$$  $$
q^{J_0}\, f_n =q^n(aq)^{1/2} \,f_n , \ \ \ \ Qf_n=q^nf_n.
$$
From the very beginning we could consider an operator
$$
A_1 :=\alpha q^{J_0/4}\,\left[(1-bQ)J_+ +Q J_- \right]\,
 q^{J_0/4} - \beta_1 q^{2J_0} + \beta_2 Q ,
$$
where $\alpha, \beta_1$, and $\beta_2$ are the same as above. This
operator is well defined, but it is not self-adjoint. Repeating
the reasoning of section 3, we find that eigenfunctions of $A_1$
are of the form
$$
\psi_\lambda(x)= \sum _{n=0}^\infty (-ab)^{-n/2}q^{-n}\left(
\frac{(aq;q)_n }{(q;q)_n} \right) ^{1/2} P_n( \lambda ;a,b;q) f_n
(x)
$$  $$
=\sum _{n=0}^\infty a^{-3n/4}(-b)^{-n/2}q^{-n} \frac{(aq;q)_n}
{(q;q)_n} P_n( \lambda ;a,b;q)x^n . \eqno (6.34)
$$
The last sum can be summed with the aid of formula (3.11.12) in
[KSw]. We thus have
$$
\psi_\lambda(x)=((-a/b^2)^{1/4}x;q)_\infty \cdot {}_2\phi_1
(bq\lambda^{-1},\, 0;\, bq;\ q, a^{-3n/4}(-b)^{-1/2}q^{-1}x\lambda
).
$$

Now we consider another operator
$$
A_2 :=\alpha q^{J_0/4}\left[ J_+ Q + J_- (1-bQ)\right] q^{J_0/4}
-\beta_1 q^{2J_0} +\beta_2 Q.
$$
This operator is adjoint to the operator $A_1: A_2^*=A_1$.
Repeating the reasoning of subsection 6.1, we find that
 eigenfunctions of $A_2$ have the form
$$
\varphi_\lambda(x)= \sum _{n=0}^\infty
(-ab)^{-n/2}q^{-n(n+1)/2}\left( \frac{(aq;q)_n
(bq;q)^2_n}{(q;q)_n} \right) ^{1/2} P_n( \lambda ;a,b;q) f_n (x)
$$  $$
=\sum _{n=0}^\infty a^{-3n/4}(-b)^{-n/2}q^{-n(n+1)/2}
\frac{(aq;q)_n (bq;q)_n} {(q;q)_n}  P_n( \lambda ;a,b;q)x^n. \eqno
(6.35)
$$
According to the formula (6.32), this function can be written as
$$
\varphi_\lambda(x)=\frac{(-abq^{2};q)_\infty}
{(a^{-3/4}(-b)^{1/2}x;q)_\infty}\, {}_2\phi_1(
aq/\lambda ,0;\, a^{3/4}(-b)^{-1/2}q/x;\;
q,\lambda /b).
$$

Let us denote by $\Psi_m(x)$, $m=0,\pm 1,\pm 2,\cdots$, the functions
$$
\Psi_m(x)=c_m\psi_{aq^{m+1}}(x),\ \ m=0,1,2,\cdots,\ \ \
\Psi_{-m}(x)=c'_{m-1}\psi_{bq^m}(x),\ \ m=1,2,\cdots,
 \eqno (6.36)
$$
and by $\Phi_m(x)$, $m=0,\pm 1,\pm 2,\cdots$, the functions
$$
\Phi_m(x)=c_m\varphi_{aq^{m+1}}(x), \ \ m=0,1,2,\cdots,\ \ \
\Phi_{-m}(x)=c'_{m-1}\varphi_{bq^{m}}(x), \ \ m=1,2,\cdots,
 \eqno (6.37)
$$
where $c_m$ and $c'_m$ are given by formulas (6.9) and (6.10).

Writing down the decompositions (6.34) and (6.35) for the
functions $\Psi_m(x)$ and $\Phi_m(x)$ (in terms of the orthonormal
basis $f_n$, $n=0,1,2,\cdots$, of the Hilbert space ${\cal H}$)
and taking into account the orthogonality relations (6.25)--(6.27)
we find that
$$
\langle \Psi_m(x),\Phi_n(x)\rangle =\delta_{mn}, \ \ \
m,n=0,\pm 1,\pm 2,\cdots .
$$
This means that we can formulate the following statement.
\medskip

{\bf Theorem 6.2.} {\it The set of functions $\Psi_m(x)$, $m=0,\pm
1,\pm 2,\cdots$, and the set of functions $\Phi_m (x)$, $m=0,\pm
1,\pm 2,\cdots$, form biorthogonal sets of functions with respect
to the scalar product in the Hilbert space ${\cal H}$.}
\bigskip

\noindent {\bf 7. ALTERNATIVE $q$-CHARLIER POLYNOMIALS AND THEIR
DUALS}

\bigskip

\noindent{\bf 7.1. Pair of operators $(B_1,J)$}
\bigskip

Let ${\cal H}$ be the same separable complex Hilbert space as
before. We have introduced into this space the orthonormal basis
$f_n$, $n = 0,1,2,\cdots$, expressed in terms of monomials in $x$.
We define on ${\cal H}$ two operators. The first one, denoted as
$Q$, acts on the basis elements as
$$
Qf_n =q^n f_n.
$$
The second one, denoted as $B_1$, is given by the formula
$$
B_1f_n  =a_n f_{n+1} +a_{n-1}f_{n-1} +b_n f_n , \eqno (7.1)
$$
with
$$
a_n=-(aq^{3n+1})^{1/2}\frac{ \sqrt{(1-q^{n+1})(1+aq^{n})}}
{(1+aq^{2n+1})\sqrt{(1+aq^{2n}) (1+aq^{2n+2})}} ,
$$  $$
b_n=q^n\left( \frac{1+aq^n}{(1+aq^{2n}) (1+aq^{2n+1})}+aq^{n-1}
\frac{1-q^{n}}{(1+aq^{2n-1}) (1+aq^{2n})}  \right) ,
$$
where $a$ is a fixed positive number. Clearly, $B_1$ is a
symmetric operator.

Since $a_n\to 0$ and $b_n\to 0$ when $n\to \infty$, the operator
$B_1$ is bounded. We assume that it is defined on the whole
Hilbert space ${\cal H}$. For this reason, $B_1$ is a self-adjoint
operator. Let us show that $B_1$ is a Hilbert--Schmidt operator.
For the coefficients $a_n$ and $b_n$ from (7.1), we have
$a_{n+1}/a_n \to q^{3/2}$ and $b_{n+1}/b_n \to q$ when $n\to
\infty$.  Since $0<q<1$, for the sum of all matrix elements of the
operator $B_1$ in the basis $f_n$, $n = 0,1,2,\cdots$, we have
$\sum _n (2a_n+b_n)< \infty$. This means that $B$ is a
Hilbert--Schmidt operator. Thus, a spectrum of $B_1$ is discrete
and has a single accumulation point at 0. Moreover, a spectrum of
$B_1$ is simple, since $B_1$ is representable by a Jacobi matrix
with $a_n\ne 0$ (see [Ber], Chapter VII).

To find eigenfunctions $\xi_\lambda $ of the operator $B_1$, $B_1
\xi_\lambda =\lambda \xi_\lambda $, we set
$$
\xi_\lambda =\sum _n\beta_n(\lambda)f_n,
$$
where $\beta_n(\lambda)$ are appropriate numerical coefficients.
Acting by the operator $B_1$ upon both sides of this relation, one
derives that
$$
\sum _{n=0}^{\infty}\, \beta_n(\lambda)\, (a_nf_{n+1}
+a_{n-1}f_{n-1} +b_nf_n )= \lambda \sum_{n=0}^{\infty}\,
\beta_n(\lambda) f_n,
$$
where $a_n$ and $b_n$ are the same as in (7.1). Collecting in this
identity all factors, which multiply $f_n$ with fixed $n$, one
derives the recurrence relation for the coefficients
$\beta_n(\lambda)$:
$$
\beta_{n+1}(\lambda)a_n +\beta_{n-1}(\lambda)a_{n-1}+
\beta_{n}(\lambda)b_n= \lambda \beta_{n}(\lambda).
$$
The substitution
$$
\beta_{n}(\lambda)=\left( \frac{(-a;q)_n\,(1+aq^{2n})} {(q;q)_n\,
(1+a)(a/q)^n}\right) ^{1/2} q^{-n(n+3)/4} \beta'_{n}(\lambda)
$$
reduces this relation to the following one
$$
-A_n \beta'_{n+1}(\lambda)- C_n \beta'_{n-1}(\lambda)
+(A_n+C_n)\beta'_{n}(\lambda)=\lambda \beta'_{n}(\lambda),
$$
$$
A_n=q^n\frac{1+aq^{n}}{(1+aq^{2n}) (1+aq^{2n+1})}, \ \ \ \
C_n=aq^{2n-1}\frac{1-q^{n}}{(1+aq^{2n-1}) (1+aq^{2n})}.
$$
This is the recurrence relation for the alternative $q$-Charlier
polynomials
$$
K_n(\lambda ;a;q):={}_2\phi_1 (q^{-n}, -aq^{n};\; 0; \; q,q\lambda
)
$$
(see, formulas (3.22.1) and (3.22.2) in [KSw]). Therefore,
$\beta'_n(\lambda )=K_n(\lambda ;a;q)$ and
$$
\beta_n(\lambda )=\left( \frac{(-a;q)_n\,(1+aq^{2n})} {(q;q)_n\,
(1+a)a^n}\right) ^{1/2} q^{-n(n+1)/4} K_n(\lambda ;a;q). \eqno
(7.2)
$$
For the eigenvectors $\xi _\lambda$ we thus have the expression
$$
\xi _\lambda =\sum_{n=0}^\infty\, \left(
\frac{(-a;q)_n\,(1+aq^{2n})} {(q;q)_n\, (1+a)a^n}\right) ^{1/2}
q^{-n(n+1)/4} K_n(\lambda ;a;q) f_n .
 \eqno (7.3)
$$
Since the spectrum of the operator $B_1$ is discrete, only for a
discrete set of values of $\lambda$ these vectors belong to the
Hilbert space ${\cal H}$.

Now we look for a spectrum of the operator $B_1$ and for a set of
polynomials, dual to alternative $q$-Charlier polynomials. To this
end we use the action of the operator
$$
J:= Q^{-1} - a\,Q
$$
upon the eigenvectors $\xi _\lambda$, which belong to the Hilbert
space ${\cal H}$. In order to find how this operator acts upon
these vectors, one can use the $q$-difference equation
$$
(q^{-n}-aq^{n})K_n(\lambda)=- aK_{n} (q\lambda
)+\lambda^{-1}K_n(\lambda)-\lambda^{-1}
(1-\lambda)K_n(q^{-1}\lambda) \eqno(7.4)
$$
for the alternative $q$-Charlier polynomials $K_n(\lambda)\equiv
K_n(\lambda ;a;q)$ (see formula (3.22.5) in [KSw]). Multiply both
sides of (7.4) by $k_n\,f_n$ and sum up over $n$, where $k_n$ are
the coefficients of $K_n(\lambda ;a;q)$ in the expression (7.2)
for $\beta_n(\lambda)$. Taking into account the formula (7.3) and
the fact that $Jf_n=(q^{-n}-aq^{n})f_n$, one obtains the relation
$$
J\,\xi _{\lambda}= -a\,\xi _{q\lambda}+ \lambda^{-1}\, \xi
_{\lambda}- \lambda^{-1}(1-\lambda)\, \xi_{q^{-1}\lambda}, \eqno
(7.5)
$$
which will be used in the next subsection.
\bigskip

\noindent{\bf 7.2. Spectrum of $B_1$ and orthogonality of
alternative $q$-Charlier polynomials}
\medskip

The aim of this section is to find, by using the operators $B_1$
and $J$, a basis in the Hilbert space ${\cal H}$, which consists
of eigenfunctions of the operator $B_1$ in a normalized form, and
to derive explicitly the unitary matrix $U$, connecting this basis
with the basis $f_n$, $n=0,1,2,\cdots$, in ${\cal H}$. This matrix
leads directly to the orthogonality relation for alternative
$q$-Charlier polynomials. For this purpose we first find a
spectrum of $B_1$.

We proceed as in the previous cases. First we analyze a form of
the spectrum of the operator $B_1$ from the point of view of the
spectral theory of Hilbert--Schmidt operators. If $\lambda$ is a
spectral point of the operator $B_1$, then (as it is easy to see
from (7.5)) a successive action by the operator $J$ upon the
function (eigenfunction of $B_1$) $\xi_\lambda$ leads to the
eigenfunctions $\xi_{q^m\lambda}$, $m=0,\pm 1, \pm 2,\cdots$.
However, since $B_1$ is a Hilbert--Schmidt operator, not all of
these points may belong to the spectrum of $B_1$, since
$q^{-m}\lambda \to\infty$ when $m\to +\infty$ once $\lambda\ne 0$.
This means that the coefficient $1-\lambda' $ of $\xi
_{q^{-1}\lambda'}$ in (7.5) must vanish for some eigenvalue
$\lambda'$. Clearly, it vanishes when $\lambda' =1$. Moreover,
this is the only possibility for the coefficient of $\xi
_{q^{-1}\lambda'}$ in (7.5) to vanish, that is, the point $\lambda
=1$ is a spectral point for the operator $B_1$. Let us show that
the corresponding eigenfunction $\xi _{1}\equiv \xi_{q^{0}}$
belongs to the Hilbert space ${\cal H}$.

By formula (II.6) of Appendix II in [GR], one has $K_n(1
;a;q)={}_2\phi_1 (q^{-n}, -aq^{n};\; 0; \; q,q )=(-a)^nq^{n^2}$.
Therefore,
$$
\langle \xi_1,\xi_1\rangle = \sum_{n=0}^\infty
\frac{(-a;q)_n(1+aq^{2n})}{(1+a)(q;q)_na^nq^{n(n+1)/2}}K^2_n(1
;a;q) = \sum_{n=0}^\infty
\frac{(-a;q)_n(1+aq^{2n})a^n}{(1+q)(q;q)_n q^{-n(3n-1)/2}} .
 \eqno (7.6)
$$
In order to calculate this sum, we take the limit $d,e\to \infty$
in the equality
$$
 \sum_{n=0}^\infty
\frac{(1+aq^{2n})(-a;q)_n(d;q)_n(e;q)_n}{(1+a)(-aq/d;q)_n(-aq/e;q)_n(q;q)_n}
 \left( \frac{aq}{de}\right) ^n q^{n(n-1)/2} =\frac{(-aq;q)_\infty
 (-aq/de;q)_\infty}{(-aq/d;q)_\infty (-aq/e;q)_\infty}
 $$
(see formula in Exercise 2.12, Chapter 2 of [GR]).  Since
$$
\lim_{d,e\to \infty}(d;q)_n(e;q)_n (aq/de)^n=q^{n(n-1)}(aq)^n ,
$$
we obtain from here that the sum in (7.6) is equal to $
(-aq;q)_\infty$, that is, $\langle \xi_1,\xi_1\rangle <\infty$ and
$\xi_1$ belongs to the Hilbert space ${\cal H}$. Thus, the point
$\lambda =1$ does belong to the spectrum of the operator $B_1$.

Let us find other spectral points of the operator $B_1$ (recall
that a spectrum of $B_1$ is discrete). Setting $\lambda = 1$ in
(7.5), we see that the operator $J$ transforms $\xi _{q^0}$ into a
linear combination of the vectors $\xi _{q}$ and $\xi_{q^0}$.
Moreover, $\xi_q$ belongs to the Hilbert space ${\cal H}$, since
the series
$$
\langle \xi_q,\xi_q\rangle = \sum_{n=0}^\infty
\frac{(-a;q)_n\,(1+aq^{2n})}{(1+a)(q;q)_n\,a^n}q^{-n(n+1)/2}
\,K^2_n(q;a;q)
$$
is majorized by the corresponding series (7.6) for $\xi_{q^0}$.
Therefore, $\xi _{q}$ belongs to the Hilbert space ${\cal H}$ and
the point $q$ is an eigenvalue of the operator $B_1$. Similarly,
setting $\lambda=q$ in (7.5), we find that $\xi _{q^2}$ is an
eigenvector of $B_1$ and the point $q^2$ belongs to the spectrum
of $B_1$. Repeating this procedure, we find that all $\xi _{q^n}$,
$n=0,1,2,\cdots$, are eigenvectors of $B_1$ and the set $q^n$,
$n=0,1,2,\cdots$, belongs to the spectrum of $B_1$. So far, we do
not know yet whether other spectral points exist or not.

The functions $\xi _{q^n}$, $n=0,1,2,\cdots$, are linearly
independent elements of the Hilbert space ${\cal H}$ (since they
correspond to different eigenvalues of the self-adjoint operator
$B_1$). Suppose that values $q^n$, $n=0,1,2,\cdots$, constitute a
whole spectrum of $B_1$. Then the set of vectors $\xi _{q^n}$,
$n=0,1,2,\cdots$, is a basis in the Hilbert space ${\cal H}$.
 Introducing the notation $\Xi _k:=\xi_{q^k}$, $k=0,1,2,\cdots$, we
find from (7.5) that
$$
J \,\Xi _k = - a \,\Xi _{k+1} + q^{-k}\, \Xi _k - q^{-k}(1-q^k)\,
\Xi _{k-1} .
$$
As we see, the matrix of the operator $J$ in the basis $\Xi _k$,
$k=0,1,2,\cdots$, is not symmetric, although in the initial basis
$|n\rangle$, $n=0,1,2,\cdots$, it was symmetric. The reason is
that the matrix $(a_{mn})$ with entries $a_{mn}:=\beta_m(q^n)$,
$m,n=0,1,2,\cdots$, where $\beta_m(q^n)$ are the coefficients
(7.2) in the expansion $\xi _{q^n}=\sum _m \,\beta_m(q^n)f_n$, is
not unitary. This fact is equivalent to the statement that the
basis $\Xi _n=\xi_{q^n}$, $n=0,1,2,\cdots$, is not normalized. To
normalize it, one has to multiply $\Xi _n$ by corresponding
numbers $c_n$. Let $\hat\Xi _n = c_n\Xi _n$, $n=0,1,2,\cdots$, be
a normalized basis. Then the matrix of the operator $J$ is
symmetric in this basis. Since  $J$ has in the basis $\{ \hat\Xi
_n\}$ the form
$$
J\, \hat\Xi _n = -c_{n+1}^{-1}c_n a\, \hat\Xi _{n+1} + q^{-n}\,
\hat\Xi _n - c_{n-1}^{-1}c_n q^{-n}(1-q^n) \,\hat\Xi_{n-1} ,
$$
then its symmetricity means that
$c_{n+1}^{-1}c_na=c_{n}^{-1}c_{n+1} q^{-n-1}(1-q^{n+1})$, that is,
$c_{n}/c_{n-1} =\sqrt{aq^n/(1-q^n)}$. Therefore,
$$
c_n= c(a^nq^{n(n+1)/2}/(q;q)_n )^{1/2},
$$
where $c$ is a constant.

The expansions
$$
\hat\xi _{q^n}(x)\equiv \hat\Xi _n(x)= \sum _m
c_n\beta_m(q^n)|m\rangle \equiv \sum _m {\hat a}_{mn}|m\rangle
\eqno (7.7)
$$
connect two orthonormal bases in the Hilbert space ${\cal H}$.
This means that the matrix $({\hat a}_{mn})$, $m,n=0,1,2,\cdots$,
with entries
$$
{\hat a}_{mn}=c_n\beta _m(q^n)= c\left(
\frac{a^nq^{n(n+1)/2}}{(q;q)_n}
\frac{(-a;q)_m\,(1+aq^{2m})}{(1+a)(q;q)_m\,a^m q^{m(m+1)/2}}
\right)^{1/2} K_m(q^n ;a;q) \eqno (7.8)
$$
is unitary, provided  that the constant $c$ is appropriately
chosen. In order to calculate this constant, we use the relation
$\sum_{m=0}^\infty |{\hat a}_{mn}|^2=1$ for $n=0$. Then this sum
is a multiple of the sum in (7.6) and, consequently,
$$
c=(-aq;q)^{-1/2}_\infty.
$$

The matrix $({\hat a}_{mn})$ is real and orthogonal, that is,
$$
\sum _n {\hat a}_{mn}{\hat a}_{m'n}=\delta_{mm'},\ \ \ \ \sum _m
{\hat a}_{mn}{\hat a}_{mn'}=\delta_{nn'} .   \eqno (7.9)
$$
Substituting into the first sum over $n$ in (7.9) the expressions
for ${\hat a}_{mn}$, we obtain the identity
$$
\sum_{n=0}^\infty
\frac{a^nq^{n(n+1)/2}}{(q;q)_n}\,K_m(q^n;a;q)\,K_{m'}(q^n;a;q) =
\frac{(-aq^m;q)_\infty \, a^m\,(q;q)_m} {(1+aq^{2m})}\,
q^{m(m+1)/2} \delta_{mm'}\,, \eqno (7.10)
$$
which must yield the orthogonality relation for alternative
$q$-Charlier polynomials. An only gap, which remains to be
clarified, is the following. We have assumed that the points
$q^n$, $n=0,1,2,\cdots$, exhaust the whole spectrum of $B_1$. Let
us show that this is the case.

Recall that the self-adjoint operator $I_1$ is represented by a
Jacobi matrix in the basis $f_n$, $n=0,1,2,\cdots$. According to
the theory of operators of such type, eigenvectors $\xi_\lambda$
of $B_1$ are expanded into series in the basis $f_n$,
$n=0,1,2,\cdots$, with coefficients, which are polynomials in
$\lambda$. These polynomials are orthogonal with respect to some
positive measure $d\mu (\lambda)$ (moreover, for self-adjoint
operators this measure is unique). The set (a subset of ${\Bbb
R}$), on which these polynomials are orthogonal, coincides with
the spectrum of the operator under consideration and the spectrum
is simple.

We have found that the spectrum of $B_1$ contains the points
$q^n$, $n=0,1,2,\cdots$. If the operator $B_1$ had other spectral
points $x$, then on the left-hand side of (7.10) there would be
other summands $\mu_{x_k}\, K_m({x_k};a;q)\,K_{m'}({x_k};a;q)$,
corresponding to these additional points. Let us show that these
additional summands do not appear. We set $m=m'=0$ in the relation
(7.10) with the additional summands. Since $K_0(x;a;q)=1$, we have
the equality
$$
\sum_{n=0}^\infty \frac{a^nq^{n(n+1)/2}}{(q;q)_n} + \sum_k
\mu_{x_k} =(-aq;q)_\infty  .
$$
According to the formula for the $q$-exponential function $E_q(a)$
(see formula (II.2) of Appendix II in [GR]), we have
$\sum_{n=0}^\infty a^nq^{n(n+1)/2}(q;q)^{-1}_n =(-aq;q)_\infty$.
Hence, $\sum_k \mu_{x_k} =0$. This means that additional summands
do not appear in (7.10) and it does represent the orthogonality
relation for alternative $q$-Charlier polynomials.

Due to the orthogonality relation for the alternative $q$-Charlier
polynomials, we arrive at the following statement:
\medskip

{\bf Proposition 8.1.} {\it The spectrum of the operator $B_1$
coincides with the set of points $q^{n}$, $n=0,1,2,\cdots$. The
spectrum is simple and has one accumulation point at 0.}
\bigskip

\noindent {\bf 7.3. Dual alternative $q$-Charlier polynomials}
\bigskip

Now we consider the second identity in (7.9), which gives the
orthogonality relation for the matrix elements ${\hat a}_{mn}$,
considered as functions of $m$. Up to multiplicative factors these
functions coincide with
$$
F_n(x;a|q)={}_2\phi_1 (x,-a/x;\; 0;\; q,q^{n+1}),  \eqno (7.11)
$$
considered on the set $x\in \{ q^{-m}\, |\, m=0,1,2,\cdots \}$.
Consequently,
$$
{\hat a}_{mn}= \left( \frac{a^nq^{n(n+1)/2}}{(q;q)_n}
\frac{(1+aq^{2m})}{(-aq^m;q)_\infty(q;q)_m\,a^m q^{m(m+1)/2}}
\right)^{1/2}   F_n(q^{-m} ;a|q)
$$
and the second identity in (7.9) gives the orthogonality relation
for $F_n(q^{-m} ;a|q)$:
$$
\sum_{m=0}^\infty \frac{(1+aq^{2m})}{a^m(-aq^m;q)_\infty (q;q)_m
q^{m(m+1)/2}} F_n(q^{-m};a|q)F_{n'}(q^{-m};a|q) =\frac{
(q;q)_n}{a^{n}q^{n(n+1)/2}}\, \delta_{nn'}. \eqno (7.12)
$$

The functions $F_n(x;a,b|q)$ can be represented in another form.
Indeed, taking in the relation (III.8) of Appendix III in [GR] the
limit $c\to \infty$, one derives the relation
$$
{}_2\phi_1 (q^{-m},-aq^m;\; 0;\; q,q^{n+1})=(-a)^mq^{m^2}
{}_3\phi_0 (q^{-m},-aq^m, q^{-n}\; -\, ;\; q,-q^{n}/a) .
$$
Therefore, we have
$$
F_{n}(q^{-m} ;a|q)=(-a)^mq^{m^2}{}_3\phi_0
(q^{-m},-aq^{m},q^{-n};\; -\, ; \; q,-q^n/a) .
 \eqno (7.13)
$$
The basic hypergeometric function ${}_3\phi_0$ in (7.13) is a
polynomial of degree $n$ in the variable $\mu(m): = q^{-m}-
a\,q^{m}$, which represents a $q$-quadratic lattice; we denote it
by
$$
d_n(\mu (m); a;q):= {}_3\phi_0(q^{-m},-a\,q^{m},q^{-n};\; -\, ;\;
q,-q^n/a)\,.  \eqno (7.14)
$$
Then formula (7.12) yields the orthogonality relation
$$
\sum_{m=0}^\infty \frac{(1+aq^{2m})a^m}{(-aq^m;q)_\infty
(q;q)_m}\, q^{m(3m-1)/2} d_n(\mu(m); a;q)
d_{n'}(\mu(m);a;q)=\frac{ (q;q)_n}{a^nq^{n(n+1)/2}}\, \delta_{nn'}
\eqno (7.15)
$$
for the polynomials (7.14) when $a>0$. We call the polynomials
$d_n(\mu (m); a;q)$ {\it dual alternative $q$-Charlier
polynomials}. Thus, the following theorem holds.
 \medskip

{\bf Theorem 7.1.} {\it The polynomials $d_n(\mu (m); a;q)$, given
by formula (7.14), are orthogonal on the set of points $\mu(m): =
q^{-m}- a\,q^{m}$, $m=0,1,2,\cdots$, and the orthogonality
relation is given by formula (7.15).}
 \medskip

Let ${\frak l}^2$ be the Hilbert space of functions on the set
$m=0,1,2,\cdots$ with the scalar product
$$
\langle F_1,F_2\rangle = \sum _{m=0}^\infty\,
\frac{(1+aq^{2m})a^m}{(-aq^m;q)_\infty(q;q)_m}\, q^{m(3m-1)/2}
\,F_1(m)\,\overline{F_2(m)} ,   \eqno (7.16)
$$
where the weight function is taken from (7.15). The polynomials
(7.14) are in one-to-one correspondence with the columns of the
orthogonal matrix $({\hat a}_{mn})$ and the orthogonality relation
(7.15) is equivalent to the orthogonality of these columns. Due to
(7.9) the columns of the matrix $({\hat a}_{mn})$ form an
orthonormal basis in the Hilbert space of sequences ${\bf a}=\{
a_n\, |\, n=0,1,2,\cdots \}$ with the scalar product $\langle {\bf
a},{\bf a}'\rangle=\sum_n a_n{\overline{a'_n}}$. This scalar
product is equivalent to the scalar product (7.16) for the
polynomials $d_n(\mu (m);a;q)$. For this reason, the set of
polynomials $d_n(\mu (m);a;q)$, $n=0,1,2,\cdots$, form an
orthogonal basis in the Hilbert space ${\frak l}^2$. This means
that {\it either the dual alternative $q$-Charlier polynomials
$d_n(\mu (m);a;q)$ correspond to determinate moment problem or the
point measure in (7.15) is extremal if these polynomials
correspond to indeterminate moment problem}. This question will
not be further pursued here.

A recurrence relation for the polynomials $d_n(\mu (m);a;q)$ is
derived from (7.4). It has the form
$$
(q^{-m}- a q^{m})d_n(\mu (m))= - a d_{n+1}(\mu (m)) +
q^{-n}d_{n}(\mu (m)) - q^{-n}(1-q^{n}) d_{n-1}(\mu (m)), \eqno
(7.17)
$$
where $d_{n}(\mu (m))\equiv d_{n}(\mu(m);a;q)$. A $q$-difference
equation for $d_{n}(\mu(m);a;q)$ can be obtained from the
three-term recurrence relation for alternative $q$-Charlier
polynomials.

Note that for the polynomials $d_n(\mu (m);a;q^{-1})$ with $q<1$
we have the expression
$$
d_n(\mu (m);a;q^{-1})={}_3\phi_2(q^{-m},-a\,q^{m},q^{-n};\; 0,0\,
;\; q,q) .  \eqno (7.18)
$$
However, the recurrence relation for these polynomials (which can
be obtained from the relation (7.17)), does not satisfy the
positivity condition $A_nC_{n+1}>0$, that is, they are not
orthogonal polynomials for $a>0$ (as it is the case for
alternative $q$-Charlier polynomials). This positivity condition
holds only if we require $a<0$. In this case, the polynomials
(7.18) are the continuous big $q$-Hermite polynomials $H_n(x;a|q)$
(for an explicit form of these polynomials see, for example,
[KSw], formula (3.18.1)), which are orthogonal on a certain
continuous set.

\bigskip

\noindent{\bf 8. DUALITY OF AL-SALAM--CARLITZ I AND\\
$q$-CHARLIER POLYNOMIALS}
 \bigskip

\noindent{\bf 8.1. Pair of operators $(B_2,Q^{-1})$}
\bigskip

Let $a$ be a real number such that $a<0$. Let ${\cal L}$ be the
separable complex Hilbert space with the orthonormal basis
$|n\rangle$, $n = 0,1,2,\cdots$. We define on ${\cal L}$ the
operator $B_2$, which is given by the formula
$$
B_2 |n\rangle  =a_n |n+1\rangle +a_{n-1}|n-1\rangle -b_n |n\rangle
, \eqno (8.1)
$$
with
$$
a_n=(-a)^{1/2}q^{n/2} \sqrt{(1-q^{n+1})}  ,\ \ \ \ b_n=(a+1)q^n.
$$
Clearly, $B_2$ is a bounded symmetric operator. Therefore, we
assume that it is defined on the whole Hilbert space ${\cal L}$.
For this reason, $B_2$ is a self-adjoint operator. As in the
previous cases, it is easy to show that $B_2$ is a
Hilbert--Schmidt operator. Thus, a spectrum of $B_2$ is discrete
and has a single accumulation point at 0. Moreover, a spectrum of
$B_2$ is simple, since $B_2$ is representable by a Jacobi matrix
with $a_n\ne 0$.

To find eigenfunctions $\xi_\lambda $ of the operator $B_2$, $B_2
\xi_\lambda =\lambda \xi_\lambda $, we set $\xi_\lambda =\sum _n
\beta_n(\lambda)|n \rangle$, where $\beta_n(\lambda)$ are
appropriate numerical coefficients. Acting by the operator $B_2$
upon both sides of this relation, one derives that
$$
\sum _{n=0}^{\infty}\, \beta_n(\lambda)\, [a_n|n+1\rangle
+a_{n-1}|n-1\rangle -b_n|n\rangle ]= \lambda \sum_{n=0}^{\infty}\,
\beta_n(\lambda) |n\rangle,
$$
where $a_n$ and $b_n$ are the same as in (8.1). Collecting in this
identity all factors, which multiply $|n\rangle $ with fixed $n$,
one derives the recurrence relation for the coefficients
$\beta_n(\lambda)$:
$$
\beta_{n+1}(\lambda)a_n +\beta_{n-1}(\lambda)a_{n-1}-
\beta_{n}(\lambda)b_n= \lambda\, \beta_{n}(\lambda).
$$
Making the substitution
$$
\beta_{n}(\lambda)=\left( \frac{q^{-n(n-1)/2}} {(q;q)_n\,
(-a)^n}\right) ^{1/2}  \beta'_{n}(\lambda) ,
$$
we reduce this relation to the following one
$$
 \beta'_{n+1}(\lambda)+ (-a)q^{n-1}(1-q^n) \beta'_{n-1}(\lambda)
-(a+1)q^n\beta'_{n}(\lambda)=\lambda\, \beta'_{n}(\lambda).
$$
This is the recurrence relation for the Al-Salam--Carlitz I
polynomials
$$
U^{(a)}_n(\lambda ;q):=(-a)^nq^{n(n-1)/2} {}_2\phi_1 (q^{-n},
\lambda^{-1};\; 0; \; q,q\lambda/a )   \eqno (8.2)
$$
(see, formulas (3.24.1) in [KSw]). Therefore, $\beta'_n(\lambda
)=U^{(a)}_n(\lambda ;q)$ and
$$
\beta_n(\lambda )=\left( \frac{q^{-n(n-1)/2}} {(q;q)_n\,
(-a)^n}\right) ^{1/2} U^{(a)}_n(\lambda ;q). \eqno (8.3)
$$
For the eigenvectors $\xi _\lambda$ we thus have the expression
$$
\xi _\lambda =\sum_{n=0}^\infty\, \left( \frac{q^{-n(n-1)/2}}
{(q;q)_n\, (-a)^n}\right) ^{1/2} U^{(a)}_n(\lambda ;q) |n\rangle .
 \eqno (8.4)
$$
Since the spectrum of the operator $B_2$ is discrete, only for a
discrete set of values of $\lambda$ these vectors belong to the
Hilbert space ${\cal L}$. This discrete set of eigenvectors
determines a spectrum of $B_2$.

Let us find a spectrum of the operator $B_2$ and a set of
polynomials, dual to Al-Salam--Carlitz I polynomials. For this
purpose we use the operator $Q^{-1}$, which is diagonal in the
basis $\{|n\rangle \}$, and is given as
$$
Q^{-1}|n\rangle =q^{-n}|n\rangle .
$$
We have to find how the operator $Q^{-1}$ acts upon the
eigenvectors $\xi _\lambda$, which belong to the Hilbert space
${\cal L}$. To this end, one can use the $q$-difference equation
for Al-Salam--Carlitz I polynomials, which can be written as
$$
q^{-n} U^{(a)}_n(\lambda ;q)= aq^{-1}\lambda^{-2}
U^{(a)}_n(q\lambda ;q) -d_\lambda U^{(a)}_n(\lambda
;q)+a\lambda^{-2} (1-\lambda)(1-\lambda/a)U^{(a)}_n(q^{-1}\lambda
;q), \eqno(8.5)
$$
where $d_\lambda=a\lambda^{-2}(1+q^{-1}-\lambda-\lambda/q)$.

Multiply both sides of (8.5) by $k_n\,|n\rangle$ and sum up over
$n$, where $k_n$ are the coefficients of $U^{(a)}_n(\lambda ;q)$
in the expression (8.3) for $\beta_n(\lambda)$. Taking into
account formula (8.4) and the fact that $Q^{-1}|n\rangle
=q^{-n}|n\rangle$, one obtains the relation
$$
Q^{-1}\,\xi _{\lambda}= aq^{-1}\lambda^{-2}\xi
_{q\lambda}-d_\lambda  \xi _{\lambda}+
a\lambda^{-2}(1-\lambda/a)(1-\lambda)\, \xi_{q^{-1}\lambda}, \eqno
(8.6)
$$
which is used in the next subsection.

\medskip

\noindent{\bf 8.2. Spectrum of $B_2$ and orthogonality of
Al-Salam--Carlitz polynomials}
\medskip

Let us analyze a form of the spectrum of $B_2$. If $\lambda$ is a
spectral point of the operator $B_2$, then (as it is easy to see
from (8.6)) a successive action by the operator $Q^{-1}$ upon the
vector (eigenvector of $B_2$) $\xi_\lambda$ leads to the
eigenvectors $\xi_{q^m\lambda}$, $m=0,\pm 1, \pm 2,\cdots$.
However, since $B_2$ is a Hilbert--Schmidt operator, not all of
these points may belong to the spectrum of $B_2$, since
$q^{-m}\lambda \to\infty$ when $m\to +\infty$ if $\lambda\ne 0$.
This means that the coefficient of $\xi _{q^{-1}\lambda'}$ in
(8.6) must vanish for some eigenvalue $\lambda'$. There are two
such values of $\lambda$: $\lambda=1$ and $\lambda=a$. Let us show
that both of these points are spectral points of $B_2$. Observe
that $U^{(a)}_n(1 ;q)=(-a)^nq^{n(n-1)/2}$ and $U^{(a)}_n(a
;q)=(-1)^nq^{n(n-1)/2}$. Hence, for the scalar product $\langle
\xi_1,\xi_1\rangle$ we have the expression
$$
 \sum_{n=0}^\infty
\frac{q^{-n(n-1)/2}}{(q;q)_n(-a)^n}(-a)^{2n}q^{n(n-1)}  =
 \sum_{n=0}^\infty
\frac{q^{-n(n-1)/2}}{(q;q)_n}(-a)^{n} =(a;q)_\infty .
 \eqno (8.7)
$$
Similarly, for $\langle \xi_a,\xi_a\rangle$ one has the expression
$$
\langle \xi_a,\xi_a\rangle=  \sum_{n=0}^\infty
\frac{q^{-n(n-1)/2}}{(q;q)_n}(-a)^{-n}  = (1/a;q)_\infty .
 \eqno (8.8)
$$
Thus, the values $\lambda=1$ and $\lambda=a$ are spectral points
of the operator $B_2$.

Let us find other spectral points of $B_2$. Setting $\lambda = 1$
in (8.6), we see that the operator $Q^{-1}$ transforms $\xi
_{q^0}$ into a linear combination of the vectors $\xi _{q}$ and
$\xi_{1}$. We have to show that $\xi _{q}$ belongs to the Hilbert
space ${\cal L}$, that is, that
$$
\langle \xi_q,\xi_q\rangle = \sum_{n=0}^\infty
\frac{q^{n(n-1/2}}{(q;q)_n\,(-a)^n}(-a)^{2n}q^{n(n-1)/2}\,U^{(a)}_n(q
;q)^2<\infty .
$$
It is made in the same way as in the case of the scalar product
$\langle \psi_{aq^2},\psi_{aq^2} \rangle$ in subsection 4.2.
Therefore, $\xi _{q}$ belongs to the Hilbert space ${\cal L}$ and
the point $q$ is an eigenvalue of the operator $B_2$. Similarly,
setting $\lambda=q$ in (8.6), we find that $\xi _{q^2}$ is an
eigenvector of $B_2$ and the point $q^2$ belongs to the spectrum
of $B_2$. Repeating this procedure, we find that all $\xi _{q^n}$,
$n=0,1,2,\cdots$, are eigenvectors of $B_2$ and the set $q^n$,
$n=0,1,2,\cdots$, belongs to the spectrum of $B_2$. Likewise, one
conclude that the elements $\xi _{aq^n}$, $n=0,1,2,\cdots$, are
eigenvectors of $B_2$ and the set $aq^n$, $n=0,1,2,\cdots$,
belongs to the spectrum of $B_2$. So far, we do not know yet
whether other spectral points exist or not.

The vectors $\xi _{q^n}$ and $\xi _{aq^n}$, $n=0,1,2,\cdots$, are
linearly independent elements of the Hilbert space ${\cal L}$.
Suppose that values $q^n$ and $aq^n$,  $n=0,1,2,\cdots$,
constitute a whole spectrum of $B_2$. Then the set of vectors $\xi
_{q^n}$ and $\xi _{aq^n}$, $n=0,1,2,\cdots$, is a basis in the
Hilbert space ${\cal L}$. Introducing the notations $\Xi
_k:=\xi_{q^k}$ and $\Xi'_k:=\xi_{aq^k}$, $k=0,1,2,\cdots$, we find
from (8.6) that
$$
Q^{-1} \,\Xi _n =  aq^{-2n-1} \,\Xi _{n+1} - d_n\, \Xi _n +
aq^{-2n}(1-q^n)(1-q^n/a)\, \Xi _{n-1} ,
$$   $$
Q^{-1} \,\Xi' _n =  a^{-1}q^{-2n-1} \,\Xi' _{n+1} - d'_n\, \Xi' _n
+ a^{-1}q^{-2n}(1-q^n)(1-aq^n)\, \Xi' _{n-1} ,
$$
where
$$
d_n=aq^{-2n}(q^{-1}+1-q^n-q^{n-1}), \ \ \ \
d'_n=a^{-1}q^{-2n}(q^{-1}+1-aq^n-aq^{n-1}).
$$
As we see, the matrix of the operator $Q^{-1}$ in the basis $\Xi
_n$, $\Xi'_n$, $n=0,1,2,\cdots$, is not symmetric, although in the
initial basis $|m\rangle$, $m=0,1,2,\cdots$, it was symmetric. The
reason is that the matrix $M:=((a_{mn})_{m,n=0}^\infty\ \;
(a'_{mn})_{m,n=0}^\infty)\equiv ((a_{mn})\ \; (a'_{mn}))$ with
entries
$$
a_{mn}:=\beta_m(q^n),\ \ \ a'_{mn}:=\beta_m(aq^n), \ \ \ \
m,n=0,1,2,\cdots,
$$
where  $\beta_m(dq^n)$, $d=1,a$,  are the coefficients (8.3) in
the expansion $\xi _{dq^n}=\sum _m \,\beta_m(dq^n)|n\rangle$, is
not unitary. This fact is equivalent to the statement that the
basis $\Xi _n=\xi_{q^n}$, $\Xi'_n=\xi_{aq^n}$, $n=0,1,2,\cdots$,
is not normalized. To normalize it, one has to multiply $\Xi _n$
by corresponding numbers $c_n$ and $\Xi'_n$ by corresponding
numbers $c'_n$. Let $\hat\Xi _n = c_n\Xi _n$ and $\hat\Xi' _n =
c'_n\Xi' _n$, $n=0,1,2,\cdots$, be a normalized basis. Then the
matrix of the operator $Q^{-1}$ is symmetric in this basis. Since
$Q^{-1}$ has in the basis $\{ \hat\Xi _n, \hat\Xi' _n\}$ the form
$$
Q^{-1} \,\hat\Xi _n =c_{n+1}^{-1}c_n  aq^{-2n-1} \,\hat\Xi _{n+1}
- d_n\, \hat\Xi _n +c_{n-1}^{-1}c_n aq^{-2n}(1-q^n)(1-q^n/a)\,
\hat\Xi _{n-1} ,
$$   $$
Q^{-1} \,\hat\Xi' _n ={c'}_{n+1}^{-1}c'_n  a^{-1}q^{-2n-1}
\,\hat\Xi' _{n+1} - d'_n\, \hat\Xi' _n +{c'}_{n-1}^{-1}c'_n
a^{-1}q^{-2n}(1-q^n)(1-aq^n)\, \hat\Xi' _{n-1} .
$$
then its symmetricity means that
$$
c_{n+1}^{-1}c_naq^{-2n-1}=c_{n}^{-1}c_{n+1}a q^{-2n-2}
(1-q^{n+1})(1-q^{n+1}/a),
$$  $$
{c'}_{n+1}^{-1}c'_na^{-1}q^{-2n-1}={c'}_{n}^{-1}c'_{n+1}a^{-1}
q^{-2n-2} (1-aq^{n+1})(1-q^{n+1}),
$$
that is,
$$
\frac{c_{n}}{c_{n-1}} =\sqrt{\frac{q}{(1-q^n)(1-q^n/a)}},\ \ \ \
\frac{c'_{n}}{c'_{n-1}} =\sqrt{\frac{q}{(1-q^n)(1-aq^n)}}\; .
$$
Therefore, for the coefficients $c_n$ and $c'_n$ we have the
expressions
$$
c_n= c\frac{q^{n/2}}{(q/a;q)_n ^{1/2}(q;q)_n ^{1/2}},\ \ \ \ c'_n=
c'\frac{q^{n/2}}{(aq;q)_n ^{1/2}(q;q)_n ^{1/2}},
$$
where $c$ and $c'$ are some constant.

Thus, in the expansions
$$
\hat\xi _{q^n}\equiv \hat\Xi _n= \sum _m c_n\beta_m(q^n)|m\rangle
\equiv \sum _m {\hat a}_{mn}|m\rangle , \eqno (8.9)
$$   $$
\hat\xi _{aq^n}\equiv \hat\Xi' _n= \sum _m
c'_n\beta_m(aq^n)|m\rangle \equiv \sum _m {\hat a}'_{mn}|m\rangle
, \eqno (8.10)
$$
the matrix $\hat M :=({\hat a}_{mn}\ {\hat a}'_{mn})$ with entries
$$
{\hat a}_{mn}=c_n\beta_m(q^n)=c\left( \frac{q^{n}}{(q/a;q)_n
(q;q)_n )} \frac{q^{-m(m-1)/2}} {(q;q)_m\, (-a)^m}\right)^{1/2}
U_m^{(a)}(q^n;q), \eqno (8.11)
$$  $$
{\hat a}'_{mn}=c'_n\beta_m(aq^n)=c'\left( \frac{q^{n}}{(aq;q)_n
(q;q)_n )} \frac{q^{-m(m-1)/2}} {(q;q)_m\, (-a)^m}\right)^{1/2}
U_m^{(a)}(aq^n;q), \eqno (8.12)
$$
is unitary, provided that the constants $c$ and $c'$ are
appropriately chosen. In order to calculate these constants, we
use the relations $\sum_{m=0}^\infty |{\hat a}_{mn}|^2=1$ and
$\sum_{m=0}^\infty |{{\hat a}'}_{mn}|^2=1$ for $n=0$. Then these
sums are  multiples of the sums in (8.7) and (8.8), so we find
that
$$
c=(a;q)^{-1/2}_\infty,\ \ \ \  c'=(1/a;q)^{-1/2}_\infty . \eqno
(8.13)
$$
The coefficients $c_n$ and $c'_n$ are thus real and equal to
$$
c_n=\frac{q^{n/2}}{(a;q)_\infty^{1/2}(q/a;q)_n^{1/2}(q;q)_n^{1/2}},
\ \ \ \
c'_n=\frac{q^{n/2}}{(1/a;q)_\infty^{1/2}(aq;q)_n^{1/2}(q;q)_n^{1/2}}.
$$
The orthogonality of the matrix $M$ means that
$$
\sum _m {\hat a}_{mn}{\hat a}_{mn'}=\delta_{nn'},\ \ \ \sum _m
{\hat a}'_{mn}{\hat a}'_{mn'}=\delta_{nn'},\ \ \ \sum _m {\hat
a}_{mn}{\hat a}'_{mn'}=0, \eqno (8.14)
$$  $$
\sum _n ({\hat a}_{mn}{\hat a}_{m'n}+{\hat a}'_{mn}{\hat
a}'_{m'n})=\delta_{mm'} .   \eqno (8.15)
$$
Substituting into (8.15) the expressions for ${\hat a}_{mn}$ and
${\hat a}'_{mn}$, we obtain
$$
\frac{1}{(a;q)_\infty}\sum_{m=0}^\infty
\frac{q^{m}}{(q/a;q)_m(q;q)_m}\,U_n^{(a)}(q^m;q)\,U_{n'}^{(a)}(q^{m};q)
\qquad\qquad\qquad\qquad\qquad
$$  $$
+\frac{1}{(1/a;q)_\infty}\sum_{m=0}^\infty
\frac{q^{m}}{(aq;q)_m(q;q)_m}\,U_n^{(a)}(aq^m;q)\,U_{n'}^{(a)}(aq^{m};q)
$$   $$  \qquad\qquad\qquad\qquad\qquad\qquad
= (-a)^n(q;q)_nq^{n(n-1)/2}\delta_{nn'} , \eqno (8.16)
$$
which must yield the orthogonality relation for Al-Salam--Carlitz
I polynomials. A problem, which remains to be clarified, is the
following. We have assumed that the points $q^n$ and $aq^n$,
$n=0,1,2,\cdots$, exhaust the whole spectrum of $B_2$. Let us show
that this is the case.

If the operator $B_2$ had other spectral points $x_k$, then on the
left-hand side of (8.16) there would be other summands
$\mu_{x_k}\, U_n^{(a)}(x_k;q)\,U_{n'}^{(a)}(x_k;q)$, corresponding
to these additional points. Let us show that these additional
summands do not appear. For this we set $n=n'=0$ in the relation
(8.16) with the additional summands. This results in the equality
$$
\frac{1}{(a;q)_\infty}\sum_{m=0}^\infty
\frac{q^{m}}{(q/a;q)_m(q;q)_m}+
\frac{1}{(1/a;q)_\infty}\sum_{m=0}^\infty
\frac{q^{m}}{(aq;q)_m(q;q)_m}+\sum _k \mu_{x_k}=1 . \eqno (8.17)
$$
In order to show that $\sum _k \mu_{x_k}=0$, take into account the
relation
$$
\frac{(Aq/C,Bq/C;q)_\infty}{(q/C,ABq/C;q)_\infty}
{}_2\phi_1(A,B;C;q,q)+ \frac{(A,B;q)_\infty}{(C/q,ABq/C;q)_\infty}
{}_2\phi_1(Aq/C,Bq/C;q^2/C;q,q)=1
$$
(see formula (2.10.13) in [GR]). Putting here $A=0$, $B=0$ and
$C=q/a$, we obtain relation (8.17) without the summand $\sum _k
\mu_{x_k}$. Therefore, one has $\sum_k \mu_{x_k} =0$. This means
that additional summands do not appear in (8.16) and it does
represent the orthogonality relation for the Al-Salam--Carlitz
polynomials. Due to this orthogonality, we arrive at the following
statement:
\medskip

{\bf Proposition 8.1.} {\it The spectrum of the operator $B_2$
coincides with the set of points $q^{n}$ and $aq^n$,
$n=0,1,2,\cdots$. This spectrum is simple and has one accumulation
point at 0.}
\bigskip

\noindent {\bf 8.3. Duals to Al-Salam--Carlitz I polynomials}
\bigskip

Now we consider the identities (8.14), which give the
orthogonality relations for the matrix elements ${\hat a}_{mn}$
and ${\hat a}'_{mn}$, considered as functions of $m$. Up to
multiplicative factors they coincide with
$$
F_n(x;a;q)={}_2\phi_1 (x,q^{-n};\; 0;\; q,q^{n+1}/a),\ \ \
n=0,1,2,\cdots , \eqno (8.18)
$$  $$
F'_n(x;a;q)={}_2\phi_1 (x,aq^{-n};\; 0;\; q,q^{n+1}),
n=0,1,2,\cdots , \eqno (8.19)
$$
considered on the set of points $q^{-m}$, $m=0,1,2,\cdots$.
Namely, we have
$$
{\hat a}_{mn}=c\left( \frac{q^{n}}{(q/a;q)_n (q;q)_n )}
\frac{q^{-m(m-1)/2}} {(q;q)_m\, (-a)^m}\right)^{1/2}
(-a)^mq^{m(m-1)/2}F_n(q^{-m};a;q),
$$   $$
{\hat a}'_{mn}=c'\left( \frac{q^{n}}{(aq;q)_n (q;q)_n )}
\frac{q^{-m(m-1)/2}} {(q;q)_m\, (-a)^m}\right)^{1/2}
(-a)^mq^{m(m-1)/2}F'_n(q^{-m};a;q),
$$
where $c$ and $c'$ are given by (8.13). The relations (8.14) lead
to the following orthogonality relations for the functions (8.18)
and (8.19):
$$
(a;q)^{-1}_\infty \sum_{m=0}^\infty \rho(m)
F_n(q^{-m};a;q)F_{n'}(q^{-m};a;q) =(q/a;q)_n(q;q)_nq^{-n}
\delta_{nn'}, \eqno (8.20)
$$  $$
(1/q;q)_\infty \sum_{m=0}^\infty \rho(m)
F'_n(q^{-m};a;q)F'_{n'}(q^{-m};a;q) =(aq;q)_n(q;q)_nq^{-n}
\delta_{nn'}, \eqno (8.21)
$$  $$
 \sum_{m=0}^\infty \rho(m) F_n(q^{-m};a;q)F'_{n'}(q^{-m};a;q)
=0, \eqno (8.22)
$$
where
$$
\rho(m)=\frac{(-a)^mq^{m(m-1)/2}}{(q;q)_m}.
$$

Comparing the expression (8.18) for the functions
$F_n(q^{-m};a;q)$ with the expression
$$
C_n(q^{-m};a';q):={}_2\phi_1 (q^{-n},q^{-m};\; 0;\; q,-q^{n+1}/a')
$$
for the $q$-Charlier polynomials, we conclude that
$$
F_n(q^{-m};a;q)=C_n(q^{-m};-a;q).  \eqno (8.23)
$$
Applying to the expression for the functions $F'_n(q^{-m};a;q)$
the transformation formula
$$
{}_2\phi_1 (q^{-n},b;\; 0;\; q,z)=(bz/q)^n {}_2\phi_1
(q^{-n},q/z;\; 0;\; q,q/b)
$$
(see (III.6) from Appendix III in [GR]), we derive that
$$
F'_n(q^{-m};a;q)=(-a)^{-m}C_n(q^{-m};-1/a;q).  \eqno (8.24)
$$
Substituting the expressions (8.23) and (8.24) into the relations
(8.20) and (8.21), we obtain the orthogonality relations for the
$q$-Charlier polynomials $C_n(q^{-m};-a;q)$ and
$C_n(q^{-m};-1/a;q)$, where $a<0$. For $C_n(q^{-m})\equiv
C_n(q^{-m};a';q)$, $a'>0$, it has the form
$$
\sum_{m=0}^\infty
\frac{{a'}^mq^{m(m-1)/2}}{(q;q)_m}C_n(q^{-m})C_{n'}(q^{-m})=
(-a';q)_\infty\, q^{-n}(-q/a';q)_n(q;q)_n\delta_{nn'}.
$$
It coincides with the orthogonality relation known from the
literature (see, for example, Chapter 7 in [GR]).

Thus, we have shown that duals of the family of Al-Salam--Carlitz
I polynomials $U_n^{(a)}(q^{-m};q)$ are two sets of $q$-Charlier
polynomials, one taken with the parameter $-a$ and the second one
with the parameter $-1/a$.

The relation (8.22) leads to the following equality for
$q$-Charlier polynomials:
$$
\sum_{m=0}^\infty
\frac{q^{m(m-1)/2}}{(q;q)_m}C_n(q^{-m};-a;q)C_{n'}(q^{-m};-1/a;q))=0,
$$
where $a<0$.

The set of functions (8.18) and (8.19) form an orthonormal basis
in the Hilbert space ${\frak l}^2$ of functions, defined on the
set of points $m=0,1,2,\cdots$, with the scalar product
$$
\langle f_1,f_2\rangle  =\sum_{m=0}^\infty \rho(m)
f_1(m)\overline{f_2(m)} ,
$$
where $\rho(m)$ is the same as in formulas (8.20)--(8.22). One can
deduce from this fact that {\it the $q$-Charlier polynomials
$C_n(q^{-m};a';q)$, $a'>0$, correspond to indeterminate moment
problem and the orthogonality measure for them, obtained above, is
not extremal.}
\bigskip

\noindent{\bf 9. DUALITY OF LITTLE $q$-LAGUERRE AND\\
Al-SALAM--CARLITZ II POLYNOMIALS}
 \bigskip

\noindent{\bf 9.1. Pair of operators $(B_3,Q^{-1})$}
 \bigskip

Let ${\cal H}\equiv {\cal H}_a$ be a separable complex Hilbert
space of functions, used in sections 3--7, with the polynomial
basis $f_n$, $n=0,1,2,\cdots$, in it.
We fix a real number $a$ such that $0<a<q^{-1}$. Let $B_3$ be the
operator on ${\cal H}\equiv {\cal H}_a$, acting upon the basis
elements $f_n$ as
$$
B_3 f_n =a_nf_{n+1}+a_{n-1}f_{n-1} +b_nf_n , \eqno (9.1)
$$
with
$$
a_n=-a^{1/2}q^{n+1/2}\sqrt{(1-q^{n+1})(1-aq^{n+1})},\ \ \ \
b_n=q^n(1+a)-aq^{2k}(1+q).
$$
Clearly, $B_3$ is a symmetric operator.

Since $a_n\to 0$ and $b_n\to 0$ when $n\to \infty$, the operator
$B_3$ is bounded. We assume that it is defined on the whole
Hilbert space ${\cal H}$ and, therefore, it is a self-adjoint
operator. Exactly as in the previous cases one can show that $B_3$
is a Hilbert--Schmidt operator. This means that the spectrum of
$B_3$ is discrete and has a single accumulation point at 0.
Moreover, a spectrum of $B_3$ is simple, since $B_3$ is
representable by a Jacobi matrix with $a_n\ne 0$.

To find eigenfunctions $\xi_\lambda $ of the operator $B_3$, $B_3
\xi_\lambda =\lambda \xi_\lambda $, we set
$$
\xi_\lambda =\sum _n \beta_n(\lambda)f_n,
$$
where $\beta_n(\lambda)$ are appropriate numerical coefficients.
Acting by the operator $B_3$ upon both sides of this relation, one
derives that
$$
\sum _{n=0}^{\infty}\, \beta_n(\lambda)\, (a_nf_{n+1}
+a_{n-1}f_{n-1} +b_nf_n )= \lambda \sum_{n=0}^{\infty}\,
\beta_n(\lambda) f_n ,
$$
where $a_n$ and $b_n$ are the same as in (9.1). Collecting in this
identity all factors, which multiply $f_n$ with fixed $n$, one
derives the recurrence relation for the coefficients
$\beta_n(\lambda)$:
$$
\beta_{n+1}(\lambda)a_n +\beta_{n-1}(\lambda)a_{n-1}+
\beta_{n}(\lambda)b_n= \lambda \beta_{n}(\lambda).
$$
The substitution
$$
\beta_{n}(\lambda)=\left( \frac{(aq;q)_n} {(aq)^{n}(q;q)_n}
\right) ^{1/2} \beta'_{n}(\lambda)
$$
reduces this relation to the following one
$$
- q^{n}(1-aq^{n+1})\beta'_{n+1}(\lambda)- aq^{n}(1-q^{n})
\beta'_{n-1}(\lambda)  + (q^n-aq^{2n+1}+
aq^{n}-aq^{2n})\beta'_n(\lambda) = \lambda \beta'_n(\lambda).
$$
This is the recurrence relation for the little $q$-Laguerre (Wall)
polynomials
$$
p_n(\lambda;a|q):={}_2 \phi_1 (q^{-n},0;\ aq;\ q;q\lambda)
=(a^{-1}q^{-n};q)^{-1}_n\; {}_2 \phi_0 (q^{-n},\lambda^{-1};\ -\,
;\ q;\lambda/a).  \eqno (9.2)
$$
Thus, we have $\beta'_k(\lambda) =p_n(\lambda;a|q)$  and,
consequently,
$$
\beta_n(\lambda)=\left(\frac{(aq;q)_n}{(aq)^n(q;q)_n}\right)^{1/2}
\,p_n(\lambda;a|q). \eqno (9.3)
$$
This means that eigenfunctions of the operator $B_3$ are of the
form
$$
\xi_\lambda (x) = \sum _{k=0}^\infty \left(
\frac{(aq;q)_k}{(aq)^k(q;q)_k}\right) ^{1/2}p_k(\lambda;a|q)\, f_k
 =\sum _{k=0}^\infty a^{-k/4}\,
\frac{(aq;q)_k}{(q;q)_k}\,p_k(\lambda;a|q)\, x^k .  \eqno (9.4)
$$
The expression for the eigenfunctions can be summed up. To show
this one needs to know a generating function
$$
F(x;\, t;\, a|q):= \sum_{n=0}^\infty \frac{(aq;q)_n}{(q;q)_n}\,
p_n(x;\, a|q)\,t^n  \eqno (9.5)
$$
for the little $q$-Laguerre polynomials. To evaluate (9.5), we
start with the second expression in (9.2) in terms of the basic
hypergeometric series ${}_2 \phi_0 $. Substituting it into (9.5)
and using the relation
$$
\frac{(q^{-n};q)_k}{(q;q)_k}=(-1)^kq^{-kn+k(k-1)/2}
\frac{(q;q)_n}{(q;q)_k(q;q)_{n-k}},
$$
one obtains that
$$
F(x;\, t;\, a|q)=\sum_{n=0}^\infty (-aqt)^nq^{n(n-1)/2} \sum
_{k=0}^n \frac{(x^{-1};q)_k}{(q;q)_k(q;q)_{n-k}} (q^{-n}x/a)^k .
\eqno (9.6)
$$
Interchanging the order of summations in (9.6) leads to the
desired expression
$$
F(x;\, t;\, a|q)=E_q(-aqt)\, {}_2 \phi_0 (x^{-1},0;\ - ;\ q;xt),
\eqno (9.7)
$$
where $E_q(z)=(-z;q)_\infty$ is the $q$-exponential function of
Jackson.

Similarly, if one substitutes into (9.5) the explicit form of the
little $q$-Laguerre polynomials in terms of ${}_2 \phi_1$ from
(9.2), this yields an expression
$$
F(x;\, t;\, a|q)=\frac{E_q(-aqt)}{E_q(-t)}\,{}_2 \phi_1 (0,0;\ q/t
;\ q;qx). \eqno (9.8)
$$

Using in (9.4) the explicit form of the generating function (9.7)
for the little $q$-Laguerre polynomials, we arrive at
$$
\xi _\lambda (x)= E_q(-qa^{-3/4}x) \, {}_2\phi_0 (\lambda^{-1},\,
0;\, - \, ;\ q;\, a^{-1/4}\lambda x) .
$$
Another expression for $\xi_\lambda (x)$ can be written by using
formula (9.8).

Since the spectrum of the operator $B_3$ is discrete, only for a
discrete set of values of $\lambda$ the functions (9.4) belong to
the Hilbert space ${\cal H}$. This discrete set of eigenvectors
determines a spectrum of $B_3$.

Now we look for a spectrum of $B_3$ and for a set of polynomials,
dual to little $q$-Laguerre polynomials. To this end we use the
operator $Q^{-1}$, where $Q$ is given by the formula
$Qf_n=q^nf_n$. In order to find how the operator $Q^{-1}$ acts
upon eigenfunctions of $B_3$, one can use the $q$-difference
equation
$$
q^{-n}\lambda p_n(\lambda)=- ap_{n}
(q\lambda)+(2+a-\lambda)p_n(\lambda)-
(1-\lambda)p_n(q^{-1}\lambda) \eqno(9.9)
$$
for the little $q$-Laguerre polynomials $p_n(\lambda)\equiv
p_n(\lambda ;a|q)$ (see formula (3.20.4) in [KSw]). Multiply both
sides of (9.9) by $d_n\,|n\rangle$ and sum up over $n$, where
$d_n$ are the coefficients of $p_n(\lambda ;a|q)$ in the
expression (9.3) for $\beta_n(\lambda)$. Taking into account
formula (9.9) and the fact that $Q^{-1}f_n=q^{-n}f_n$, one obtains
the relation
$$
Q^{-1}\,\xi _{\lambda}= -a\lambda^{-1}\,\xi _{q\lambda}+
\lambda^{-1}(2+a-\lambda)\, \xi _{\lambda}-
\lambda^{-1}(1-\lambda)\, \xi_{q^{-1}\lambda}, \eqno (9.10)
$$
which is used in the next section.
 \bigskip

\noindent{\bf 9.2. The spectrum of $B_3$ and orthogonality of
little $q$-Laguerre polynomials}
 \bigskip

Let us find, by using the operators $B_3$ and $Q^{-1}$, a basis in
the Hilbert space ${\cal H}$, which consists of eigenfunctions of
the operator $B_3$ in a normalized form, and the unitary matrix
$A$, connecting this basis with the initial basis $f_n$,
$n=0,1,2,\cdots$, in ${\cal H}$. First we have to find a spectrum
of $B_3$.

Let us first look at a form of the spectrum of $B_3$. If $\lambda$
is a spectral point of the operator $B_3$, then (as it is easy to
see from (9.10)) a successive action by the operator $Q^{-1}$ upon
the function (eigenfunction of $B_3$) $\xi_\lambda$ leads to the
eigenfunctions $\xi_{q^m\lambda}$, $m=0,\pm 1, \pm 2,\cdots$.
However, since $B_3$ is a Hilbert--Schmidt operator, not all of
these points belong to the spectrum of $B_3$, since $q^{-m}\lambda
\to\infty$ when $m\to +\infty$. This means that the coefficient
${\lambda'}^{-1}-1 $ of $\xi _{q^{-1}\lambda'}$ in (9.10) must
vanish for some eigenvalue $\lambda'$. Clearly, it vanishes when
$\lambda' =1$. Moreover, this is the only possibility for the
coefficient of $\xi _{q^{-1}\lambda'}$ in (9.10) to vanish, that
is, the point $\lambda =1$ is a spectral point for the operator
$B_3$. Let us show that the corresponding eigenfunction $\xi
_{1}\equiv \xi_{q^{0}}$ belongs to the Hilbert space ${\cal H}$.

One has the following equality
$$
p_n(1 ;a|q)={}_2\phi_1 (q^{-n}, 0;\; aq; \; q,q
)=\frac{(aq)^nq^{n(n-1)/2}}{(aq;q)_n}.
$$
Therefore, for the scalar product $\langle
\xi_1,\xi_1\rangle$ in ${\cal H}$ we have
$$
\langle \xi_1,\xi_1\rangle = \sum_{n=0}^\infty
\frac{(aq;q)_n}{(aq)^n(q;q)_n}p^2_n(1 ;a|q) = (-1;q)_\infty.
$$
Thus, the point $\lambda =1$ does belong to the spectrum of $B_3$.

Let us find other spectral points of the operator $B_3$. Setting
$\lambda = 1$ in (9.10), we see that the operator $Q^{-1}$
transforms $\xi _{q^0}$ into a linear combination of the vectors
$\xi _{q}$ and $\xi_{q^0}$. Moreover, $\xi_q$ belongs to the
Hilbert space ${\cal H}$, since the series
$$
\langle \xi_q,\xi_q\rangle = \sum_{n=0}^\infty
\frac{(aq;q)_n}{(aq)^n(q;q)_n}\,p^2_n(q ;a|q)
$$
is majorized by the corresponding series for $\xi_{q^0}$.
Therefore, $\xi _{q}$ belongs to the Hilbert space ${\cal H}$ and
the point $q$ is an eigenvalue of the operator $B_3$. Similarly,
setting $\lambda=q$ in (9.10), one finds likewise that $\xi
_{q^2}$ is an eigenvector of $B_3$ and the point $q^2$ belongs to
the spectrum of $B_3$. Repeating this procedure, we find that all
$\xi _{q^n}$, $n=0,1,2,\cdots$, are eigenvectors of $B_3$ and the
set $q^n$, $n=0,1,2,\cdots$, belongs to the spectrum of $B_3$. So
far, we do not know yet whether other spectral points exist or
not.

The functions $\xi _{q^n}$, $n=0,1,2,\cdots$, are linearly
independent elements of the Hilbert space ${\cal H}$. Suppose that
values $q^n$, $n=0,1,2,\cdots$, constitute a whole spectrum of
$B_3$. Then the set of functions $\xi _{q^n}$, $n=0,1,2,\cdots$,
is a basis in the Hilbert space ${\cal H}$. Introducing the
notation $\Xi _k:=\xi_{q^k}$, $k=0,1,2,\cdots$, we find from
(9.10) that
$$
Q^{-1} \,\Xi _k = - aq^{-k} \,\Xi _{k+1} + q^{-k}(2+a-q^k)\, \Xi
_k - q^{-k}(1-q^k)\, \Xi _{k-1} .   \eqno (9.12)
$$
As we see, the matrix of the operator $Q^{-1}$ in the basis $\Xi
_k$, $k=0,1,2,\cdots$, is not symmetric, although in the initial
basis $|n\rangle$, $n=0,1,2,\cdots$, it was symmetric. The reason
is that the matrix $(a_{mn})$ with entries $a_{mn}:=\beta_m(q^n)$,
$m,n=0,1,2,\cdots$, where $\beta_m(q^n)$ are the coefficients
(9.3) in the expansion $\xi _{q^n}=\sum _m \,\beta_m(q^n)f_n$, is
not unitary. This fact is equivalent to the statement that the
basis $\Xi _n=\xi_{q^n}$, $n=0,1,2,\cdots$, is not normalized. To
normalize it, one has to multiply $\Xi _n$ by corresponding
numbers $c_n$. Let $\hat\Xi _n = c_n\Xi _n$, $n=0,1,2,\cdots$, be
a normalized basis. Then the matrix of the operator $Q^{-1}$ is
symmetric in this basis. It follows from (9.12) that $Q^{-1}$ has
in the basis $\{ \hat\Xi _n\}$ the form
$$
Q^{-1}\, \hat\Xi _n = -c_{n+1}^{-1}c_nq^{-n} a\, \hat\Xi _{n+1} +
q^{-n}(2+a-q^n)\, \hat\Xi _n - c_{n-1}^{-1}c_n q^{-n}(1-q^n)
\,\hat\Xi_{n-1} .
$$
The symmetricity of $Q^{-1}$ in the basis $\{ \hat\Xi _n\}$ means
that $c_{n+1}^{-1}c_naq^{-n}=c_{n}^{-1}c_{n+1}
q^{-n-1}(1-q^{n+1})$, that is, $c_{n}/c_{n-1} =\sqrt{aq/(1-q^n)}$.
Therefore,
$$
c_n= c((aq)^n/(q;q)_n )^{1/2},
$$
where $c$ is a constant.

The expansions
$$
 \hat\xi _{q^n}(x)\equiv \hat\Xi _n(x)= \sum _m
c_n\beta_m(q^n)f_m \equiv \sum _m {\hat a}_{mn}f_m
$$
connect two orthonormal bases in the Hilbert space ${\cal H}$.
This means that the matrix $({\hat a}_{mn})$, $m,n=0,1,2,\cdots$,
with entries
$$
{\hat a}_{mn}=c_n\beta _m(q^n)= c\left( \frac{(aq)^n}{(q;q)_n}
\frac{(aq;q)_m}{(aq)^m(q;q)_m} \right)^{1/2} p_m(q^n ;a|q) ,
\eqno(9.13)
$$
is unitary, provided  that the constant $c$ is appropriately
chosen. In order to calculate this constant, we use the relation
$\sum_{n=0}^\infty |{\hat a}_{0n}|^2=\sum_{n=0}^\infty
c_n^2\beta_0^2(q^n)=1$. Since $\beta_0^2(q^n)=1$ and
$$
\sum_{n=0}^\infty \frac{(aq)^n}{(q;q)_n}=(aq;q)^{-1}_\infty ,
$$
we have
$$
c=(aq;q)^{1/2}_\infty .
$$

The matrix $A:=({\hat a}_{mn})$ is real and orthogonal. Thus, if
$\hat\Xi _n$, $n=0,1,2,\cdots$, is a complete basis in ${\cal H}$,
then $AA^{-1}=A^{-1}A=E$, that is,
$$
\sum _n {\hat a}_{mn}{\hat a}_{m'n}=\delta_{mm'},\ \ \ \ \sum _m
{\hat a}_{mn}{\hat a}_{mn'}=\delta_{nn'} .   \eqno (9.14)
$$
Substituting into the first sum over $n$ the expressions for
${\hat a}_{mn}$, we obtain the identity
$$
\sum_{n=0}^\infty
\frac{(aq)^n}{(q;q)_n}\,p_m(q^n;a|q)\,p_{m'}(q^n;aq) =
\frac{(aq)^m (q;q)_m } {(aq;q)_\infty (aq;q)_m }\, \delta_{mm'}\,,
\eqno (9.15)
$$
which must yield the orthogonality relation for little
$q$-Laguerre polynomials. However, we have assumed that the points
$q^n$, $n=0,1,2,\cdots$, exhaust the whole spectrum of $B_3$. Let
us show that this is the case. The reasoning here is the same as
in the previous sections. Namely, we have found that the spectrum
of $B_3$ contains the points $q^n$, $n=0,1,2,\cdots$. If the
operator $B_3$ had other spectral points $x_k$, then on the
left-hand side of (9.15) there would be other summands
$\mu_{x_k}\, p_m({x_k};a|q)\,p_{m'}({x_k};a|q)$ with positive
$\mu_{x_k}$, corresponding to these additional points. Let us show
that these additional summands do not appear. We set $m=m'=0$ in
the relation (9.15) with the additional summands. Since
$p_0(x;a|q)=1$, we have the equality
$$
\sum_{n=0}^\infty \frac{(aq)^n}{(q;q)_n} + \sum_k \mu_{x_k}
=(aq;q)^{-1}_\infty  .
$$
This formula is true if $\sum_k \mu_{x_k} =0$. This means that
additional summands do not appear in (9.15) and thus (9.15) does
represent the orthogonality relation for little $q$-Laguerre
polynomials. Consequently, the following proposition is true:
\medskip

 {\bf Proposition 9.1.} {\it The spectrum of the operator $B_3$
coincides with the set of points $q^{n}$, $n=0,1,2,\cdots$. This
spectrum is simple and the functions $\xi_{q^n}$,
$n=0,1,2,\cdots$, form a complete set of eigenfunctions of $B_3$.
The matrix $({\hat a}_{mn})$ with entries (9.13) relates the
initial basis $\{ f_n\}$ with the normalized basis $\{
\hat\Xi_n\}$. }
 \bigskip

\noindent{\bf 9.3. Al-Salam--Carlitz II polynomials as duals to
little $q$-Laguerre polynomials}
 \bigskip

Now we consider the second relation in (9.14). Taking into account
the explicit expression for $\hat a_{mn}$, one obtains the
orthogonality relation for the functions
$$
F_n(q^{-m};\, a| q):= (a^{-1}q^{-m};q)^{-1}_m\,
{}_2\phi_0(q^{-m},\, q^{-n};\, -\, ;q,q^n/a). \eqno (9.16)
$$
This relation has the form
$$
\sum_{m=0}^\infty  (aq)^{-m} \frac{(aq;q)_m}{(q;q)_m}
F_n(q^{-m};\, a| q) F_{n'}(q^{-m};\, a| q) =(aq)^{-n}
\frac{(q;q)_n}{(aq;q)_\infty} \delta _{nn'}.  \eqno (9.17)
$$
Comparing (9.16) with the Al-Salam--Carlitz II polynomials
$$
V^{(a)}_n(x;q)=(-a)^n q^{-n(n-1)/2} {}_2\phi_0 (q^{-n},\; x;\; -
\; ;\; q, q^n/a).
$$
we see that they are related to the functions (9.16), and (9.17)
therefore leads to the orthogonality relation for the
Al-Salam--Carlitz II polynomials
$$
\sum_{m=0}^\infty   \frac{q^{m^2}a^m}{(q;q)_m(aq;q)_m}
V^{(a)}_n(q^{-m}; q) V^{(a)}_{n'}(q^{-m}; q)  =
\frac{a^n(q;q)_n}{(aq;q)_\infty q^{n^2}} \delta _{nn'},
$$
known from the literature. Thus, Al-Salam--Carlitz II polynomials
are duals to the little $q$-Laguerre polynomials.
 \bigskip

\noindent {\bf Appendix}
\bigskip

In this appendix we prove the summation formula
$$
\sum _{n=0}^\infty \frac{(abq,bq;q)_n}{(aq,q;q)_n}
\frac{1-abq^{2n+1}}{1-abq}
a^nq^{n^2}=\frac{(abq^2;q)_\infty}{(aq;q)_\infty}. \eqno (A.1)
$$
First of all, observe that when $b=0$ this relation reduces to
$$
\sum_ {n=0}^\infty \frac{a^nq^{n^2}}{(aq,q;q)_n}=
\frac{1}{(aq;q)_\infty} ,
$$
which is a well-known limiting form of Jacobi's triple product
identity (see [4], formula (1.6.3)).

One can employ an easily verified relation
$$
\frac{(aq,-aq;q)_n}{(a,-a;q)_n}=\frac{1-a^2q^{2n}}{1-a^2}
 \eqno (A.2)
$$
in order to express the infinite sum in (A.1) in terms of a
very-well-poised ${}_4\phi_5$ basic hypergeometric series. This
results in
$$
\sum _{n=0}^\infty \frac{(abq,bq;q)_n}{(aq,q;q)_n}
\frac{1-abq^{2n+1}}{1-abq} a^nq^{n^2} ={}_4\phi_5 \left( \left.
{abq,\; bq,\; q\sqrt{abq},\; -q\sqrt{abq} \atop  aq,\;
\sqrt{abq},\; -\sqrt{abq},\; 0,\; 0} \right|
 q,aq \right) . \eqno (A.3)
$$
The next step is to utilize a limiting case of Jackson's sum of a
terminating very-well-poised balanced ${}_8\phi_7$ series,
$$
{}_6\phi_5 \left( \left. {a,\;  q\sqrt{a},\; -q\sqrt{a},\; b,\;
c,\; d \atop
             \sqrt{a},\;  -\sqrt{a},\; aq/b,\; aq/c,\; aq/d} \right|
 q,\frac{aq}{bcd} \right)
 = \frac{(aq,aq/bc,aq/bd,aq/cd;q)_\infty}
{(aq/b,aq/c,aq/d,aq/bcd;q)_\infty} , \eqno (A.4)
$$
which represents a $q$-analogue of Dougall's formula for a
very-well-poised 2-balanced ${}_7F_6$ series. When the parameters
$c$ and $d$ tend to infinity, from (A.4) it follows that
$$
{}_4\phi_5 \left( \left. {a,\; q\sqrt{a},\; -q\sqrt{a},\; b \atop
    \sqrt{a},\;  -\sqrt{a},\; aq/b,\;  0,\; 0} \right|
 q,\frac{aq}{b} \right)  = \frac{(aq;q)_\infty}{(aq/b;q)_\infty} .
\eqno (A.5)
$$
To verify this, one needs only to use the limit relation
$$
\lim_{c,d\to \infty} (c,d;q)_n\left( \frac{aq}{bcd}\right)^n =
q^{n(n-1)} \left( \frac{aq}{b}\right) ^n .
$$
With the substitutions $a\to abq$ and $b\to bq$ in (A.5), one
recovers the desired identity (A.1).

Similarly, when $d\to \infty$ we derive from (A.4) the identities
$$
\sum_{n=0}^\infty\,\frac{(1-abq^{2n+1})(aq,abq/c,abq;q)_n}
{(1-abq)(bq,cq,q;q)_n (-a/c)^n} q^{n(n-1)/2}=
\frac{(abq^2,c/a;q)_\infty}{(bq,cq;q)_\infty}\, ,    \eqno (A.6)
$$  $$
\sum_{n=0}^\infty\,\frac{(1-abq^{2n+1})(abq,bq,cq;q)_n}
{(1-abq)(aq,abq/c,q;q)_n (-c/a)^n} q^{n(n-1)/2}
=\frac{(abq^2,a/c;q)_\infty}{(aq,abq/c;q)_\infty}\, .  \eqno (A.7)
$$
They have been employed in section 7.

We conclude this appendix with the following remark. There is
another proof of the identity (4.33), based on vital use of the
same summation formula (A.4). Actually, a relation may be derived,
which is somewhat more general than (4.33). Indeed, consider the
function
$$
\eta_k(a;q):=\sum_{n=0}^\infty (-1)^nq^{n(n-1)/2}
\frac{1-aq^{2n+1}}{1-aq} \frac{(aq;q)_n}{(q;q)_n} \mu^k(n;a)
 \eqno (A.8)
$$
for arbitrary nonnegative integers $k$, where the $q$-quadratic
lattice $\mu(n; a)$ is defined as above:
$$
\mu(n;a):= q^{-n}+aq^{n+1}.  \eqno (A.9)
$$
We argue that all $\eta_k(a;q)=0$, $k=0,1,2,\cdots$. To verify
that, begin with the case when $k=0$ and employ relation (A.2) to
show that
$$
\eta_0(a;q)=\left. {}_3\phi_3 \left( {q\sqrt{aq},\; -q\sqrt{aq},\;
aq     \atop \sqrt{aq},\;  -\sqrt{aq},\; 0} \right| q, 1 \right) .
$$
The summation formula (A.4) in the limit as $d\to \infty$ takes
the form
$$
\left. {}_5\phi_5 \left( {a,\; q\sqrt{a},\; -q\sqrt{a},\; b,\; c
  \atop \sqrt{a},\;  -\sqrt{a},\; aq/b,\; aq/c,\; 0} \right| q,\;
  \frac{aq}{bc}
  \right) =\frac{(aq,aq/bc;q)_\infty}{(aq/b,aq/c;q)_\infty}\ .
 \eqno (A.10)
$$
In the particular case when $bc=aq$ this sum reduces to
$$
\left. {}_3\phi_3 \left( {a,\; q\sqrt{a},\; -q\sqrt{a}
  \atop 0,\; \sqrt{a},\;  -\sqrt{a}} \right| q, 1
  \right) =\frac{(aq,1;q)_\infty}{(b,c;q)_\infty}=0 \,,
$$
since $(z;q)_\infty =0$ for $z=1$. Consequently, the function
$\eta_0(a;q)$ does vanish.

For $k=1,2,3,\cdots$, one can proceed inductively. Employ the
relation $q\mu (n+1;a)=\mu (n;q^2a)$ to show that
$$
\eta_{k+1}(a;q)=(1+aq)\eta_k(a;q)-q^{-k-1}(1-aq^2)(1-aq^3)
\eta_k(aq^2;q).
$$
So, one obtains that indeed all $\eta_k(a;q)$, $k=0,1,2,\cdots$,
vanish. The identity (4.33) is now an easy consequence of this
statement if one takes into account that a product of the two
polynomials $D_n(\mu(m); a,b,c|q)$ and $D_{n'}(\mu(m);
b,a,abq/c|q)$ in (4.33) is some polynomial in $\mu(m)$ of degree
$n+n'$. This completes the proof of (4.33), which is independent
of the one, given in section 4.

\end{document}